\documentclass[a4paper, 11pt]{amsart}

\usepackage[T1]{fontenc}
\usepackage[latin1]{inputenc}
\usepackage{amssymb}
\usepackage{euscript}
\usepackage{amsxtra}
\usepackage[all]{xy}
\usepackage{enumerate}
\usepackage{mathrsfs}

\usepackage{graphicx}
\usepackage[lite]{amsrefs}

\newcommand*{\brd}{-\hspace{0pt}}
\newcommand*{\nbd}{\nobreakdash-\hspace{0pt}}

\newcommand*{\blank}{\text{\textvisiblespace}}

\newcommand*{\Mult}{{\mathcal M}}     
\newcommand*{\Cred}{C^*_{\mathrm r}}  
\newcommand*{\Comp}{{\mathbb K}}      
\newcommand*{\Bound}{{\mathbb B}}     

\newcommand*{\abs}[1]{\lvert#1\rvert}   
\newcommand*{\norm}[1]{\lVert#1\rVert}  
\newcommand*{\gen}[1]{\langle#1\rangle} 

\newcommand*{\rcross}{\mathbin{\rtimes_{\mathrm{r}}}} 
\newcommand*{\cross}{\mathbin{\rtimes}}               
\newcommand*{\defeq}{\mathrel{:=}}                    
\newcommand*{\congto}{\overset{\cong}\to}             
\newcommand*{\hot}{\mathbin{\hat\otimes}}

\DeclareMathOperator{\Ad}{Ad}      
\DeclareMathOperator{\Isom}{Isom}  
\DeclareMathOperator{\CAT}{CAT}    
\DeclareMathOperator{\PSL}{PSL}    
\DeclareMathOperator{\Rep}{Rep}    
\DeclareMathOperator{\supp}{supp}  

\newcommand*{\N}{{\mathbb N}}      
\newcommand*{\Z}{{\mathbb Z}}      
\newcommand*{\R}{{\mathbb R}}      
\newcommand*{\C}{{\mathbb C}}      
\newcommand*{\F}{{\mathbb F}}      
\newcommand*{\Ztwo}{{\mathbb Z/2}} 
\newcommand*{\Hyp}{{\mathbb H}}    

\newcommand*{\ID}{{\mathrm{id}}}   
\newcommand*{\Hc}{\mathrm{H}_{\mathrm{c}}} 
\newcommand*{\pnt}{\mathrm{pnt}}   
\newcommand*{\alg}{{\mathrm{alg}}} 
\newcommand*{\HH}{\mathrm{HH}}     
\newcommand*{\KK}{\mathrm{KK}}     
\newcommand*{\RKK}{\mathrm{RKK}}   
\newcommand*{\K}{\mathrm{K}}       
\newcommand*{\Ktop}{\mathrm{K}^{\mathrm{top}}}
\newcommand*{\Eul}{\mathrm{Eul}}    
\newcommand*{\dR}{{\mathrm{dR}}}    
\newcommand*{\comb}{{\mathrm{cmb}}} 
\newcommand*{\op}{{\mathrm{op}}}    
\newcommand*{\Cliff}{{\mathbb{C}\mathrm{liff}}} 

\newcommand*{\ADir}{{\mathsf P}}   
\newcommand*{\Dirac}{{\mathsf D}}  
\newcommand*{\REP}{{\mathsf R}}    
\newcommand*{\Spinor}{\mathsf{Spinor}} 
\newcommand*{\Dslash}{{\mathsf /\!\!\!\!D}}  

\newcommand*{\EG}{{\mathcal EG}}   
\newcommand*{\Hilm}{{\mathcal E}}  
\newcommand*{\Hecke}{{\mathcal H}} 

\newcommand*{\Diagonal}{\Delta}    

\newcommand*{\cl}[1]{\bar #1}      
\newcommand*{\bd}[1]{\partial #1}  
\newcommand*{\vis}{\infty}         

\newcommand*{\sour}{{\mathscr P}}
\newcommand*{\comul}{\nabla}
\newcommand*{\PD}{\mathrm{PD}}

\newcommand*{\collapse}{{\mathscr C}}

\newcommand*{\poset}{\mathcal{S}(\no)}
\newcommand*{\no}{{\underline{\mathbf{n}}}}
\newcommand*{\nok}{{\underline{\mathbf{k}}}}

\newcommand*{\teich}{\mathscr{T}}          
\newcommand*{\map}{\mathrm{Mod}(\Sigma_g)} 

\theoremstyle{plain}
\newtheorem{theorem}{Theorem}
\newtheorem{lemma}[theorem]{Lemma}
\newtheorem{proposition}[theorem]{Proposition}
\newtheorem{corollary}[theorem]{Corollary}
\theoremstyle{definition}
\newtheorem{definition}[theorem]{Definition}
\theoremstyle{remark}
\newtheorem{remark}[theorem]{Remark}

\newtheorem{example}[theorem]{Example}

\begin{document}

\title[Euler characteristics and boundary actions]
{Euler characteristics and Gysin sequences for group actions on boundaries}

\author{Heath Emerson}
\email{hemerson@math.uni-muenster.de}

\author{Ralf Meyer}
\email{rameyer@math.uni-muenster.de}

\address{Mathematisches Institut der
         Westfälischen Wilhelms-Universität Münster\\
         Einsteinstraße 62\\
         48149 Münster\\
         Deutschland}

\begin{abstract}
  Let~\(G\) be a locally compact group, let~\(X\) be a universal proper
  \(G\)\nbd{}space, and let~\(\cl{X}\) be a \(G\)\nbd{}equivariant
  compactification of~\(X\) that is \(H\)\nbd{}equivariantly contractible for
  each compact subgroup \({H\subseteq G}\).  Let \(\bd{X}=\cl{X}\setminus X\).
  Assuming the Baum-Connes conjecture for~\(G\) with coefficients \(\C\) and
  \(C(\bd{X})\), we construct an exact sequence that computes the map on
  \(\K\)\nbd{}theory induced by the embedding \(\Cred G\to {C(\bd{X})\rcross
    G}\).  This exact sequence involves the equivariant Euler characteristic
  of~\(X\), which we study using an abstract notion of Poincaré duality in
  bivariant \(K\)-theory.  As a consequence, if~\(G\) is torsion-free and the
  Euler characteristic \(\chi(G\backslash X)\) is non-zero, then the unit
  element of \(C(\bd{X})\rcross G\) is a torsion element of order
  \(\abs{\chi(G\backslash X)}\).  Furthermore, we get a new proof of a theorem
  of Lück and Rosenberg concerning the class of the de Rham operator in
  equivariant \(\K\)\nbd{}homology.
\end{abstract}
 

\thanks{This research was supported by the EU-Network \emph{Quantum
  Spaces and Noncommutative Geometry} (Contract HPRN-CT-2002-00280)
  and the \emph{Deutsche Forschungsgemeinschaft} (SFB 478).}

\maketitle

\section{Introduction}
\label{sec:intro}

Let~\(G\) be a locally compact group, let~\(X\) be a proper \(G\)\nbd{}space,
and let~\(\cl{X}\) be a compact \(G\)-space containing~\(X\) as a
\(G\)\nbd{}invariant open subset.  Suppose that both \(X\) and~\(\cl{X}\) are
\(H\)\nbd{}equivariantly contractible for all compact subgroups~\(H\)
of~\(G\); we briefly say that they are \emph{strongly contractible} and call
the action of~\(G\) on \(\bd{X} \defeq \cl{X}\setminus X\) a \emph{boundary
  action}.

For example, the group \(G = \PSL(2,\Z)\) admits the following two distinct
boundary actions.  On the one hand, since~\(G\) is a free product of finite
cyclic groups, it acts properly on a tree~\(X\) (\cite{Serre:Trees}).  Let
\(\bd{X}\) be its set of ends and let~\(\cl{X}\) be its ends compactification,
which is a Cantor set.  Then \(X\) and~\(\cl{X}\) are strongly contractible
and the action of~\(G\) on~\(\bd{X}\) is a boundary action.  On the other
hand, \(\PSL(2,\Z) \subseteq\PSL(2,\R)\) acts by Möbius transformations on the
hyperbolic plane~\(\Hyp^2\).  We compactify~\(\Hyp^2\), as usual, by a circle
at infinity.  Again, \(\Hyp^2\) and~\(\cl{\Hyp^2}\) are strongly contractible
and the action on the circle \(\bd{\Hyp^2}\) is a boundary action.  Other
examples are: a word-hyperbolic group acting on its Gromov boundary; a group
of isometries of a \(\CAT(0)\) space~\(X\) acting on the visibility boundary
of~\(X\); a mapping class group of a Riemann surface acting on the Thurston
boundary of Teichmüller space; a discrete subgroup of \(\Isom(\Hyp^n)\) acting
on its limit set.  We discuss these examples in Section~\ref{sec:examples}.

The purpose of this article is to describe the map on \(\K\)-theory induced by
the obvious inclusion \(u\colon \Cred G \to C(\bd{X}) \rcross G\), where
\(G\times \bd{X} \to \bd{X}\) is a boundary action and \(C(\bd{X})\rcross G\)
is its reduced crossed product \(C^*\)-algebra.  Our result is analogous to
the classical Gysin sequence, which we recall first.

Let~\(M\) be a smooth \(n\)\nbd{}dimensional manifold.  Let \(TM\), \(BM\)
and~\(SM\) be its tangent, disk and sphere bundles, respectively.  Thus~\(TM\)
is an open subset of~\(BM\) with \(SM=BM\setminus TM\).  Let \(\Hc^*\) denote
cohomology with compact supports.  Since the bundle projection \(BM\to M\) is
a proper homotopy equivalence, we have \(\Hc^*(BM)\cong\Hc^*(M)\) and
\(\K^*(BM)\cong\K^*(M)\).  We assume now that~\(M\) is oriented or
\(\K\)\nbd{}oriented, respectively.  Then we get Thom isomorphisms
\(\Hc^{*-n}(M)\cong \Hc^*(TM)\) or \(\K^{*-n}(M)\cong \K^*(TM)\), and excision
for the pair \((BM,SM)\) provides us with long exact sequences of the form
\[\xymatrix@-2mm{
  \dotsb \ar[r] &
  \Hc^{*-n}(M) \ar[r]^-{\varepsilon^*} &
  \Hc^{*}(M) \ar[r]^-{\pi^*} &
  \Hc^{*}(SM) \ar[r]^-{\delta} &
  \Hc^{*-n+1}(M) \ar[r] &
  \dotsb,
}
\]
\[\xymatrix@-2mm{
  \dotsb \ar[r] &
  \K^{*-n}(M) \ar[r]^-{\varepsilon^*} &
  \K^*(M) \ar[r]^-{\pi^*} &
  \K^*(SM) \ar[r]^-{\delta} &
  \K^{*-n+1}(M) \ar[r] &
  \dotsb.
}
\]
These are the classical \emph{Gysin sequences}.

The map \(\varepsilon^*\colon \Hc^{*-n}(M)\to \Hc^*(M)\) vanishes unless
\(*=n\), for dimension reasons; if~\(M\) is not compact, then \(\Hc^0(M)=0\)
and hence \(\varepsilon^*=0\).  If~\(M\) is compact, then
\(\varepsilon^n\colon \Hc^0(M)\to\Hc^n(M)\) is given by \(x\mapsto x\wedge
e_M\), where \(e_M\in \Hc^n(M)\) is the Euler characteristic of~\(M\) (see
\cite{Bott-Tu}).

The map \(\varepsilon^*\colon \K^{*-n}(M)\to \K^*(M)\) is given by
\(\varepsilon^*(x)=x\otimes \Spinor\), where \(\Spinor\) denotes the spinor
bundle of~\(M\) (see \cite{Karoubi:K-theory}*{IV.1.13}).  Theorem
\ref{the:Gysin_classical} gives a more precise description: \(\varepsilon^*\)
vanishes on \(\K^1(M)\) and is given by \(x\mapsto \chi(M)\dim(x)\cdot \pnt!\)
on \(K^0(M)\); here \(\chi(M)\in\Z\) is the Euler characteristic of~\(M\), the
functional \(\dim\colon \K^0(M)\to\Z\) sends a vector bundle to its dimension,
and \(\pnt!\in\KK_{-n}(\C,C_0(M))\cong\K^n(M)\) is the wrong way element
corresponding to the inclusion of a point in~\(M\), which is a
\(\K\)\nbd{}oriented map.  Notice that \(\dim=0\) unless~\(M\) is compact.

We shall construct an analogue of the \(\K\)\nbd{}theoretic Gysin sequence for
boundary actions.  We have to assume that~\(X\), in addition to being strongly
contractible and proper, is a finite-dimensional simplicial complex equipped
with a simplicial action of~\(G\) or a smooth Riemannian manifold equipped
with an isometric action of~\(G\).  Although our arguments work for locally
compact groups, we concentrate on the case of torsion-free discrete groups in
the introduction, for expository purposes.  If~\(G\) is torsion-free,
then~\(X\) is a universal free proper \(G\)\nbd{}space, so that \(G\backslash
X\) is a model for the classifying space \(BG\).  We warn the reader that
\(\K^*(G\backslash X)\) depends on the particular choice of~\(BG\) because
\(\K\)\nbd{}theory is only functorial for \emph{proper} maps; we discuss this
for \(\PSL(2,\Z)\) at the end of Section~\ref{sec:examples} (see also Example
\ref{exa:PSL2Z}).

\begin{theorem}
  \label{the:bigtheorem}
  Let~\(G\) be a torsion-free discrete group and let \(G\times \bd{X} \to
  \bd{X}\) be a boundary action, where~\(X\) is a finite-dimensional
  simplicial complex with a simplicial action of~\(G\).  Assume that~\(G\)
  satisfies the Baum-Connes conjecture with coefficients \(\C\) and
  \(C(\bd{X})\).  Let \(u\colon \Cred G\to C(\bd{X})\rcross G\) be the
  embedding induced by the unit map \(\C\to C(\bd{X})\).

  If \(G\backslash X\) is compact and \(\chi(G\backslash X)\neq0\), then there
  are exact sequences
  \[\xymatrix@C-1.1em{
    0 \ar[r] &
    \gen{\chi(G\backslash X)[1_{\Cred G}]} \ar[r]^-{\subseteq} &
    \K_0(\Cred G) \ar[r]^-{u_*} &
    \K_0(C(\bd{X}) \rcross G) \ar[r] &
    \K^1(G\backslash X) \ar[r] & 0,
  }
  \]
  \[\xymatrix@C-.3em{
    0 \ar[r] &
    \K_1(\Cred G) \ar[r]^-{u_*} &
    \K_1(C(\bd{X}) \rcross G) \ar[r] &
    \K^0(G\backslash X) \ar[r]^-{\dim} &
    \Z \ar[r] & 0.
  }
  \]
  Here \(\gen{\chi(G\backslash X)[1_{\Cred G}]}\) denotes the free cyclic
  subgroup of \(\K_0(\Cred G)\) generated by \(\chi(G\backslash X)[1_{\Cred
    G}]\) and \(\dim\) maps a vector bundle to its dimension.

  If \(G\backslash X\) is not compact or if \(\chi(G\backslash X)=0\), then
  there are exact sequences
  \[\xymatrix@C-.3em{
    0 \ar[r] &
    \K_0(\Cred G) \ar[r]^-{u_*} &
    \K_0(C(\bd{X}) \rcross G) \ar[r] &
    \K^1(G\backslash X) \ar[r] & 0,
  }
  \]
  \[\xymatrix@C-.3em{
    0 \ar[r] &
    \K_1(\Cred G) \ar[r]^-{u_*} &
    \K_1(C(\bd{X}) \rcross G) \ar[r] &
    \K^0(G\backslash X) \ar[r] & 0.
  }
  \]
\end{theorem} 

\begin{corollary}
  \label{cor:torsion_unit}
  The class of the unit element in \(\K_0(C(\bd{X})\rcross G)\) is a torsion
  element of order \(\abs{\chi(G\backslash X)}\) if \(G\backslash X\) is
  compact and \(\chi(G\backslash X)\neq0\), and not a torsion element
  otherwise.
\end{corollary}

Several authors have already noticed various instances of this corollary
(\cites{Connes:Cyclic_transverse, Connes:NCG, Delaroche, Emerson,
  Natsume:Torsion, Spielberg:Free-product, Robertson:Torsion}): for lattices
in \(\PSL(2,\R)\) and \(\PSL(2,\C)\), acting on the boundary of hyperbolic
\(2\)- or \(3\)-space, respectively; for closed subgroups of \(\PSL(2,\F)\)
for a non-Archimedean local field~\(\F\) acting on the projective space
\(\mathbb{P}^1(\F)\), where~\(X\) is the Bruhat-Tits tree of \(\PSL(2,\F)\);
for free groups acting on their Gromov boundary.

Comparing the classical and non-commutative Gysin sequences, we see that the
inclusion \(u\colon \Cred G \to C(\bd{X})\rcross G\) plays the role of the
embedding \(C(M) \to C(SM)\) induced by the bundle projection \(SM \to M\).
Therefore, if we view \(\Cred G\) as the algebra of functions on a
non-commutative space~\(\hat{G}\), then \(C(\bd{X})\rcross G\) plays the role
of the algebra of functions on the sphere bundle of~\(\hat{G}\).  Such an
analogy has already been advanced by Alain Connes and Marc Rieffel in
\cites{Connes:Geometry_spectral, Rieffel:Group_metric} for rather different
reasons (and for a different class of boundary actions).

For groups with torsion and, more generally, for locally compact groups, the
Euler characteristic of \(G\backslash X\) is replaced by an \emph{equivariant}
Euler characteristic in \(\KK^G_0(C_0(X),\C)\).  To define it, we use a
general notion of Poincaré duality in bivariant Kasparov theory.  An
\emph{abstract dual} for a space~\(X\) consists of a
\(G\)-\(C^*\)-algebra~\(\sour\) and a class \(\Theta\in\RKK^G_n(X;\C,\sour)\)
for some \(n\in\Z\) such that the map
\[
\RKK^G_{*-n}(Y;A\hot\sour,B) \to \RKK^G_*(X\times Y;A,B),
\qquad f\mapsto \Theta\hot_\sour f,
\]
is an isomorphism for all pairs of \(G\)-\(C^*\)-algebras \(A\), \(B\) and all
\(G\)\nbd{}spaces~\(Y\) (compare \cite{Kasparov:Novikov}*{Theorem 4.9}).

Let~\(X\) be any \(G\)\nbd{}space that has such an abstract dual.  The
diagonal embedding \(X\to X\times X\) yields classes
\[
\Diagonal_X\in\RKK^G_0(X;C_0(X),\C),
\qquad
\PD^{-1}(\Diagonal_X)\in\KK^G_{-n}(C_0(X)\hot\sour,\C).
\]
Let \(\bar\Theta\in\KK^G_n(C_0(X),C_0(X)\hot\sour)\) be obtained from~\(\Theta\)
by forgetting the \(X\)\nbd{}linearity.  We define the \emph{abstract
  equivariant Euler characteristic} by
\[
\Eul_X\defeq \bar\Theta \hot_{C_0(X)\hot\sour} \PD^{-1}(\Diagonal_X)
\in\KK^G_0(C_0(X),\C).
\]
Examples show that this class deserves to be called an Euler characteristic.
We were led to this definition by the consideration of the Gysin sequence.

In order to compute \(\Eul_X\), we need an explicit formula for \(\PD^{-1}\).
Therefore, it is useful to consider a richer structure than an abstract dual,
which we call a \emph{Kasparov dual}.  Gennadi Kasparov constructs the
required structure for a smooth Riemannian manifold in
\cite{Kasparov:Novikov}*{Section 4}, using for~\(\sour\) the algebra of
\(C_0\)\nbd{}sections of the Clifford algebra bundle on~\(X\).  A fairly
simple computation shows that the associated equivariant Euler characteristic
is the class in \(\KK^G_0(C_0(X),\C)\) of the de Rham operator on~\(X\), which
we denote by \(\Eul_X^\dR\) and call the \emph{equivariant de-Rham-Euler
  characteristic of~\(X\)}.

If~\(X\) is a simplicial complex and~\(G\) acts simplicially, then a Kasparov
dual for~\(X\) is constructed in \cite{Kasparov-Skandalis:Buildings}.  Since
the description of~\(\Theta\) in \cite{Kasparov-Skandalis:Buildings} is too
indirect for our purposes, we give a slightly different construction where we
can write down~\(\Theta\) very concretely.  We describe the equivariant Euler
characteristic that we get from this combinatorial dual; the result may be
computed in terms of counting \(G\)\nbd{}orbits on the set of simplices and is
called the \emph{equivariant combinatorial Euler characteristic}
\(\Eul_X^\comb\in\KK^G_0(C_0(X),\C)\).

We show that the (abstract) Euler characteristic of~\(X\) does not depend on
the choice of the abstract dual.  Therefore, if~\(X\) admits a smooth
structure and a triangulation at the same time, then
\(\Eul_X^\dR=\Eul_X=\Eul_X^\comb\).  This result is due to Wolfgang Lück and
Jonathan Rosenberg (\cites{Luck-Rosenberg:Euler, Rosenberg:Euler}).  In the
non-equivariant case, the assertion is that for a connected smooth
manifold~\(M\), we have \(\Eul_M^\dR= \chi(M)\cdot\dim\), where
\(\dim\in\K_0(M)\) is the class of the point evaluation homomorphism.

Now we outline the proof of Theorem \ref{the:bigtheorem} and its analogues for
general locally compact groups.  The starting point is the extension of
\(G\)-\(C^*\)-algebras
\[
0 \to C_0(X) \to C(\cl{X}) \to C(\bd{X}) \to 0,
\]
which yields a six term exact sequence for the functor \(\Ktop_*(G,\blank)\).
The strong contractibility of~\(\cl{X}\) implies that
\(\Ktop_*(G,C(\cl{X}))\cong \Ktop_*(G)\).  The resulting map \(\Ktop_*(G)\cong
\Ktop_*(G,C(\cl{X})) \to\Ktop_*(G,C(\bd{X}))\) in the exact sequence is
induced by the unital inclusion \(\C\to C(\bd{X})\).  A purely formal argument
shows that the map \(\Ktop_*(G,C_0(X))\to
\Ktop_*(G,C(\cl{X}))\cong\Ktop_*(G)\) in the exact sequence is given by the
Kasparov product with the equivariant \emph{abstract} Euler characteristic
\(\Eul_X\in \KK^G_0(C_0(X),\C)\).  The heart of the computation is the explicit
description of \(\Eul_X\) as \(\Eul_X^\dR\) or \(\Eul_X^\comb\).

Our interest in the problem of calculating the \(\K\)-theory of boundary
crossed products was sparked by discussions with Guyan Robertson at a meeting
in Oberwolfach in 2004.  We would like to thank him for drawing our attention
to this question.  We also thank Wolfgang Lück for helpful suggestions
regarding Euler characteristics.

\subsection{General setup}
\label{sec:setup}

Throughout this article, topological spaces and groups are assumed locally
compact, Hausdorff, and second countable.  Let~\(G\) be such a group.  A
\emph{\(G\)\nbd{}space} is a locally compact space equipped with a continuous
action of~\(G\).  A \emph{\(G\)\nbd{}\(C^*\)\brd{}algebra} is a separable
\(C^*\)\brd{}algebra equipped with a strongly continuous action of~\(G\) by
automorphisms.  We denote \emph{reduced crossed products} by \(A\rcross G\).
We always equip~\(\C\) with the trivial action of~\(G\).  Thus \(\C\rcross G\)
is the \emph{reduced group \(C^*\)\brd{}algebra} \(\Cred G\) of~\(G\).

A \(G\)\nbd{}space is \emph{proper} if the set of \(g\in G\) with \(gK\cap
K\neq\emptyset\) is compact for any compact subset \(K\subseteq X\).  A
\emph{universal proper \(G\)\nbd{}space} is a proper \(G\)\nbd{}space \(\EG\)
with the property that for any proper \(G\)\nbd{}space~\(X\) there is a
continuous \(G\)\nbd{}equivariant map \(X\to\EG\) which is unique up to
\(G\)\nbd{}equivariant homotopy.  Such a \(G\)\nbd{}space exists for any~\(G\)
by \cite{Kasparov-Skandalis:Bolic}, and any two of them are
\(G\)\nbd{}equivariantly homotopy equivalent.

\begin{definition}
  \label{def:strong_contractible}
  We call a \(G\)\nbd{}space \emph{strongly contractible} if it is
  \(H\)\nbd{}equivariantly contractible for any compact subgroup \({H\subseteq
    G}\).
\end{definition}

A proper \(G\)\nbd{}space is strongly contractible if and only if it is
universal.  It is easy to see that universality implies strong
contractibility: look at maps \(G/H\times\EG\to\EG\).  The converse
implication is proved in full generality in \cite{Meyer-Nest:Ideals}.  The
case of a \(G\)-CW-complex is easier: in that case, it already suffices to
assume contractibility of fixed-point subsets
(\cite{Luck:Classifying_Survey}).

\begin{definition}
  \label{def:boundary_action}
  Let~\(X\) be a \emph{proper} \(G\)\nbd{}space and let~\(\cl{X}\) be a
  \emph{compact} \(G\)\nbd{}space that contains~\(X\) as an open
  \(G\)\nbd{}invariant subspace.  Then \(\bd{X}\defeq \cl{X}\setminus X\) is
  another compact \(G\)\nbd{}space.  We call the induced action of~\(G\)
  on~\(\bd{X}\) a \emph{boundary action} if both \(X\) and~\(\cl{X}\) are
  strongly contractible.  The \(G\)\nbd{}spaces \(X\) and~\(\cl{X}\) are part
  of the data of a boundary action.
\end{definition}

In all our examples, \(\cl{X}\) is a compactification of~\(X\), that is, \(X\)
is dense in~\(\cl{X}\).  A boundary action yields an extension of
\(G\)\nbd{}\(C^*\)\brd{}algebras
\begin{equation}
  \label{eq:G_space_ext}
  0 \to C_0(X) \overset{\iota}\to C(\cl{X}) \overset{\pi}\to C(\bd{X}) \to 0.
\end{equation}
Let \(\cl\upsilon\colon \C\to C(\cl{X})\) and \(\bd\upsilon\colon \C\to
C(\bd{X})\) be the unit maps.  The map~\(\bd\upsilon\) induces the obvious
embedding \(\Cred G\to C(\bd{X})\rcross G\), which we also denote by~\(u\).
We are going to study the induced map
\begin{equation}
  \label{eq:interesting_map_Kcross}
  u_* = \bd\upsilon_*\colon \K_*(\Cred G)\to\K_*(C(\bd{X})\rcross G).
\end{equation}

\section{Examples of boundary actions}
\label{sec:examples}

Many examples of boundary actions are special cases of two general
constructions: the visibility compactification of a \(\CAT(0)\) space and the
Gromov compactification of a hyperbolic space (see \cite{Bridson-Haefliger}
for both).

Let~\(X\) be a second countable, locally compact \(\CAT(0)\) space and
let~\(G\) act properly and isometrically on~\(X\).  We consider geodesic rays
\(\R_+\to X\) parametrised by arc length.  Two such rays are equivalent if
they are at bounded distance from each other.  The \emph{visibility boundary}
of~\(X\) is the set \(\bd{X}_\vis\) of equivalence classes of geodesic rays
in~\(X\) and the \emph{visibility compactification}~\(\cl{X}_\vis\) is \(X\cup
\bd{X}_\vis\).  This is a compactification of~\(X\) for a canonical compact
metrisable topology on \(\cl{X}_\vis\); it has the property that \(r(t)\)
converges towards \([r]\) for \(t\to\infty\) for any geodesic ray \(r\colon
\R_+\to X\).  The obvious action of~\(G\) on \(\cl{X}_\vis\) is continuous.

Let \(H\subseteq G\) be a compact subgroup.  Then~\(H\) has a
fixed-point~\(\xi_H\) in~\(X\).  For any \(x\in X\) there is a unique geodesic
segment connecting \(x\) and~\(\xi_H\).  We may contract the space~\(X\) along
these geodesics, so that~\(X\) is strongly contractible.  Similarly,
\(\cl{X}_\vis\) is strongly contractible because any point in \(\bd{X}_\vis\)
is represented by a unique geodesic ray emanating from~\(\xi_H\).

For instance, if~\(X\) is a simply connected Riemannian manifold of
non-positive sectional curvature, then~\(X\) is \(\CAT(0)\).  If \(\dim X=n\),
then there is a homeomorphism from \(\cl{X}_\vis\) onto a closed
\(n\)\nbd{}cell that identifies~\(X\) with the open \(n\)\nbd{}disk and
\(\bd{X}_\vis\) with \(S^{n-1}\).

Let~\(G\) be an almost connected Lie group whose connected component is linear
and reductive, and let \({K\subseteq G}\) be a maximal compact subgroup.  Then
the homogeneous space \(G/K\) with any \(G\)\nbd{}invariant Riemannian metric
is a \(\CAT(0)\) space and has non-positive sectional curvature
(\cite{Bridson-Haefliger}).  If~\(G\) is semi-simple and has rank~\(1\), then
the visibility boundary is equivalent to the Fürstenberg boundary \(G/P\)
of~\(G\), where \(P\) is a minimal parabolic subgroup.  If~\(G\) has higher
rank, then the Fürstenberg boundary of~\(G\) is not a boundary action in our
sense.  However, there are points in \(\bd{X}_\vis\) that are fixed by~\(P\).
Hence the unit map \(\C\to C(G/P)\) factors through \(\bd\upsilon\colon \C\to
C(\bd{X}_\vis)\).

Euclidean buildings and trees are \(\CAT(0)\) spaces as well.  For instance,
let~\(G\) be a reductive \(p\)\nbd{}adic group.  Then its affine Bruhat-Tits
building is a \(\CAT(0)\) space, on which~\(G\) acts properly and
isometrically (\cite{Tits:Corvallis}).  The visibility compactification of the
Bruhat-Tits building is equivalent to the Borel-Serre compactification in this
case (see \cite{Borel-Serre} and \cite{Schneider-Stuhler}*{Lemma IV.2.1}).
The relationship between the visibility boundary of the building and the
Fürstenberg boundary of~\(G\) is exactly as in the Lie group case.

Let~\(X\) be a \(\CAT(0)\) space on which~\(G\) acts properly and
isometrically and let \(\cl{X}_\vis\) be its visibility compactification.  The
\emph{limit set} \(\Lambda X_\vis\subseteq \bd{X}_\vis\) is the set of all
accumulation points in the boundary of \(G\cdot x\) for some \(x\in X\).  Its
definition is most familiar for classical hyperbolic space~\(\Hyp^n\) (see
\cite{Ratcliffe}*{Section 12}).  The limit set is independent of the choice
of~\(x\) and therefore \(G\)\nbd{}invariant.  Its complement \(X'\defeq
\cl{X}_\vis\setminus \Lambda X_\vis\) is called the \emph{ordinary set} of
\(G\).  It is strongly contractible for the same reason as \(\cl{X}_\vis\).
If the action of~\(G\) on~\(X'\) is proper, then \(\Lambda X_\vis\) is a
boundary action.  This is always the case for classical hyperbolic space
(\cite{Ratcliffe}).

Let~\(X\) be a (quasi-geodesic) hyperbolic metric space.  Two quasi-geodesic
rays in~\(X\) are considered equivalent if they have bounded distance.  The
\emph{Gromov boundary} \(\bd{X}\) of~\(X\) is defined as the set of
equivalence classes of quasi-geodesic rays in~\(X\).  There is a canonical
compact metrisable topology on \(\cl{X}\defeq X\cup\bd{X}\) so that this
becomes a compactification of~\(X\).  The construction is natural with respect
to quasi-isometric equivalence.  That is, if \(X\) and~\(Y\) are hyperbolic
metric spaces and \(f\colon X \to Y\) is a quasi-isometric equivalence,
then~\(f\) extends in a unique way to a map \(\cl{f}\colon \cl{X} \to
\cl{Y}\), whose restriction to the boundary \(\bd{f}\colon \bd{X}\to\bd{Y}\)
is a homeomorphism.

If~\(G\) is a word-hyperbolic group, then one may apply this construction to
the metric space underlying~\(G\) (with a word metric).  Since the action
of~\(G\) by left translation on itself is isometric, this action extends to an
action of~\(G\) by homeomorphisms of its boundary.  Now let \(X = P_d(G)\) be
the Rips complex of~\(G\) with parameter~\(d\).  This space is strongly
contractible for sufficiently large~\(d\) (\cite{Meintrup-Schick}).  We may
equip~\(X\) with a \(G\)\nbd{}invariant metric for which any orbit map \(G \to
X\) is a quasi-isometric equivalence.  Since hyperbolicity and the Gromov
boundary are invariant under quasi-isometric equivalence, \(X\) is itself
hyperbolic and there is a canonical \(G\)\nbd{}equivariant homeomorphism
\(\bd{G}\cong\bd{X}\).  It is shown in \cite{Rosenthal-Schutz} that~\(\cl{X}\)
is strongly contractible.  Hence \(G\times\bd{G}\to\bd{G}\) is a boundary
action.

Now we come to a completely different example of a boundary action.
Let~\(\Sigma_g\) be a closed Riemann surface of genus \(g\ge2\) and let~\(G\)
be a torsion-free subgroup of the \emph{mapping class group},
\[
\map\defeq \mathrm{Diff}(\Sigma_g)/\mathrm{Diff}_0(\Sigma_g).
\]
Let~\(\teich\) be the \emph{Teichmüller space} for~\(\Sigma_g\), and let
\(\bd{\teich}\) and \(\cl{\teich} \defeq \teich \cup \bd{\teich}\) be its
\emph{Thurston boundary} and \emph{Thurston compactification}, respectively
(\cite{fathi}).  Recall that~\(\bd{\teich}\) is the space of projective
geodesic laminations associated to a fixed hyperbolic metric on~\(\Sigma_g\).
It is well-known that \(\teich\) and~\(\cl{\teich}\) are contractible.  This
means that \(G\times \bd{\teich} \to \bd{\teich}\) is a boundary action
because we require~\(G\) to be torsion-free.

It seems plausible that \(\teich\) and~\(\cl{\teich}\) are strongly
contractible with respect to \(\map\) itself, so that \(\bd{\teich}\) would be
a boundary action of \(\map\).  That~\(\teich\) is strongly contractible
(hence a universal proper space) follows from the proof of the \emph{Nielsen
  realisation problem} by Steven Kerckhoff in \cite{kerck}.  The main result
is that every finite subgroup~\(H\) of \(\map\) fixes a point of~\(\teich\).
Moreover, any two points of~\(\teich\) are connected by a unique
\emph{earthquake path} (see \cites{thurston, kerck}).  If \(x\in\teich\) is
fixed by~\(H\), then the contraction of Teichmüller space along earthquake
paths emanating from~\(x\) provides an \(H\)\nbd{}equivariant contracting
homotopy for~\(\teich\).  It seems likely that this contracting homotopy
extends to one for the compactification~\(\cl{\teich}\), but a proof does not
seem to exist in the literature so far.

A group may admit more than one boundary action.  For example, the group
\(G=\PSL(2,\Z)\) has two natural models for \(\EG\), namely, the hyperbolic
plane~\(\Hyp^2\), on which it acts by Möbius transformations, and the
tree~\(X\) which corresponds to the free-product decomposition \(G \cong \Z/2
* \Z/3\) (see \cite{Serre:Trees}*{\S I.4}).  Whereas \(G\backslash X\) is
compact, \(G\backslash\Hyp^2\) is not.  Both models for \(\EG\) admit boundary
actions, which are completely different: we may compactify~\(\Hyp^2\) by the
circle at infinity \(\mathbb{P}^1(\R)\), and~\(X\) by its set of ends, which
is a Cantor set.  These two compactifications agree with the visibility
compactification or the Gromov compactification of~\(\Hyp^2\) and with the
Gromov compactification of~\(X\), respectively.  One can also verify directly
that they are boundary actions.

These two boundary actions of \(\PSL(2,\Z)\) are related by a well-known
\(G\)\nbd{}equivariant embedding of the tree~\(X\) in~\(\Hyp^2\) (see
\cite{Serre:Trees}*{\S I.4}).  This embedding extends to a continuous map
between the compactifications and hence yields a map \(\bd{X}\to\bd{\Hyp^2}\),
which is known to be two-to-one.  In general, it is not clear why different
boundary actions should be related in a similar fashion.

\section{Applying the Baum-Connes conjecture}
\label{sec:apply_BC}

We recall the definition of the Baum-Connes assembly map
(see~\cite{Baum-Connes-Higson}).  Let \(\EG\) be a second countable, locally
compact model for the universal proper \(G\)\brd{}space.  Write
\(\EG=\bigcup_{n\in\N} \EG_n\) for some increasing sequence of
\(G\)\nbd{}compact subsets \(\EG_n\subseteq\EG\).  The maps
\(\EG_n\to\EG_{n+1}\) are proper, so that we get an associated projective
system of \(G\)\brd{}\(C^*\)\brd{}algebras
\(\bigl(C_0(\EG_n)\bigr)_{n\in\N}\).  Let
\[
\Ktop_*(G,A) \defeq \varinjlim \KK^G_*(C_0(\EG_n),A).
\]
The \emph{Baum-Connes assembly map} is the composite map
\[
\mu_{G,A}\colon \Ktop_*(G,A)
\to \varinjlim \KK_*(C_0(\EG_n)\rcross G, A\rcross G)
\to \K_*(A\rcross G),
\]
where the first map is descent and the second map is the Kasparov product with
a certain natural class in \(\varprojlim \K_0(C_0(\EG_n)\rcross G)\).

Let \(X\), \(\cl{X}\) and \(\bd{X}\) be as in Definition
\ref{def:boundary_action}.  The Baum-Connes conjecture always holds for proper
coefficient algebras and, especially, for \(C_0(X)\).  \emph{We assume from
  now on that~\(G\) satisfies the Baum-Connes conjecture with coefficients
  \(\C\) and \(C(\bd{X})\)}.  That is, the vertical maps in the diagram
\[
\xymatrix{
  \Ktop_*(G,\C) \ar[r]^-{\bd\upsilon_*} \ar[d]^{\mu_{G,\C}} &
  \Ktop_*(G,C(\bd{X})) \ar[d]^{\mu_{G,C(\bd{X})}} \\
  \K_*(\Cred G) \ar[r]^-{\bd\upsilon_*} &
  \K_*(C(\bd{X})\rcross G)
}
\]
are isomorphisms.  Thus the map~\eqref{eq:interesting_map_Kcross} that we are
interested in is equivalent to the map
\begin{equation}
  \label{eq:interesting_map_Ktop}
  \bd\upsilon_*\colon \Ktop_*(G,\C)\to\Ktop_*(G,C(\bd{X})).
\end{equation}

Our assumption on the Baum-Connes conjecture is known to be valid in many
examples.  Closed subgroups of \(\Isom(\Hyp^n)\) satisfy it because they even
satisfy the Baum-Connes conjecture with arbitrary coefficients
(\cite{Kasparov:Warsaw}).  The same holds for closed subgroups of other
semi-simple Lie groups of rank~\(1\) by \cites{Kasparov-Julg:SUn1, Julg:Spn1}.
Word-hyperbolic groups satisfy the assumption as well: the Baum-Connes
conjecture with trivial coefficients is proved
in~\cite{Mineyev-Yu:BC_Hyperbolic}, and the Baum-Connes conjecture for the
coefficients \(C(\bd{G})\) follows from~\cite{Tu:Amenable} because the action
of~\(G\) on~\(\bd{G}\) is amenable.

We will exclusively deal with the map in~\eqref{eq:interesting_map_Ktop} in
the following.  We \emph{only} need the Baum-Connes conjecture to relate it
to~\eqref{eq:interesting_map_Kcross}.

It is shown in~\cite{Kasparov-Skandalis:Bolic} that \(\Ktop_*\) satisfies
excision for arbitrary extensions of \(G\)\nbd{}\(C^*\)\brd{}algebras.
Hence~\eqref{eq:G_space_ext} gives rise to a six term exact sequence
\begin{equation}
  \label{eq:K_boundary_ext}
  \begin{aligned}
  \xymatrix{
    \Ktop_0(G,C_0(X)) \ar[r]^{\iota_*} &
    \Ktop_0(G,C(\cl{X})) \ar[r]^{\pi_*} &
    \Ktop_0(G,C(\bd{X})) \ar[d]^{\delta} \\
    \Ktop_1(G,C(\bd{X})) \ar[u]_{\delta} &
    \Ktop_1(G,C(\cl{X})) \ar[l]_{\pi_*} &
    \Ktop_1(G,C_0(X)).  \ar[l]_{\iota_*}
  }
  \end{aligned}
\end{equation}
We are going to modify it in several steps.

\begin{lemma}
  \label{lem:Ktop_clX}
  The map \(\cl\upsilon_*\colon\Ktop_*(G,\C)\to \Ktop_*(G,C(\cl{X}))\) is an
  isomorphism if~\(\cl{X}\) is strongly contractible.
\end{lemma}

\begin{proof}
  Since~\(\cl{X}\) is strongly contractible, \([\cl\upsilon]\) is invertible
  in \(\KK^H_0(\C,C(\cl{X}))\) for all compact subgroups \(H\subseteq G\).
  That is, \(\cl\upsilon\) is a \emph{weak equivalence} in the notation of
  \cite{Meyer-Nest}.  It is shown in
  \cites{Chabert-Echterhoff-Oyono:Going_down,Meyer-Nest} that such maps induce
  isomorphisms on \(\Ktop_*(G,\blank)\).
\end{proof}

Plugging the isomorphism of Lemma~\ref{lem:Ktop_clX}
into~\eqref{eq:K_boundary_ext}, we get an exact sequence
\begin{equation}  \label{eq:K_boundary_ext_II}
  \begin{aligned}
  \xymatrix{
    \Ktop_0(G,C_0(X)) \ar[r]^-{\cl\upsilon_*^{-1}\iota_*} &
    \Ktop_0(G,\C) \ar[r]^-{\bd\upsilon_*} &
    \Ktop_0(G,C(\bd{X})) \ar[d]^{\delta} \\
    \Ktop_1(G,C(\bd{X})) \ar[u]_{\delta} &
    \Ktop_1(G,\C) \ar[l]_-{\bd\upsilon_*} &
    \Ktop_1(G,C_0(X)). \ar[l]_-{\cl\upsilon_*^{-1}\iota_*}
  }
  \end{aligned}
\end{equation}
It contains the map~\eqref{eq:interesting_map_Ktop} that we are interested in.

\begin{example}
  \label{exa:fixed_boundary_point}
  Suppose that the action of~\(G\) on~\(\bd{X}\) admits a fixed-point~\(\xi\).
  Then evaluation at~\(\xi\) provides a section for \(\cl\upsilon\).  Since
  this evaluation homomorphism annihilates \(C_0(X)\subseteq C(\cl{X})\), we
  get \(\cl\upsilon_*^{-1}\iota_*=0\).  Therefore, the long exact
  sequence~\eqref{eq:K_boundary_ext_II} splits into two short exact sequences.
\end{example}

\begin{proposition}
  \label{pro:fixed_point_bd}
  Let~\(G\) be a locally compact group with non-compact centre.  Suppose
  that~\(X\) is \(G\)\nbd{}compact and that~\(\cl{X}\) is admissible in the
  sense of~\cite{Higson:Bivariant}, that is, compatible with the coarse
  geometric structure of~\(X\).  Then~\(\bd{X}\) contains a fixed-point
  for~\(G\).  Hence \(\cl\upsilon_*^{-1}\iota_*=0\)
  in~\eqref{eq:K_boundary_ext_II}.
\end{proposition}

\begin{proof}
  Let \(x\in X\) and let \((g_i)_{i\in\N}\) be a sequence in the centre
  of~\(G\) that leaves any compact subset of~\(G\).  Since~\(\cl{X}\) is
  compact, we may assume that the sequence \(g_ix\) converges towards some
  \(\xi\in\cl{X}\).  Since the sequences \((g_ix)\) and \((g\cdot
  g_ix)=(g_igx)\) are uniformly close for any \(g\in G\), the compatibility of
  the coarse structure with the compactification implies that they have the
  same limit point in~\(\cl{X}\).  Thus \(g\xi=\xi\) for all \(g\in G\).
\end{proof}

Of course, the map \(\cl\upsilon_*^{-1}\iota_*\)
in~\eqref{eq:K_boundary_ext_II} is non-zero is general.  The following
notation is needed in order to describe it.

Let~\(Y\) be a locally compact \(G\)\nbd{}space.  We define a graded Abelian
group \(\RKK^G_*(Y;A,B)\) for a pair of \(G\)\nbd{}\(C^*\)\brd{}algebras
\(A,B\) as in~\cite{Kasparov:Novikov}; its cycles are cycles \((\Hilm, F)\)
for \(\KK^G_*(C_0(Y)\hot A, C_0(Y)\hot B)\) that satisfy the additional
condition that the left and right \(C_0(Y)\)\brd{}actions on~\(\Hilm\) agree.
We may think of these cycles as \(G\)\nbd{}equivariant continuous families of
\(\KK^G_*(A,B)\)-cycles parametrised by~\(Y\).  This makes it clear that there
is a natural map
\begin{equation}
  \label{eq:forget_Y}
  \RKK^G_*(Y; A,B) \to \KK^G_*(C_0(Y)\hot A, C_0(Y)\hot B)
\end{equation}
which forgets the \(Y\)\nbd{}structure and a natural map
\begin{equation}
  \label{eq:inflate_Y}
  p_Y^*\colon\KK^G_*(A,B) \to \RKK^G_*(Y; A,B)  
\end{equation}
which sends a Kasparov cycle \((\Hilm, F)\) to the constant family of Kasparov
cycles \((C_0(Y)\hot \Hilm, 1\hot F)\).

\begin{definition}
  \label{def:diagonalrestrictionclass}
  We let \(\Diagonal_Y \in \RKK^G_0(Y; C_0(Y),\C)\) be the class of the
  \(Y\cross G\)\brd{}equivariant \(*\)\nbd{}homomorphism \(C_0(Y\times Y)\to
  C_0(Y)\) that is induced by the diagonal embedding \(Y\to Y\times Y\).
\end{definition}

Recall that the space~\(X\) is a universal proper \(G\)\nbd{}space.  It is
shown in \cite{Meyer-Nest}*{Section 7} that the map in~\eqref{eq:inflate_Y}
for \(Y=X\) is an isomorphism whenever~\(A\) is a proper
\(G\)\nbd{}\(C^*\)\brd{}algebra.  Hence
\(\Diagonal_X\in\RKK^G_0(X;C_0(X),\C)\) has a pre-image in
\(\KK^G_0(C_0(X),\C)\).  Anticipating a little, we denote this pre-image by
\(\Eul_X\).  This agrees with our official definition of the abstract Euler
characteristic in Definition \ref{def:abstract_Euler} by Lemma
\ref{lem:pX_invertible} and Proposition \ref{pro:abstract_dual_for_EG}.

\begin{proposition}
  \label{pro:abstract_Gysin}
  Let~\(G\) be a locally compact group and let \(\bd{X}=\cl{X}\setminus X\) be
  a boundary action of~\(G\); that is, \(X\) is a strongly contractible proper
  \(G\)\nbd{}space and~\(\cl{X}\) is a strongly contractible compact
  \(G\)\nbd{}space containing~\(X\) as an open \(G\)\nbd{}invariant subspace.
  Then there is an exact sequence
  \[\xymatrix{
    \Ktop_0(G,C_0(X)) \ar[r]^-{\Eul_X} &
    \Ktop_0(G,\C) \ar[r]^-{\bd\upsilon_*} &
    \Ktop_0(G,C(\bd{X})) \ar[d]^{\delta} \\
    \Ktop_1(G,C(\bd{X})) \ar[u]_{\delta} &
    \Ktop_1(G,\C) \ar[l]_-{\bd\upsilon_*} &
    \Ktop_1(G,C_0(X)), \ar[l]_-{\Eul_X}
  }
  \]
  where \(\Eul_X\) denotes the Kasparov product with
  \(\Eul_X\in\KK^G_0(C_0(X),\C)\).
\end{proposition}

\begin{proof}
  It only remains to identify the maps \(\cl\upsilon_*^{-1}\iota_*\)
  in~\eqref{eq:K_boundary_ext_II} with \(\Eul_X\).  Recall that the map
  \(p_X^*\) in~\eqref{eq:inflate_Y} is an isomorphism if~\(A\) is proper.
  Therefore, \(p_X^*\) induces an isomorphism
  \[
  \Ktop_*(G,A) \congto \varinjlim \RKK^G_*(X;C_0(\EG_n),A).
  \]
  Thus the Kasparov product in \(\RKK^G_*(X)\) gives rise to natural bilinear
  maps
  \begin{equation}
    \label{eq:bullet_action}
    \Ktop_*(G,A) \otimes \RKK^G_0(X;A,B) \to \Ktop_*(G,B),
    \qquad x\otimes y \mapsto x\bullet y.
  \end{equation}

  It is shown in \cite{Meyer-Nest}*{Section 7} that \(p_X^*(f)\) for
  \(f\in\KK^G_*(A,B)\) is invertible if and only if~\(f\) is a weak
  equivalence.  Therefore, \(p_X^*[\cl\upsilon]\) is invertible in
  \(\RKK^G_0(X;\C,C(\cl{X}))\).  The homomorphism \(C_0(X)\to
  C_0(X\times\cl{X})\) induced by the coordinate projection \(X\times\cl{X}\to
  X\) is a representative for \(p_X^*[\cl\upsilon]\).  Since the diagonal
  embedding \(X\to X\times X\subseteq X\times\cl{X}\) is a section for the
  coordinate projection, the element in \(\RKK^G_0(X;C(\cl{X}),\C)\) associated
  to the diagonal embedding is inverse to \(p_X^*[\cl\upsilon]\).  Hence
  \(p_X^*[\iota]\hot_{X,C(\cl{X})} p_X^*[\cl\upsilon]^{-1}=\Delta_X\).  It
  follows easily from the definition of the product
  in~\eqref{eq:bullet_action} that
  \[
  \iota_*(x)= x \hot_{C_0(X)} [\iota] = x\bullet p_X^*[\iota],
  \qquad
  \cl\upsilon_*(y)= y \hot [\cl\upsilon] = y\bullet p_X^*[\cl\upsilon],
  \]
  for all \(x\in\Ktop_*(G,C_0(X))\), \(y\in\Ktop_*(G,\C)\).  This implies
  \[
  \cl\upsilon_*^{-1}\iota_*(x)
  = x\bullet (p_X^*[\iota] \hot_{X,C(\cl{X})} p_X^*[\cl\upsilon]^{-1})
  = x\bullet \Diagonal_X
  = x\hot_{C_0(X)} p_X^{-1}(\Diagonal_X)
  \]
  because~\(p_X\) is invertible and the map in~\eqref{eq:bullet_action} is
  natural.
\end{proof}

\section{Abstract Euler characteristics via Kasparov duality}
\label{sec:Kasparov_duality}

The element \(\Eul_X\in\KK^G_0(C_0(X),\C)\) that appears in the Gysin sequence
in Proposition \ref{pro:abstract_Gysin} is so far only defined if~\(X\) is a
universal proper \(G\)\nbd{}space; furthermore, it is not clear how it should
be computed.  In this section, we extend its definition to a more general
class of \(G\)\nbd{}spaces, using a formulation of Poincaré duality due to
Gennadi Kasparov (\cite{Kasparov:Novikov}*{Section 4}).  In the following
sections, we will compute \(\Eul_X\) using this alternative definition.

\begin{definition}
  \label{def:abstract_dual}
  Let~\(X\) be any locally compact \(G\)\nbd{}space (we require neither
  properness nor strong contractibility).  Let \(n\in\Z\).  Let~\(\sour\) be a
  (possibly \(\Ztwo\)\brd{}graded) \(G\)\nbd{}\(C^*\)-algebra, and let
  \(\Theta\in\RKK^G_n(X;\C,\sour)\).

  We call \((\sour,\Theta)\) an (\(n\)\nbd{}dimensional) \emph{abstract dual}
  for~\(X\) if the map
  \begin{equation}
    \label{eq:def_PD}
    \PD\colon \RKK^G_{*-n}(Y;A_1\hot\sour,A_2) \to \RKK^G_*(X\times Y;A_1,A_2),
    \quad f\mapsto \Theta\hot_\sour f,
  \end{equation}
  is an isomorphism for all \(G\)\nbd{}spaces~\(Y\) and all
  \(G\)-\(C^*\)-algebras \(A_1,A_2\).
\end{definition}

Here we use the Kasparov product
\begin{multline*}
  \hot_\sour\colon \RKK^G_i(X;A,B\hot\sour) \times \RKK^G_j(Y;A'\hot\sour,B')
  \\ \to \RKK^G_{i+j}(X\times Y;A\hot A',B\hot B')
\end{multline*}
(see~\cite{Kasparov:Novikov}).  Observe that~\eqref{eq:def_PD} is the most
general form for a natural transformation \(\RKK^G_{*-n}(Y;A\hot\sour,B) \to
\RKK^G_*(X\times Y;A,B)\) that is compatible with Kasparov products in the
sense that
\begin{equation}
  \label{eq:PD_and_product}
  \PD(f_1\hot_B f_2) = \PD(f_1)\hot_B f_2
  \qquad \text{in \(\RKK^G_{i+j}(Y;A_1\hot A_3,A_2\hot A_4)\)}
\end{equation}
for all \(f_1\in \RKK^G_{i-n}(Y;A_1\hot\sour,A_2\hot B)\),
\(f_2\in\RKK^G_j(Y;A_3\hot B,A_4)\).  Since exterior products are graded
commutative, we also get
\begin{equation}
  \label{eq:PD_and_product_II}
  \PD(f_1\hot_B f_2) = (-1)^{in} f_1\hot_B \PD(f_2)
  \qquad \text{in \(\RKK^G_{i+j}(Y;A_1\hot A_3,A_2\hot A_4)\)}
\end{equation}
for all \(f_1\in \RKK^G_i(Y;A_1,A_2\hot B)\), \(f_2\in\RKK^G_{j-n}(Y;A_3\hot
B\hot\sour,A_4)\); both sides are equal to \((\Theta\hot f_1)\hot_{\sour\hot
  B} f_2\).

The space~\(Y\) does not play any serious role.  We have put it into our
definitions because Kasparov works in this generality in
\cite{Kasparov:Novikov}*{Theorem 4.9}.  The dimension~\(n\) is not
particularly important either because we can always reduce to the case \(n=0\)
by a suspension.

We are mainly interested in the case of complex \(C^*\)\nbd{}algebras and
therefore only formulate Gysin sequences in this case.  However, the purely
formal arguments in this section are independent of Bott periodicity and
therefore also work in the real and ``real'' cases.  This is why we are
careful to distinguish \(\KK_n\) and \(\KK_{-n}\) in our notation.  Of course,
in the real case one has to replace~\(\C\) by~\(\R\) everywhere and use
real-valued function spaces.

Later, we shall introduce further structure in order to write down the inverse
map \(\PD^{-1}\) more concretely, which is important for applications.
However, this additional structure involves some choices.  Since the abstract
Euler characteristic is supposed to be independent of the dual, we discuss the
formal aspects of the duality in the situation of Definition
\ref{def:abstract_dual}.

\begin{remark}
  \label{rem:space_without_dual}
  There exist spaces that do not possess an abstract dual, even for
  trivial~\(G\).  If~\(X\) is compact then
  \(\RKK_*(X;A,B)\cong\KK_*(A,C(X)\hot B)\) for all \(A,B\).  Hence an
  abstract dual for~\(X\) is nothing but a \(\KK\)-dual for \(C(X)\).  Since
  \(C(X)\) belongs to the bootstrap category, it has a \(\KK\)-dual if and
  only if \(\K_*(C(X))\) is finitely generated (recall that all spaces are
  assumed second countable).  This fails, for example, if~\(X\) is a Cantor
  set.
\end{remark}

\begin{definition}
  \label{def:abstract_Euler}
  Let~\(X\) be a \(G\)\nbd{}space with an abstract dual \((\sour,\Theta)\).
  Let \(\bar\Theta\in\KK^G_n(C_0(X),C_0(X)\hot\sour)\) be the image
  of~\(\Theta\) under the forgetful map~\eqref{eq:forget_Y}; let
  \(\Diagonal_X\in\RKK^G_0(X;C_0(X),\C)\) be induced by the diagonal embedding
  as in Definition \ref{def:diagonalrestrictionclass}; thus
  \(\PD^{-1}(\Diagonal_X)\in\KK^G_{-n}(C_0(X)\hot\sour,\C)\).  We call
  \[
  \Eul_X \defeq \bar\Theta\hot_{C_0(X)\hot\sour} \PD^{-1}(\Diagonal_X)
  \in \KK^G_0(C_0(X),\C)
  \]
  the \emph{\(G\)\nbd{}equivariant abstract Euler characteristic of~\(X\)}.
\end{definition}

This name will be justified by the examples in the following sections.

Our first task is to analyse the uniqueness of abstract duals and to show that
\(\Eul_X\) does not depend on their choice.  We consider the slightly more
complicated issue of functoriality right away.

Let \(X\) and~\(X'\) be two \(G\)\nbd{}spaces with abstract duals
\((\sour,\Theta)\) and \((\sour',\Theta')\) of dimension \(n\) and~\(n'\),
respectively, and let \(\PD\) and \(\PD'\) be the associated duality
isomorphisms.  Let \(f\colon X\to X'\) be a continuous \(G\)\nbd{}map; we do
\emph{not} require~\(f\) to be proper.  Then~\(f\) induces natural maps
\begin{equation}
  \label{eq:finduced}
  f^*\colon \RKK^G_*(X'\times Y;A,B) \to \RKK^G_*(X\times Y;A,B)
\end{equation}
for all \(Y,A,B\).  Hence we get \(f^*\Theta'\in\RKK^G_{n'}(X;\C,\sour')\)
and
\[
\alpha_f\defeq \PD^{-1}(f^*\Theta') \in \KK^G_{n'-n}(\sour,\sour').
\]
Equivalently, \(\Theta\hot_\sour \alpha_f = f^*\Theta'\); this property
characterises~\(\alpha_f\) uniquely and implies
\begin{equation}
  \label{eq:PD_functorial}
  \PD(\alpha_f \hot_{\sour'} h)
  = f^*\Theta' \hot_{\sour'} h
  = f^*\bigl(\PD'(h)\bigr)
  \qquad\text{in \(\RKK^G_*(X;A,B)\)}
\end{equation}
for all \(h\in\KK^G_{*-n'}(A\hot\sour',B)\).  The map \(f\mapsto\alpha_f\) is
a covariant functor in the following sense.  If \(f=\ID_X\) and
\((\sour,\Theta)=(\sour',\Theta')\), then \(\alpha_\ID=1_\sour\).  Given
composable maps \(f\colon X\to X'\), \(f'\colon X'\to X''\) and abstract duals
for the \(G\)\nbd{}spaces \(X,X',X''\), we get \(\alpha_{f'\circ
  f}=\alpha_f\hot_{\sour'} \alpha_{f'}\).

If two maps \(f_1,f_2\colon X\to X'\) are \(G\)\nbd{}equivariantly homotopic,
then they induce the same maps \(f_1^*=f_2^*\) in~\eqref{eq:finduced}.  Hence
\(\alpha_{f_1}=\alpha_{f_2}\).  Moreover, if~\(f\) is a \(G\)\nbd{}homotopy
equivalence then~\eqref{eq:finduced} is bijective for all \(Y,A,B\).  By
functoriality of~\(\alpha_f\), we conclude that~\(\alpha_f\) is invertible
if~\(f\) is a \(G\)\nbd{}homotopy equivalence.  In the special case \(f=\ID\),
we get a canonical \(\KK^G\)\nbd{}equivalence between two duals
\((\sour,\Theta)\), \((\sour',\Theta')\) for the same space~\(X\).

\begin{proposition}
  \label{pro:Eul_functorial}
  Let \(f\colon X\to X'\) be a \(G\)\nbd{}homotopy equivalence and proper; we
  do not assume its homotopy inverse to be proper.  We denote the class of the
  induced map \(f^*\colon C_0(X')\to C_0(X)\) in \(\KK^G_0(C_0(X'),C_0(X))\)
  by \([f^*]\).  Then
  \[
  [f^*] \hot_{C_0(X)} \Eul_X=\Eul_{X'}.
  \]
  The abstract Euler characteristic is independent of the choice of the
  abstract dual.
\end{proposition}

\begin{proof}
  Let \(\Diagonal_{X'}\colon C_0(X'\times X')\to C_0(X')\) be the diagonal
  restriction homomorphism.  Then
  \(f^*(\Diagonal_{X'})\in\RKK^G_0(X;C_0(X'),\C)\) is the class of the
  \(G\)\nbd{}equivariant \(*\)\nbd{}homomorphism induced by the map
  \((\ID,f)\colon X\to X\times X'\).  Thus \(f^*(\Diagonal_{X'})=
  [f^*]\hot_{C_0(X)} \Diagonal_X\).  Since the map in~\eqref{eq:finduced} is
  bijective and natural with respect to the Kasparov product, this is
  equivalent to
  \[
  \Diagonal_{X'}
  = (f^*)^{-1}([f^*]\hot_{C_0(X)} \Diagonal_X)
  = [f^*]\hot_{C_0(X)} (f^*)^{-1}(\Diagonal_X).
  \]
  Define \(\alpha_f\in\KK^G_{n'-n}(\sour,\sour')\) as above.  Equation
  \eqref{eq:PD_functorial} is equivalent to
  \(\alpha_f^{-1}\hot_\sour\PD^{-1}(h) = (\PD')^{-1}(f^*)^{-1}(h)\) for all
  \(h\in\RKK^G_*(X;A,B)\).

  We shall use the forgetful functors defined in~\eqref{eq:forget_Y} for the
  spaces \(X\) and~\(X'\) and denote them by \(h\mapsto\overline h\).  They
  satisfy the compatibility relation
  \[
  [f^*] \hot_{C_0(X)} \overline{f^* h} = \overline{h} \hot_{C_0(X')} [f^*]
  \]
  in \(\KK^G_*(C_0(X')\hot A,C_0(X)\hot B)\) for all \(h\in\RKK^G_*(X';A,B)\);
  these Kasparov products are comparatively easy to compute because~\([f^*]\)
  is represented by a \(*\)\nbd{}homomorphism.  Now we compute
  \begin{multline*}
    [f^*] \hot_{C_0(X)} \Eul_X = [f^*] \hot_{C_0(X)} \overline{\Theta}
    \hot_{C_0(X)\hot\sour} \PD^{-1}(\Diagonal_X)
    \\ = [f^*] \hot_{C_0(X)} \overline{f^*\Theta'} \hot_{\sour'} \alpha_f^{-1}
    \hot_{C_0(X)\hot\sour} \PD^{-1}(\Diagonal_X)
    \\ = \overline{\Theta'} \hot_{C_0(X')} [f^*] \hot_{C_0(X)\hot\sour'}
    (\PD')^{-1}\circ (f^*)^{-1}(\Diagonal_X)
    \\ = \overline{\Theta'} \hot_{C_0(X')\hot\sour'}
    (\PD')^{-1}\bigl([f^*]\hot_{C_0(X)} (f^*)^{-1}(\Diagonal_X)\bigr)
    \\ = \overline{\Theta'} \hot_{C_0(X')\hot\sour'} (\PD')^{-1}(\Diagonal_{X'})
    = \Eul_{X'}.
  \end{multline*}
  We use~\eqref{eq:PD_and_product_II} in the step from the third to the fourth
  line.

  Especially, if \(f=\ID\) then the computation above shows that \(\Eul_X\)
  does not depend on the choice of the abstract dual.
\end{proof}

\begin{proposition}
  \label{pro:abstract_dual_for_EG}
  If~\(X\) is a universal proper \(G\)\nbd{}space, then~\(X\) has an abstract
  dual.
\end{proposition}

\begin{proof}
  Let \(\Dirac\in\KK^G_0(\ADir,\C)\) be a Dirac morphism in the notation of
  \cite{Meyer-Nest}.  Since~\(\Dirac\) is a weak equivalence,
  \(p_X^*(\Dirac)\) is invertible; let \(\Theta\in\RKK^G_0(X;\C,\ADir)\) be
  its inverse.  In \cite{Meyer-Nest}*{Theorem 7.1}, we may take
  \(\pi\in\KK^G_0(\tilde{A}, A)\) to be
  \(1_A\otimes\Dirac\in\KK^G_0(A\otimes\ADir,A)\).  Then
  \cite{Meyer-Nest}*{Theorem 7.1} implies that the map \(\PD\) that we get
  from~\(\Theta\) is an isomorphism.
\end{proof}

Assume now that \(X\) and~\(X'\) are \emph{\(G\)\nbd{}compact} universal
proper \(G\)\nbd{}spaces.  Then they are \(G\)\nbd{}homotopy equivalent in a
canonical way, so that their abstract duals are canonically
\(\KK^G\)\nbd{}equivalent.  Any continuous \(G\)\nbd{}map between \(X\)
and~\(X'\) is proper because both spaces are \(G\)\nbd{}compact.  Hence
\(C_0(X)\) and \(C_0(X')\) are \(\KK^G\)\nbd{}equivalent.  Moreover, we have a
natural isomorphism \(\Ktop_*(G)\defeq\Ktop_*(G,\C)\cong \KK^G_*(C_0(X),\C)\).
It follows from Proposition~\ref{pro:Eul_functorial} that the abstract Euler
characteristics of \(X\) and~\(X'\) agree as elements of \(\Ktop_0(G)\).
Hence we may give the following definition:

\begin{definition}
  \label{def:abstract_Euler_group}
  Let~\(G\) be a locally compact group that has a \(G\)\nbd{}compact universal
  proper \(G\)\nbd{}space~\(X\).  Identify
  \(\KK^G_*(C_0(X),\C)\cong\Ktop_*(G)\) and view the abstract Euler
  characteristic \(\Eul_X\) as an element of \(\Ktop_0(G)\).  We denote the
  result by \(\Eul_\EG\) and call it the \emph{abstract Euler characteristic
    of~\(G\)}.
\end{definition}

If~\(X'\) is any universal proper \(G\)\nbd{}space and~\(X\) is a
\(G\)\nbd{}compact universal proper \(G\)\nbd{}space, then any \(G\)\nbd{}map
\(f\colon X\to X'\) is proper and a \(G\)\nbd{}homotopy equivalence.  Hence
Proposition \ref{pro:Eul_functorial} yields \(\Eul_{X'} = f^*(\Eul_\EG)\).

Let~\(X\) again be an arbitrary \(G\)\nbd{}space with an abstract dual
\((\sour,\Theta)\).  Then we define \(D\in\KK^G_{-n}(\sour,\C)\) by
\begin{equation}
  \label{eq:def_D}
  \PD(D) \defeq \Theta \hot_\sour D = 1_\C
  \qquad\text{in \(\RKK^G_0(X;\C,\C)\).}
\end{equation}

\begin{lemma}
  \label{lem:pX_invertible}
  Let~\(X\) be a \(G\)\nbd{}space with an abstract dual \((\sour,\Theta)\).
  Then~\(\Theta\) is invertible if and only if
  \(p_X^*(D)\in\RKK^G_{-n}(X;\sour,\C)\) is invertible, if and only if the map
  \begin{equation}
    \label{eq:pX_invertible_I}
    p_X^*\colon \RKK^G_*(Y;\sour\hot A,B)\to\RKK^G_*(X\times Y;\sour\hot A,B)    
  \end{equation}
  is invertible for all \(Y,A,B\).  In this case, the map
  \begin{equation}
    \label{eq:pX_invertible_II}
    p_X^*\colon \RKK^G_*(Y;C_0(X)\hot A,B)\to
    \RKK^G_*(X\times Y;C_0(X)\hot A,B)
  \end{equation}
  is invertible as well, and \(\Eul_X= (p_X^*)^{-1}(\Diagonal_X)\) in
  \(\KK^G_0(C_0(X),\C)\).
\end{lemma}

\begin{proof}
  Since \(\Theta\hot_\sour D=1_\C\), \(\Theta\) is a left inverse for
  \(p_X^*(D)\) with respect to the Kasparov composition product in
  \(\RKK^G_*(X)\).  Hence one is invertible if and only if the other is, and
  they are inverse to each other in that case.  By hypothesis, \(\PD(f)\defeq
  \Theta \hot_\sour p_X^*(f)\) defines an invertible map on
  \(\RKK^G_*(Y;A\hot\sour,B)\) for all \(Y,A,B\).  If~\(\Theta\) is
  invertible, then so is the Kasparov product with~\(\Theta\) that appears in
  \(\PD\); hence the map in~\eqref{eq:pX_invertible_I} is invertible.
  Conversely, if the map in~\eqref{eq:pX_invertible_I} is invertible, then the
  Kasparov product with~\(\Theta\) is invertible as a map
  \(\RKK^G_{*-n}(X;\sour,B)\to\RKK^G_*(X;\C,B)\) for all~\(B\).  This implies
  invertibility of~\(\Theta\).

  If~\(\Theta\) is invertible, so is
  \(\bar\Theta\in\KK^G_n(C_0(X),C_0(X)\hot\sour)\).  Therefore, the map
  in~\eqref{eq:pX_invertible_II} is equivalent to one of the
  form~\eqref{eq:pX_invertible_I}; thus it is invertible as well.  Its inverse
  is computed as follows.  We have
  \[
  p_X^*(f)
  = (\Theta \hot_\sour D) \hot f
  = \PD(D\hot f)
  = \PD(\bar\Theta^{-1}\hot_{C_0(X)} f).
  \]
  for all \(f\in\KK^G_*(C_0(X)\hot A,B)\) because \(1_{C_0(X)}\hot D\)
  and~\(\bar\Theta\) are inverse to each other and \(\Theta\hot_\sour
  D=1_\C\).  Hence
  \begin{equation}
    \label{eq:pXinverse}
    (p_X^*)^{-1}(f') = \bar\Theta\hot_{C_0(X)\hot\sour} \PD^{-1}(f')    
  \end{equation}
  for all \(f'\in\RKK^G_*(X;C_0(X)\hot A,B)\).  In particular,
  \((p_X^*)^{-1}(\Diagonal_X)=\Eul_X\) as desired.
\end{proof}

Proposition \ref{pro:abstract_dual_for_EG} and Lemma \ref{lem:pX_invertible}
show that the definition of \(\Eul_X\) in Section~\ref{sec:apply_BC} is a
special case of Definition \ref{def:abstract_Euler}.

The following discussion has the purpose of motivating the definition of a
Kasparov dual by explaining how it is related to an abstract dual.  Define
\(\comul\in\KK^G_n(\sour,\sour\hot\sour)\) by
\begin{equation}
  \label{eq:def_comul}
  \PD(\comul) = \Theta \hot_\sour \comul = \Theta\hot_X\Theta
  \in\RKK^G_{2n}(X;\C,\sour\hot\sour).
\end{equation}
Let~\(\Phi_P\) for a \(G\)\nbd{}\(C^*\)-algebra~\(P\) be the flip automorphism
\[
\Phi_P\colon P\hot P\to P\hot P,
\qquad x_1\hot x_2 \mapsto
(-1)^{\abs{x_1}\cdot \abs{x_2}} x_2\hot x_1,
\]
where \(\abs{x}\in\Ztwo\) denotes the degree of~\(x\).  The sign only occurs
if~\(P\) is \(\Ztwo\)\nbd{}graded.  Recall also the class
\(D\in\KK^G_n(\sour,\C)\) defined in~\eqref{eq:def_D}.

\begin{lemma}
  \label{lem:abstract_dual_coalgebra}
  The maps \(D\) and~\(\comul\) satisfy the conditions for a graded
  cocommutative, counital coalgebra object in \(\KK^G_*\); that is,
  \begin{gather*}
    \comul\hot_{\sour\hot\sour}\Phi_\sour= (-1)^n\comul
    \\
    \comul\hot_{\sour\hot \sour}(\comul\hot 1_\sour)
    = \comul \hot_{\sour\hot\sour}(1_\sour\hot\comul)
    \\
    (-1)^n\comul \hot_{\sour\hot\sour} (D\hot 1_\sour)
    = 1_\sour
    = \comul \hot_{\sour\hot\sour} (1_\sour\hot D).
  \end{gather*}
  These equalities hold in the groups \(\KK^G_n(\sour,\sour\hot\sour)\),
  \(\KK^G_{2n}(\sour,\sour^{\hot 3})\), and \(\KK^G_0(\sour,\sour)\),
  respectively.
\end{lemma}

\begin{proof}
  It is well-known that the exterior product in \(\RKK^G_*(X)\) is graded
  commutative.  Especially, \((\Theta\hot_X\Theta)
  \hot_{\sour\hot\sour}\Phi_\sour= (-1)^n \Theta\hot_X\Theta\).  This is
  equivalent to the cocommutativity of~\(\comul\) because \(\PD\) is
  compatible with Kasparov products and bijective.  One checks easily that
  \(\PD\) maps both \(\comul\hot_{\sour\hot \sour} (\comul\hot1_\sour)\) and
  \(\comul \hot_{\sour\hot\sour} (1_\sour\hot\comul)\) to \(\Theta \hot_X
  \Theta\hot_X \Theta\) in \(\RKK^G_{3n}(X; \sour,\sour^{\hot 3})\).
  Thus~\(\comul\) is coassociative.  Similarly,
  \[
  \PD(1_\sour)
  = (-1)^n \PD(\comul \hot_{\sour\hot\sour} (D\hot1_\sour))
  = \PD(\comul \hot_{\sour\hot\sour} (1_\sour\hot D))
  = \Theta.
  \]
  Therefore, \(D\) is a counit for the coalgebra \((\sour,\comul)\).
\end{proof}

Now we define a natural transformation
\begin{equation}
  \label{eq:def_sigmaprime}
  \sigma'_{X,\sour}\colon \RKK^G_*(X\times Y;A,B) \to
  \RKK^G_*(Y;A\hot\sour,B\hot\sour)
\end{equation}
by \(\sigma'_{X,\sour}(f)\defeq \comul \hot_\sour \PD^{-1}(f)\),
where~\(\hot_\sour\) operates on the \emph{second} copy of~\(\sour\) in the
target \(\sour\hot\sour\) of \(\comul\).  We have
\begin{multline}
  \label{eq:sigmaprime_Theta}
  \PD(\sigma'_{X,\sour}(f))
  \defeq \Theta \hot_\sour \sigma'_{X,\sour}(f)
  = \Theta \hot_\sour \comul\hot_\sour \PD^{-1}(f)
  \\ = \Theta \hot_X \Theta \hot_\sour \PD^{-1}(f)
  = \Theta \hot_X f
\end{multline}
in \(\RKK^G_*(X\times Y;A, A'\hot\sour)\) for all \(f\in\RKK^G_{*-n}(X\times
Y;A,A')\).  It follows from the graded commutativity of exterior products and
Lemma \ref{lem:abstract_dual_coalgebra} that
\begin{equation}
  \label{eq:sigmaprime_D}
  (-1)^{ni} \sigma'_{X,\sour}(f) \hot_\sour D
  = \comul \hot_{\sour\hot\sour} D\hot \PD^{-1}(f)
  = \PD^{-1}(f)
\end{equation}
for \(f\in\RKK^G_i(X\times Y;A,B)\).  This may seem useless for computing
\(\PD^{-1}\) because the definition of \(\sigma'_{X,\sour}\) itself involves
\(\PD^{-1}\).  The point of the notion of a Kasparov dual is that we require
\(\sigma'_{X,\sour}\) to agree with another map that is easy to compute.

Recall that an \(X\cross G\)-\(C^*\)-algebra is a
\(G\)\nbd{}\(C^*\)-algebra~\(P\) equipped with a \(G\)\nbd{}equivariant
essential \(*\)\nbd{}homomorphism from \(C_0(X)\) into the centre of the
multiplier algebra of~\(P\).  This is equivalent to a \(G\)\nbd{}equivariant
essential \(*\)\nbd{}homomorphism \(m\colon C_0(X)\hot P\to P\), which we call
the \emph{\(X\)\nbd{}structure map for~\(P\)}.  Given any \(X\cross
G\)-\(C^*\)-algebra~\(P\), we get natural maps
\begin{equation}
  \label{eq:def_sigma}
  \sigma_{X,P}\colon \RKK^G_*(X;A,B) \to \KK^G_*(P\hot A,P\hot B),  
\end{equation}
which send the class of a cycle \((\Hilm,F)\) to \([(P\hot_{C_0(X)}\Hilm,
1\hot_{C_0(X)} F)]\) (see~\cite{Kasparov:Novikov}).  It is clear from the
definition that
\begin{equation}
  \label{eq:sigma_XP_pX}
  \sigma_{X,P} \bigl(p_X^*(f)\bigr) = 1_P \hot f  
\end{equation}
for all \(f\in\KK^G_*(A,B)\) and all~\(P\).

\begin{definition}
  \label{def:Kasparov_dual}
  Let~\(X\) be a locally compact \(G\)\nbd{}space.  An
  (\(n\)\nbd{}dimensional) \emph{Kasparov dual} for~\(X\) is a triple
  \((\sour,D,\Theta)\) where~\(\sour\) is a (possibly \(\Ztwo\)-graded)
  \(X\cross G\)-\(C^*\)-algebra, \(D\in\KK^G_{-n}(\sour,\C)\), and
  \(\Theta\in\RKK^G_n(X;\C,\sour)\), such that
  \begin{enumerate}[\ref{def:Kasparov_dual}.1.]
  \item \(\Theta \hot_\sour D=1_\C\) in \(\RKK^G_0(X;\C,\C)\);

  \item \(\Theta \hot_X f = \Theta \hot_\sour \sigma_{X,\sour}(f)\) in
    \(\RKK^G_*(X\times Y;A,B\hot\sour)\) for all \(f\in\RKK^G_{*-n}(X\times
    Y;A,B)\);

  \item \(\sigma_{X,\sour}(\Theta)\hot_{\sour\hot\sour}\Phi_\sour= (-1)^n
    \sigma_{X,\sour}(\Theta)\) in \(\KK^G_n(\sour,\sour\hot\sour)\).

  \end{enumerate}
\end{definition}

This definition is abstracted from the arguments in
\cite{Kasparov:Novikov}*{Section 4}.

\begin{proposition}
  \label{pro:Kasparov_dual}
  A triple \((\sour,D,\Theta)\) as above is a Kasparov dual for~\(X\) if and
  only if \((\sour,\Theta)\) is an abstract dual for~\(X\) (Definition
  \ref{def:abstract_dual}), \(\Theta \hot_\sour D=1_\C\), and the maps
  \(\sigma'_{X,\sour}\) and \(\sigma_{X,\sour}\) defined in
  \eqref{eq:def_sigmaprime} and~\eqref{eq:def_sigma} agree.

  If \((\sour,D,\Theta)\) is a Kasparov dual, then \(\comul =
  \sigma_{X,\sour}(\Theta)\), and the duality isomorphisms
  \begin{alignat*}{2}
    \PD &\colon \KK^G_{*-n}(\sour\hot A,B) \to \RKK^G_*(X;A,B),
    \\
    \PD^{-1}&\colon \RKK^G_i(X;A,B) \to \KK^G_{i-n}(\sour\hot A,B),
  \end{alignat*}
  are given by \(\PD(f)= \Theta \hot_\sour f\) and \(\PD^{-1}(f')=
  (-1)^{ni}\sigma_{X,\sour}(f')\hot_\sour D\) for \(f'\in\RKK^G_i(X;A,B)\).
\end{proposition}

\begin{proof}
  First we show that an abstract dual with the additional properties required
  in the proposition is a Kasparov dual.  Condition \ref{def:Kasparov_dual}.1
  is clear and \ref{def:Kasparov_dual}.2 follows
  from~\eqref{eq:sigmaprime_Theta} and \(\sigma'_{X,\sour}=\sigma_{X,\sour}\).
  The formula for \(\PD\) is part of the definition of an abstract dual, the
  one for \(\PD^{-1}\) follows from~\eqref{eq:sigmaprime_D}.  Since
  \(\PD^{-1}(\Theta)=1_\sour\), we have
  \(\sigma_{X,\sour}(\Theta)=\sigma'_{X,\sour}(\Theta)=\comul\).  Therefore,
  \ref{def:Kasparov_dual}.3 follows from Lemma
  \ref{lem:abstract_dual_coalgebra}.

  Suppose conversely that we have a Kasparov dual.  Define \(\PD\) and
  \(\PD^{-1}\) as in the statement of the proposition.  We must check that
  they are inverse to each other.  The composite \(\PD\circ\PD^{-1}\) sends
  \(f\in\RKK^G_i(X;A,B)\) to
  \[
  (-1)^{ni} \Theta \hot_\sour \sigma_{X,\sour}(f) \hot_\sour D
  = (-1)^{ni} \Theta \hot_X f \hot_\sour D
  = f \hot_X \Theta \hot_\sour D
  = f
  \]
  as desired.  Let \(\comul\defeq
  \sigma_{X,\sour}(\Theta)\in\RKK^G_n(\sour,\sour\hot\sour)\).  It follows
  from \ref{def:Kasparov_dual}.1 that \(\comul\hot_{\sour\hot\sour} (\ID_\sour
  \hot D)=1_\sour\).  Using \ref{def:Kasparov_dual}.3, we also get
  \(\comul\hot_{\sour\hot\sour} (D\hot \ID_\sour)=(-1)^n 1_\sour\).
  Therefore, the composite \(\PD^{-1}\circ\PD\) sends
  \(f\in\KK^G_{i-n}(\sour\hot A,B)\) to
  \begin{multline*}
    (-1)^{ni} \sigma_{X,\sour}(\Theta\hot_\sour f) \hot_\sour D
    = (-1)^{ni} \sigma_{X,\sour}(\Theta) \hot_{\sour\hot\sour} (f\hot D)
    \\ = (-1)^n \comul \hot_{\sour\hot\sour} (D\hot f)
    = (-1)^n \comul \hot_{\sour\hot\sour} (D\hot \ID_\sour) \hot_\sour f
    = f.
  \end{multline*}
  In the second step we use graded commutativity of exterior products.  Thus
  \((\sour,\Theta)\) is an abstract dual for~\(X\).  Condition
  \ref{def:Kasparov_dual}.2 asserts that \(\PD(\sigma_{X,\sour}(f))=
  \Theta\hot_X f\) for all \(f\in\RKK^G_*(X;A,B)\).  The same equation holds
  for \(\sigma'_{X,\sour}\) by~\eqref{eq:sigmaprime_Theta}.  We get
  \(\sigma'_{X,\sour}=\sigma_{X,\sour}\) because \(\PD\) is bijective.
\end{proof}

If we have a Kasparov dual for~\(X\), then we can improve the definition of
the abstract Euler characteristic:

\begin{lemma}
  \label{lem:abstract_Euler_Kd}
  Let~\(X\) be a locally compact \(G\)\nbd{}space that admits a Kasparov dual
  \((\sour,D,\Theta)\).  Let \(m\colon C_0(X)\hot\sour\to\sour\) be the
  \(X\)\nbd{}structure map for~\(\sour\) and let
  \(\bar\Theta\in\KK^G_n(C_0(X),C_0(X)\hot\sour)\) be the image
  of~\(\Theta\) under the functor that forgets the \(X\)\nbd{}structure.  Then
  \[
  \Eul_X =  \bar\Theta\hot_{C_0(X)\hot\sour} [m] \hot_\sour D.
  \]
\end{lemma}

\begin{proof}
  We have
  \[
  \Eul_X \defeq  \bar\Theta\hot_{C_0(X)\hot\sour} \PD^{-1}(\Diagonal_X)
  = \bar\Theta\hot_{C_0(X)\hot\sour} \sigma_{X,\sour}(\Diagonal_X)
  \hot_\sour D.
  \]
  It remains to check that \(\sigma_{X,\sour}(\Diagonal_X)=[m]\).  We have
  \(\Diagonal_X(f_1\hot f_2) = f_1\cdot f_2\) for all \(f_1,f_2\in C_0(X)\).
  Hence the homomorphism \(\sigma_{X,\sour}(\Diagonal_X) \colon
  C_0(X)\hot\sour\to\sour\) maps \(f_1\hot f_2\mapsto f_1\cdot f_2\) for all
  \(f_1\in C_0(X)\), \(f_2\in\sour\).  This is the definition of the
  homomorphism~\(m\).
\end{proof}

It is useful for proofs to know that Definition~\ref{def:Kasparov_dual}.2 can
be weakened as follows:

\begin{lemma}
  \label{lem:weaken_def_Kd}
  Let~\(\sour\) be an \(X\cross G\)-\(C^*\)-algebra and let
  \(\Theta\in\RKK^G_n(X;\C,\sour)\).  Suppose that the formula \(\Theta \hot_X
  f = \Theta \hot_\sour \sigma_{X,\sour}(f)\) required in Definition
  \ref{def:Kasparov_dual}.2 holds whenever~\(f\) is the class of an \((X\times
  Y)\cross G\)-linear \(*\)\nbd{}homomorphism \(C_0(X\times Y,A)\to
  C_0(X\times Y,B)\).  Then \ref{def:Kasparov_dual}.2 holds in complete
  generality.
\end{lemma}

\begin{proof}
  If \(f=p_X^*(f')\), then~\eqref{eq:sigma_XP_pX} implies \(\Theta \hot_X f =
  \Theta \hot f' = \Theta \hot_\sour \sigma_{X,\sour}(f)\).  Similarly,
  \ref{def:Kasparov_dual}.2 holds for \(p_X^*(f_0)\hot_{X,A} f\hot_{X,B}
  p_X^*(f_1)\) with \(f_0\in \KK^G_*(A_1,A)\), \(f_1\in\KK^G_*(B,B_1)\) once
  it holds for~\(f\).  Therefore, we are done if we show that every class in
  \(\RKK^G_*(X\times Y;A,B)\) can be written as a product of this kind,
  where~\(f\) is the class of an \((X\times Y)\cross G\)-linear
  \(*\)\nbd{}homomorphism.  We deduce this from considerations related to the
  universal property of \(\KK^G\).

  \cite{Meyer:KKG}*{Proposition 5.4} identifies \(\KK^G_0(A,B)\) with the set
  of homotopy classes of \(G\)\nbd{}equivariant \(*\)\nbd{}homomorphisms
  \[
  \chi(A\hot \Comp(L^2 G))\to B\hot\Comp(L^2(G\times\N)\oplus L^2(G\times
  N)^\op);
  \]
  here \(\chi(A)\) is a certain universal algebra due to Ulrich Haag.  The
  identification sends a homomorphism~\(f\) to the Kasparov product \(f_0
  \hot_{\chi(\dotso)} f\hot_{B\hot\Comp(\dotso)} f_1\) for certain natural
  elements \(f_0\in\KK^G_0(A,\chi(\dotso))\),
  \(f_1\in\KK^G_0(B\hot\Comp(\dotso),B)\).

  The same reasoning shows that the map
  \[
  f\mapsto p_X^*(f_0) \hot_{\chi(\dotso)} f\hot_{B\hot\Comp(\dotso)} p_X^*(f_1)
  \]
  identifies \(\RKK^G_0(X\times Y;A,B)\) with the set of homotopy classes of
  \((X\times Y)\cross G\)\nbd{}equivariant \(*\)\nbd{}homomorphisms
  \[
  C_0(X\times Y)\hot \chi(A\hot\Comp(L^2G))\to C_0(X\times Y)\hot
  B\hot\Comp(L^2(G\times\N)\oplus L^2(G\times N)^\op).
  \]
  Hence every element of \(\RKK^G_0(X\times Y;A,B)\) can be decomposed in the
  required form.  Putting in some more Clifford algebras, we get the same
  assertion for elements of \(\RKK^G_*(X\times Y;A,B)\).
\end{proof}

\section{The combinatorial Euler characteristic}
\label{sec:Euler_combinatorial}

Let~\(X\) be a countable, locally finite, simplicial complex, equipped with a
simplicial, continuous action of a locally compact group~\(G\).  We are going
to define the combinatorial Euler characteristic for such~\(X\).  Although we
only write down definitions for complex \(C^*\)\nbd{}algebras, it is evident
that everything we do here works in the real case as well.

We define simplicial complexes as in, say,
\cite{Brown:Buildings}*{I.Appendix}.  However, we do not consider the empty
set as a simplex.  The geometric realisation of~\(X\) is a second countable,
locally compact space because~\(X\) is locally finite and countable.  Since we
want to denote the geometric realisation by~\(X\) as well, it is convenient to
write \(SX\) for the set of (non-empty) simplices of~\(X\).  For each simplex
\(\sigma\in SX\), we write~\(\abs{\sigma}\) for the corresponding subset
of~\(X\), and we let \(\xi_\sigma\in\abs{\sigma}\) be its barycentre.  The
resulting map
\begin{equation}
  \label{eq:def_xi}
  \xi\colon SX\to X,
  \qquad \sigma\mapsto\xi_\sigma
\end{equation}
identifies \(SX\) with a discrete \(G\)\nbd{}invariant subset of~\(X\).  We
give~\(SX\) the discrete topology and the induced action of~\(G\), which is of
course continuous.  Equivalently, the stabiliser~\(G_\sigma\) of
\(\xi_\sigma\) is open for all \(\sigma\in SX\).

\emph{We require~\(G_\sigma\) to act trivially on~\(\abs{\sigma}\)}.  This is
crucial to get a correct formula for the Euler characteristic.  However, this
assumption involves no loss of generality because we may, if necessary,
replace~\(X\) by its barycentric subdivision, which clearly satisfies this
condition.

We decompose~\(SX\) into the subsets \(S_\pm X\) of simplices of even and odd
dimension.  Let \(\ell^2(S_\pm X)\) be the \(\Ztwo\)\nbd{}graded Hilbert space
with orthonormal basis \(SX\) and even and odd subspaces \(\ell^2(S_+ X)\) and
\(\ell^2(S_-X)\), respectively.  In Section~\ref{sec:dual_simplicial}, we will
mostly be using the Hilbert space \(\ell^2(SX)\) with trivial grading.  We
write \(\ell^2\bigl(S_\pm X\bigr)\) here to emphasise the non-trivial grading.
Representing \(C_0(SX)\) by diagonal operators on \(\ell^2\bigl(S_\pm
X\bigr)\), we get a natural injective \(*\)\nbd{}homomorphism \(C_0(SX)\to
\Comp\bigl(\ell^2\bigl(S_\pm X\bigr)\bigr)\), which is \(G\)\nbd{}equivariant
and by operators of even parity.  Moreover, the map~\eqref{eq:def_xi} induces
a \(G\)\nbd{}equivariant \(*\)\nbd{}homomorphism \(\xi^*\colon C_0(X)\to
C_0(SX)\).

\begin{definition}
  \label{def:Euler_combinatorial}
  Let \(\Eul_X^\comb\in\KK^G_0(C_0(X),\C)\) be the class of the
  \(G\)\nbd{}equivariant \(*\)\nbd{}homomorphism \(C_0(X)\to C_0(SX)\to
  \Comp\bigl(\ell^2\bigl(S_\pm X\bigr)\bigr)\) described above.  We call
  this the \emph{combinatorial \(G\)\nbd{}equivariant Euler characteristic
    of~\(X\)}.
\end{definition}

We now describe \(\Eul_X^\comb\) more explicitly.  For a subgroup
\({H\subseteq G}\), we let \(X^H\subseteq X\) be the fixed-point subset.  For
each connected component~\(A\) of~\(X^H\), pick a point \(x\in A\) and let
\(\dim_{H,A}\in\KK^G_0(C_0(X),\C)\) be the class of the homomorphism
\(C_0(X)\to C_0(G/H)\subseteq \Comp(\ell^2(G/H))\) that sends \(f\in C_0(X)\)
to the operator of multiplication by the function \(gH\mapsto f(gx)\).  This
does not depend on the choice of~\(x\) by homotopy invariance.  Moreover, we
have \(\dim_{gHg^{-1},gA}=\dim_{H,A}\) for all \(g\in G\) because the
resulting Kasparov cycles are unitarily equivalent.  In particular,
\(\dim_{H,gA}=\dim_{H,A}\) if~\(g\) belongs to the normaliser \(N(H)\)
of~\(H\).  Thus we may replace the connected components of~\(X^H\) by the
connected components of \(N(H)\backslash X^H\).

For an open subgroup \(H\subseteq G\) and a connected component~\(A\) of
\(N(H)\backslash X^H\), let \(S(H,A)\subseteq SX\) be the set of all simplices
of~\(A\) whose stabiliser is \emph{exactly} equal to~\(H\); notice that
each~\(A\) is a subcomplex of~\(X\).  Let
\[
\Eul^\comb_{X,H,A} \defeq [C_0(X)\overset{\xi^*}\to C_0(SX) \to
\Comp\bigl(\ell^2\bigl(G\cdot S_\pm A)\bigr)\bigr)]
\in \KK^G_0(C_0(X),\C),
\]
where the second map is the representation by diagonal operators as above.
Then
\[
\Eul_X^\comb = \sum_{(H),A} \Eul^\comb_{X,H,A},
\]
where~\((H)\) runs through the set of conjugacy classes of those open
subgroups of~\(G\) that occur as stabilisers of simplices in~\(X\), and, for
each of these, \(A\) runs through the set of connected components of
\(N(H)\backslash X^H\).

Suppose first that \(N(H)\backslash S(H,A)\) is finite.  Let
\(\chi(X,H,A)\in\Z\) be the alternating sum of the numbers of
\(n\)\nbd{}simplices in \(N(H)\backslash S(H,A)\).  Then
\(\Eul^\comb_{X,H,A}=\chi(X,H,A)\cdot \dim_{H,A}\).  If \(S(H,A)\) is
infinite, we let \(\chi(X,H,A)\defeq0\) and claim that
\(\Eul^\comb_{X,H,A}=0\) and \(\dim_{H,A}=0\).  The reason for this is that
there is a continuous path \((x_t)_{t\in\R_+}\) in~\(A\) such that
\(\lim_{t\to\infty} \norm{f|_{Gx_t}}_\infty = 0\) for all \(f\in C_0(X)\).
Thus
\begin{equation}
  \label{eq:Eul_decompose}
  \Eul_X^\comb = \sum_{(H),A} \chi(X,H,A) \cdot \dim_{H,A},
\end{equation}
where the summation runs over the same data \((H),A\) as above.  If this sum
is infinite, we have to add the cycles, not just their classes.  The summation
in~\eqref{eq:Eul_decompose} is finite if and only if all fixed-point
subspaces~\(X^H\) have finitely many connected components and there are, up to
conjugacy, only finitely many different subgroups of~\(G\) that occur as the
stabiliser of a simplex in~\(X\).

We are mainly interested in the case where~\(X\) is strongly contractible.
Then all fixed-point subsets~\(X^H\) are contractible and \emph{a fortiori}
connected.  Hence we may write \(\dim_H\label{eq:Eul_decompose_simplify}\) and
\(\chi(X,H)\) instead of \(\dim_{H,A}\) and \(\chi(X,H,A)\).

\begin{example}
  \label{exa:Eul_comb_PSL2Z}
  Consider \(G\defeq \PSL(2,\Z)\).  The free product decomposition \(G\cong
  \Z/2*\Z/3\) gives rise to a tree~\(X\) on which~\(G\) acts in such a way
  that the fundamental domain is an edge with stabilisers \(\Z/2\) and
  \(\Z/3\) at the end points and \(\{1\}\) in the interior (see
  \cite{Serre:Trees}*{\S I.4, Theorem 7}).  The action of~\(G\) on \(SX\) has
  only three orbits in this case, two orbits on vertices and one on edges.
  We find
  \[
  \Eul^\comb_X = \dim_{\Z/2} + \dim_{\Z/3} - \dim_{\{1\}}
  \in \KK^G_0(C_0(X),\C).
  \]
\end{example}

\begin{example}
  \label{exa:free_action}
  The case where~\(G\) is discrete and acts freely on~\(X\) is particularly
  simple.  Then we have natural isomorphisms
  \begin{equation}
    \label{eq:free_nonequivariant}
    \KK^G_0(C_0(X),\C) \cong \KK_0(C_0(X)\rcross G,\C)
    \cong \KK_0(C_0(G\backslash X),\C).
  \end{equation}
  They map \(\Eul_X^\comb\in\KK^G_0(C_0(X),\C)\) to \(\Eul^\comb_{G\backslash
    X}\in\KK_0(C_0(G\backslash X),\C)\).  We have
  \[
  \Eul^\comb_{G\backslash X}
  = \sum_{A\in \pi_0(G\backslash X)} \chi(A)\cdot \dim_A,
  \]
  where \(\chi(A)\) is the usual Euler characteristic of \(A\subseteq
  G\backslash X\) and \(\dim_A\) is the class in \(\KK_0(C_0( G\backslash
  X),\C)\) of the homomorphism \(C_0(G\backslash X) \to\C\), \(f\mapsto
  f(x)\), for any \(x\in A\).  If \(G\backslash X\) is connected, we get
  \(\Eul^\comb_{G\backslash X} = \chi(G\backslash X)\cdot\dim\).
\end{example}

In \cite{Luck-Rosenberg:Euler}, the relationship between various topological
constructions of Euler characteristics is discussed.  Here we consider another
construction from representation theory that is related to the Euler
characteristic of a group.

We assume that~\(G\) is a totally disconnected locally compact group for
which there is a \(G\)\nbd{}compact universal proper \(G\)\nbd{}space~\(X\).
This holds, for instance, for reductive \(p\)\nbd{}adic groups or hyperbolic
groups, where we may take the affine Bruhat-Tits building or the Rips complex,
respectively.  We may then choose~\(X\) to be a \(G\)\nbd{}finite simplicial
complex with simplicial action of~\(G\).  As in Definition
\ref{def:abstract_Euler_group}, we identify
\(\KK^G_*(C_0(X),\C)\cong\Ktop_*(G)\), so that we can view the combinatorial
Euler characteristic of~\(X\) as an element \(\Eul_X^\comb\in\Ktop_0(G)\).
This class is independent of the choice of~\(X\); we omit the verification
because our main theorem (Theorem \ref{the:Euler_combinatorial}) yields in any
case that \(\Eul_X^\comb=\Eul_\EG\), where \(\Eul_\EG\) is as in Definition
\ref{def:abstract_Euler_group}.

We assume now that~\(G\) satisfies the Baum-Connes conjecture, so that we lose
nothing by mapping \(\Eul_X\in\Ktop_0(G)\) to \(\K_0(\Cred G)\).  This class
can be described as follows:
\[
\mu_G(\Eul_X^\comb) = \sum_{\sigma \in G\backslash SX} (-1)^{\abs{\sigma}}
[\tau(G_\sigma)] \in \K_0(\Cred G),
\]
where \(\tau(G_\sigma)\in \Cred(G_\sigma)\subseteq \Cred(G)\) is the
projection onto the trivial representation of the compact-open
subgroup~\(G_\sigma\).  As a projection in the reduced \(C^*\)-algebra,
\(\tau(G)\) is given by \(\mathrm{Vol}(G_\sigma)^{-1} \cdot 1_{G_\sigma}\),
where \(1_{G_\sigma}\) is the characteristic function of~\(G_\sigma\).

The class \(\mu_G(\Eul_X^\comb) \in \K_0(\Cred G)\) is related to the
representation theory of~\(G\).  For discrete groups, this idea goes back to
Hyman Bass (\cite{Bass:Euler}).

Recall that the \emph{Hecke algebra} of~\(G\) is the space of locally
constant, compactly supported functions \(G\to\C\) with the convolution
product.  If~\(G\) is discrete, this is nothing but the group algebra
of~\(G\).  The projections \(\tau(G_\sigma)\) all lie in \(\Hecke(G)\), so
that their alternating sum actually lies in \(\K_0^\alg(\Hecke(G))\).

Let \(\REP(G)\) be the category of smooth representations of~\(G\) (always on
complex vector spaces).  We say that a smooth representation of~\(G\) has
\emph{type (FP)} if it has a finite length resolution \((P_n,\delta_n)\) by
finitely generated projective objects of \(\REP(G)\).  This resolution is
unique up to chain homotopy equivalence.  We have \(P_n\cong
\Hecke(G)^{k_n}\cdot p_n\) for certain projections \(p_n\in
M_{k_n}(\Hecke(G))\), which yield classes in \(\K_0^\alg(\Hecke(G))\).  Let
\[
\chi(M)\defeq \sum_{n=0}^\infty (-1)^n [p_n] \in \K_0^\alg(\Hecke(G)).
\]
This is well-defined, that is, independent of the choices of the resolution
\((P_n,\delta_n)\) and the projections~\(p_n\).

We may consider the cellular chain complex of~\(X\) with coefficients~\(\C\)
as a chain complex of smooth representations of~\(G\).  Its homology vanishes
for \(*>0\) and is~\(\C\) with the trivial representation of~\(G\) for
\(*=0\).  Since~\(X\) is \(G\)\nbd{}compact, \(\C[SX]\) is a \emph{finite}
direct sum of \(\Hecke(G)\)\brd{}modules of the form \(\C[G/H]\) for certain
compact-open subgroups \(H\subseteq G\).  The latter are finitely generated
projective objects of \(\REP(G)\) because \(\C[G/H]\cong \Hecke(G)\cdot
\tau_H\) for all compact-open subgroups \(H\subseteq G\).  Thus the trivial
representation of~\(G\) has type (FP) and the natural map
\(\K_0^\alg(\Hecke(G))\to\K_0(\Cred G)\) maps \(\chi(\C)\) to
\(\mu_G(\Eul_\EG)\).

Representation theorists usually replace \(\chi(\C)\) by its Chern character,
which belongs to \(\HH_0(\Hecke(G))\defeq \Hecke(G)/[\Hecke(G),\Hecke(G)]\);
it is represented by the function
\[
\sum_{\sigma \in G\backslash SX} (-1)^{\abs{\sigma}} \tau(G_\sigma) \in
\Hecke(G).
\]
If~\(G\) is a semi-simple \(p\)\nbd{}adic group, another representative for
the same class is the \emph{Euler-Poincaré function} of Robert Kottwitz
(\cite{Kottwitz}*{Section 2}); it is computed from the cellular chain complex
of the affine Bruhat-Tits building with its natural poly-simplicial structure.
Recall that we refine this to a simplicial structure with the additional
property that~\(G_\sigma\) fixes the simplex~\(\sigma\) pointwise.  Both chain
complexes produce the same class in \(\HH_0(\Hecke(G))\), even in
\(\K_0^\alg(\Hecke(G))\), because they both come from finite projective
resolutions of the trivial representation of~\(G\).

Peter Schneider and Ulrich Stuhler construct analogous Euler-Poincaré
functions for general irreducible representations of semi-simple
\(p\)\nbd{}adic groups in \cite{Schneider-Stuhler}; it is not hard to lift
these Euler-Poincaré functions to elements of \(\K_0^\alg(\Hecke(G))\),
see~\cite{Meyer:Reductive_isocoh}.  Although the Borel-Serre compactification
plays an important role in~\cite{Schneider-Stuhler}, it is not clear to us
whether these more general Euler characteristics of irreducible
representations are related to Kasparov duals or boundary actions.

\section{A Kasparov dual for simplicial complexes}
\label{sec:dual_simplicial}

In this section, we assume~\(X\) to be a \emph{finite-dimensional}, locally
finite, countable simplicial complex equipped with a simplicial, continuous
action of a locally compact group~\(G\).  We do not require the action to be
proper.

Our main goal is to exhibit the combinatorial Euler characteristic as an
abstract Euler characteristic.  A Kasparov dual for~\(X\) has already been
constructed by Gennadi Kasparov and Georges Skandalis in
\cite{Kasparov-Skandalis:Buildings}.  However, they only describe~\(\Theta\)
indirectly, which makes it hard to compute \(\Eul_X\).  Therefore, we give an
independent and completely explicit construction for~\(\Theta\).  We also
modify their definitions of \(\sour\) and~\(D\) slightly to get a simple
formula for~\(\Theta\).

We need some preparations before we can start the actual construction.  As in
Section~\ref{sec:Euler_combinatorial}, \(SX\) denotes the set of simplices
of~\(X\), and we usually write~\(X\) both for the simplicial complex and its
geometric realisation.  Let \(S_0X\subseteq SX\) be the set of vertices, that
is, \(0\)\nbd{}simplices of~\(X\).  Suppose that~\(X\) is at most
\(n\)\nbd{}dimensional and let
\[
\no\defeq\{0,1,\dotsc,n\}.
\]
A \emph{colouring} on~\(X\) is a function \(\nu\colon S_0X\to\no\) such that
for any simplex \(\sigma\in SX\), the images under~\(\nu\) of the vertices
of~\(\sigma\) are pairwise distinct.  A \emph{coloured simplicial complex} is
a simplicial complex equipped with such a colouring.  The action of~\(G\) is
compatible with the colouring if the function~\(\nu\) is \(G\)\nbd{}invariant.
(Coloured simplicial complexes are called \emph{typed} in
\cite{Kasparov-Skandalis:Buildings}.)  Most of our constructions only involve
a single simplex in~\(X\) at a time.  The colouring allows us to piece these
local constructions together.

Let~\(X\) be any \(n\)\nbd{}dimensional simplicial complex and let~\(X^{(1)}\)
be its barycentric subdivision.  Recall that the vertex set of~\(X^{(1)}\) is
equal to the set of simplices of~\(X\); the simplices in~\(X^{(1)}\) are
labelled bijectively by strictly increasing chains in the partially ordered
set of simplices of~\(X\); here the partial order is defined by
\(\sigma\le\sigma'\) if~\(\sigma\) is a face of~\(\sigma'\).  The map
\(S_0X^{(1)}=SX\to\no\) that sends a simplex to its dimension is a canonical
colouring on~\(X^{(1)}\).  Thus it is no loss of generality to assume~\(X\) to
carry a \(G\)\nbd{}invariant colouring; we assume this in the following.

We shall use the affine Euclidean space
\begin{equation}
  \label{eq:def_E}
  E \defeq \Bigl\{ (t_0,\dotsc,t_n)\in\R^{n+1} \Bigm| \sum t_i = 1\Bigr\}.  
\end{equation}
Sometimes, we specify a point in~\(E\) by \emph{homogeneous coordinates}:
\begin{equation}
  \label{eq:homogeneous_coordinates}
  [t_0,\dotsc,t_n] \defeq \bigl(\sum t_i\bigr)^{-1} (t_0,\dotsc,t_n)  
\end{equation}
provided \(\sum t_i\neq0\).  We realise the standard \(n\)\nbd{}simplex as the
subset
\begin{equation}
  \label{eq:def_Sigma}
  \Sigma \defeq \{(t_0,\dotsc,t_n) \in E \mid \text{\(t_i \ge 0\) for all
    \(i\in\no\)}\}.
\end{equation}

Let \(\poset\) be the partially ordered set of \emph{non-empty} subsets
of~\(\no\).  We extend the colouring~\(\nu\) to a map \(SX\to\poset\) by
sending a simplex to the set of colours of its vertices.  We also define
\(\nu(\emptyset)\defeq\emptyset\).  We identify~\(\no\) with the set of
vertices of~\(\Sigma\).  Since a face of~\(\Sigma\) is determined by the set
of vertices it contains, this identifies \(\poset\) with the partially ordered
set of faces of~\(\Sigma\).  Under this identification, \(f\subseteq\no\)
corresponds to the face
\begin{multline}
  \label{eq:def_face}
  \abs{f} \defeq \{(t_0,\dotsc,t_n)\in\Sigma\mid
  \text{\(t_i=0\) for \(i\in\no\setminus f\)}\}
  \\ = \{(t_0,\dotsc,t_n)\in E\mid
  \text{\(t_i\ge0\) for \(i\in f\) and \(t_i=0\) for \(i\in\no\setminus f\)}\}.
\end{multline}

We may view the map \(\nu\colon SX\to\poset\) as a \(G\)\nbd{}invariant
simplicial map; passing to geometric realisations, we get a
\(G\)\nbd{}invariant continuous map
\begin{equation}
  \label{eq:colouring_map}
  \abs{\nu}\colon X\to\Sigma.
\end{equation}
Any point \(x\in X\) belongs to some simplex \(\sigma\in SX\).  The
restriction of~\(\abs{\nu}\) to~\(\abs{\sigma}\) is the unique affine map that
sends a vertex of colour~\(i\) to the corresponding vertex of~\(\Sigma\).
If~\(\sigma\) is of dimension~\(k\), then \(\nu(\sigma)\subseteq\no\) has
\(k+1\) elements and hence defines a \(k\)\nbd{}dimensional face
of~\(\Sigma\).  Hence the restriction of~\(\abs{\nu}\) to~\(\abs{\sigma}\) is
a homeomorphism from~\(\abs{\sigma}\) to the face \(\abs{\nu(\sigma)}\)
of~\(\Sigma\), which we denote by~\(\abs{\nu}_\sigma\).

\begin{figure}
  \includegraphics{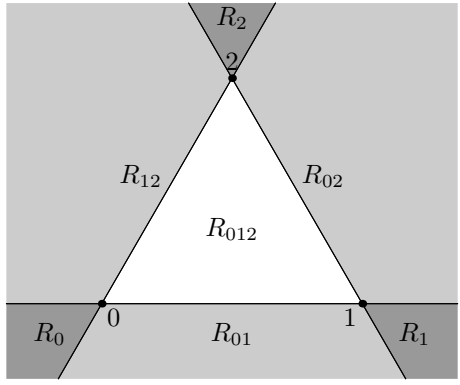}

  \caption{The regions \(R_f\)}
  \label{fig:regions}
\end{figure}

For any \(f\subseteq\no\), we define a closed convex subset \(R_f\subseteq E\)
by
\begin{equation}
  \label{eq:def_Rf}
  R_f\defeq \{(t_0,\dotsc,t_n)\in E\mid
  \text{\(t_i\ge0\) for \(i\in f\) and \(t_i\le0\) for \(i\in\no\setminus
    f\)}\}.
\end{equation}
These regions are a crucial ingredient of our construction.
Figure~\ref{fig:regions} illustrates them for \(n=2\).
We have \(R_\no=\Sigma\) and \(R_\emptyset=\emptyset\) because no point
of~\(E\) satisfies \(t_i<0\) for all \(i\in\no\).  The sets~\(R_f\) for
\(f\in\poset\) cover~\(E\) and have mutually disjoint interiors.  We also
define \(R_S\defeq \bigcup_{f\in S} R_f\) if \(S\subseteq\poset\) is a set of
faces of~\(\Sigma\).  We are mainly interested in
\begin{equation}
  \label{eq:def_Rlef}
  R_{\le f} \defeq \bigcup_{\{l\in\poset\mid l\le f\}} R_l
  = \{(t_0,\dotsc,t_n)\in E\mid \text{\(t_i\le0\) for \(i\in\no\setminus
    f\)}\}
\end{equation}
for \(f\subseteq\no\).  It follows immediately that
\begin{equation}
  \label{eq:intersect_Rlef}
  R_{\le f_1}\cap R_{\le f_2} = R_{\le (f_1\cap f_2)}
\end{equation}
for all \(f_1,f_2\in\poset\).

We define a retraction \(q\colon E \to\Sigma\) by
\begin{equation}
  \label{eq:def_qf}
  q(t_0,\dotsc,t_n)\defeq [\max(t_0,0),\dotsc,\max(t_n,0)].
\end{equation}
Inspection of \eqref{eq:def_face} and~\eqref{eq:def_Rlef} yields
\begin{equation}
  \label{eq:qinv_f}
  q^{-1}(\abs{f}) = R_{\le f}.
\end{equation}

Kasparov and Skandalis use the nearest point retraction to~\(\Sigma\) instead
of~\(q\) in~\cite{Kasparov-Skandalis:Buildings}.  We prefer~\(q\) because of
the more explicit formula.

We can now define the underlying \(C^*\)\nbd{}algebra~\(\sour\) of our
Kasparov dual.  We use the \(C^*\)\nbd{}algebra of compact operators on
\(\ell^2(SX)\).  The group~\(G\) acts on this Hilbert space via its action on
the basis~\(SX\).  We equip \(\ell^2(SX)\) with the trivial grading, as
opposed to the grading by parity that we used in
Section~\ref{sec:Euler_combinatorial}.  We describe an operator~\(T\) on
\(\ell^2(SX)\) by a matrix \((T_{\sigma\sigma'})_{\sigma,\sigma'\in SX}\).
For a function \(\varphi\colon Y\to \Bound(\ell^2(SX))\), its matrix
coefficients are also functions \(\varphi_{\sigma\sigma'}\colon Y\to\C\),
defined by \(\varphi_{\sigma\sigma'}(y)\defeq \varphi(y)_{\sigma\sigma'}\) for
\(y\in Y\).  Let
\begin{multline}
  \label{eq:def_sour_combinatorial}
  \sour \defeq \{\varphi\in C_0(E,\Comp(\ell^2(SX))) \mid
  \\ \text{\(\supp \varphi_{\sigma\sigma'}\subseteq R_{\le \nu(\sigma\cap\sigma')}\)
    for all \(\sigma,\sigma'\in SX\)}\}.
\end{multline}
Hence \(\varphi_{\sigma\sigma'}=0\) unless \(\sigma\) and~\(\sigma'\) have a
common face.  We let~\(G\) act on \(C_0(E,\Comp(\ell^2(SX)))\) by
\(g\varphi(t)\defeq \pi_g\circ\varphi(t)\circ\pi_g^{-1}\) for all \(g\in G\),
\(t\in E\).  Obviously, \(\sour\) is a closed, self-adjoint,
\(G\)\nbd{}invariant subspace of \(C_0(E,\Comp(\ell^2(SX)))\).  We have to
check that~\(\sour\) is closed under multiplication.  If
\(\varphi,\psi\in\sour\), \(\sigma,\sigma'\in SX\), then we have
\((\varphi\cdot\psi)_{\sigma\sigma'} = \sum_{\tau\in SX} \varphi_{\sigma\tau}
\psi_{\tau\sigma'}\).  Using equation~\eqref{eq:intersect_Rlef} and that the
colouring is injective on the vertices of~\(\tau\), we get
\[
\supp \varphi_{\sigma\tau}\psi_{\tau\sigma'}
\subseteq R_{\le\nu(\sigma\cap\tau)} \cap R_{\le\nu(\sigma'\cap\tau)}
= R_{\le\nu(\sigma\cap\tau\cap\sigma')}
\subseteq R_{\le\nu(\sigma\cap\sigma')}.
\]
Hence each individual summand \(\varphi_{\sigma\tau} \psi_{\tau\sigma'}\)
satisfies the support condition~\eqref{eq:def_sour_combinatorial}.
Thus~\(\sour\) is a \(G\)\nbd{}invariant \(C^*\)\nbd{}subalgebra of
\(C_0(E,\Comp(\ell^2(SX)))\).

We may interpret the algebra~\(\sour\) physically as follows.  The
simplices~\(\sigma\) are possible states of a system.  For \(t\in E\), the
system may only be in the state~\(\sigma\) if \(t_i<0\) for all \(i\in
\no\setminus \nu(\sigma)\); two such states \(\sigma,\sigma'\) may interact if
\(t_i<0\) for all \(i\in\no\setminus \nu(\sigma\cap\sigma')\).

Next we define the \(X\)\brd{}structure map \(m\colon
C_0(X)\otimes\sour\to\sour\).  Recall that the map~\(\abs{\nu}\) defined
in~\eqref{eq:colouring_map} restricts to a homeomorphism
\(\abs{\nu}_\sigma\colon \abs{\sigma}\to\abs{\nu(\sigma)}\) for each
\(\sigma\in SX\) and that \(q(R_{\le f})\subseteq \abs{f}\)
by~\eqref{eq:qinv_f}.  Hence we may define a continuous \(G\)\nbd{}equivariant
map
\begin{equation}
  \label{eq:def_barq}
  E\times SX \supseteq \bigcup_{\sigma\in SX}
  R_{\le\nu(\sigma)}\times\{\sigma\}
  \overset{\bar{q}}\longrightarrow X,
  \qquad \bar{q}(t,\sigma) \defeq \abs{\nu}_\sigma^{-1}\bigl(q(t)\bigr).
\end{equation}
It is convenient to extend this function to all of \(E\times SX\).  We may do
this by \(\bar{q}(t,\sigma)\defeq \abs{\nu}_\sigma^{-1}\circ
a_{\nu(\sigma)}\circ q(t)\), where we choose simplicial retractions
\(a_f\colon \Sigma\to\abs{f}\) for \(f\in\poset\).  With this extension
of~\(\bar{q}\) we get a \(G\)\nbd{}equivariant essential
\(*\)\nbd{}homomorphism
\[
m'\colon \xymatrix{
  C_0(X) \ar[r]^-{\bar{q}^*} & C_b(E\times SX) \ar[r] &
  C_0(E,\Comp(\ell^2(SX))),
}
\]
where the second map is the representation by diagonal operators on
\(\ell^2(SX)\).  By construction, \(\bar{q}(t,\sigma)\in\abs{\sigma}\) for all
\(t\in E\), \(\sigma\in SX\).

\begin{lemma}
  \label{lem:mprime_central}
  We have \(m'(\varphi_1)\circ\varphi_2=\varphi_2\circ m'(\varphi_1)\in\sour\)
  for all \(\varphi_1\in C_0(X)\), \(\varphi_2\in\sour\).
\end{lemma}

\begin{proof}
  It follows from the definitions that
  \begin{align*}
    (m'(\varphi_1)\circ\varphi_2)_{\sigma\sigma'}(t)
    &= \varphi_1\bigl(\bar{q}(t,\sigma)\bigr) \cdot
    \varphi_2(t)_{\sigma\sigma'},
    \\
    (\varphi_2\circ m'(\varphi_1))_{\sigma\sigma'}(t)
    &= \varphi_1\bigl(\bar{q}(t,\sigma')\bigr) \cdot
    \varphi_2(t)_{\sigma\sigma'}
  \end{align*}
  for all \(\sigma,\sigma'\in SX\), \(t\in E\).  The function
  \((\varphi_2)_{\sigma\sigma'}\) is supported in the region
  \(R_{\le\nu(\sigma\cap\sigma')}\) because \(\varphi_2\in\sour\),
  see~\eqref{eq:def_sour_combinatorial}.  Therefore,
  \((m'(\varphi_1)\varphi_2)_{\sigma\sigma'}\) and \((\varphi_2
  m'(\varphi_1))_{\sigma\sigma'}\) are supported in this region as well, so
  that \(m'(\varphi_1) \varphi_2\in\sour\) and \(\varphi_2
  m'(\varphi_1)\in\sour\).  It remains to check that the two matrix
  coefficients above agree.  This follows if
  \(\bar{q}(t,\sigma)=\bar{q}(t,\sigma')\) for all \(t\in
  R_{\le\nu(\sigma\cap\sigma')}\) because both matrix coefficients are
  supported in this region.

  If \(t\in R_{\le\nu(\sigma\cap\sigma')}\), then \(t\in R_{\le\nu(\sigma)}\)
  and \(t\in R_{\le\nu(\sigma')}\) by~\eqref{eq:intersect_Rlef}.  Thus
  \(\bar{q}(t,\sigma)\) and \(\bar{q}(t,\sigma')\) are both defined
  by~\eqref{eq:def_barq}.  We have \(q(t)\in\abs{\nu(\sigma\cap\sigma')}\)
  by~\eqref{eq:qinv_f}.  Both \(\abs{\nu}_\sigma\) and~\(\abs{\nu}_{\sigma'}\)
  extend \(\abs{\nu}_{\sigma\cap\sigma'}\), which is a homeomorphism onto the
  face \(\abs{\nu(\sigma\cap\sigma')}\).  Therefore,
  \(\bar{q}(t,\sigma)=\bar{q}(t,\sigma\cap\sigma')=\bar{q}(t,\sigma')\) as
  desired.
\end{proof}

Hence there is a unique \(*\)\nbd{}homomorphism \(m\colon
C_0(X)\otimes\sour\to\sour\) with \(m(\varphi_1\otimes\varphi_2)\defeq
m'(\varphi_1)\circ\varphi_2\); since \(C_0(X)\) is nuclear, it does not matter
which tensor product we choose here.  The map~\(m\) is \(G\)\nbd{}equivariant
and essential, so that~\(\sour\) becomes an \(X\cross G\)-\(C^*\)-algebra.  If
we view \(\varphi\in C_0(X)\otimes\sour\) as a function \(X\times
E\to\Comp(\ell^2(SX))\), we can describe~\(m\) explicitly in terms of matrix
coefficients:
\begin{equation}
  \label{eq:def_m}
  m(\varphi)_{\sigma\sigma'}(t)
  = \varphi_{\sigma\sigma'}(\bar{q}(t,\sigma\cap\sigma'),t)
  = \varphi_{\sigma\sigma'}(\abs{\nu}_{\sigma\cap\sigma'}^{-1}\circ q(t), t).
\end{equation}
The last expression has to be taken with a grain of salt because
\(\abs{\nu}_{\sigma\cap\sigma'}^{-1}\circ q(t)\) is only defined for \(t\in
R_{\le\nu(\sigma\cap\sigma')}\); for other values of~\(t\), we have
\(\varphi_{\sigma\sigma'}(x,t)=0\) regardless of the value of~\(x\) because of
the definition of~\(\sour\) in~\eqref{eq:def_sour_combinatorial}.

Next, we define \(D\in\KK^G_n(\sour,\C)\).  Let \([i]\in\KK^G_0\bigl(\sour,
C_0(E)\bigr)\) be the class of the inclusion map \(i\colon \sour \to
C_0(E)\otimes\Comp(\ell^2(SX))\).  Since \(E\cong\R^n\), we have canonical
invertible elements
\[
\beta_E\in\KK_n(C_0(E),\C),
\qquad
\hat{\beta}_E\in\KK_{-n}\bigl(\C,C_0(E)\bigr),
\]
such that
\begin{alignat*}{2}
  \beta_E \otimes \hat{\beta}_E &= 1_{C_0(E)}
  &\qquad& \text{in \(\KK^G_0(C_0(E),C_0(E))\),}
  \\
  \hat{\beta}_E \otimes_{C_0(E)} \beta_E &= 1_\C
  &\qquad& \text{in \(\KK^G_0(\C,\C)\).}
\end{alignat*}
We set
\begin{equation}
  \label{eq:def_D_comb}
  D \defeq [i]\otimes_{C_0(E)} \beta_E\in\KK^G_n(\sour,\C).
\end{equation}

We will construct \(\Theta\in\RKK^G_{-n}(X;\C,\sour)\) as \(\Theta\defeq
\hat{\beta}_E\otimes_{C_0(E)} [\vartheta]\), where \([\vartheta]\in\RKK^G_0(X;
C_0(E),\sour)\) is the class of an \(X\cross G\)\nbd{}equivariant
\(*\)\nbd{}homomorphism \(\vartheta\colon C_0(X)\otimes C_0(E)\to
C_0(X)\otimes \sour\).  The latter is, of course, equivalent to a
\(G\)\nbd{}equivariant continuous family of \(*\)\nbd{}homomorphisms
\(\vartheta_x\colon C_0(E)\to\sour\) for \(x\in X\).  Its construction is
rather involved.  This is the point where we deviate most seriously
from~\cite{Kasparov-Skandalis:Buildings}.

\begin{figure}
  \includegraphics{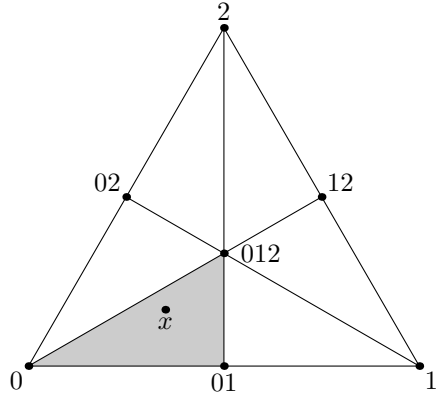}\centering

  \caption{barycentric subdivision of the standard \(2\)-simplex}
  \label{fig:barycentric}
\end{figure}

The first ingredient for~\(\vartheta\) is a certain \(G\)\nbd{}equivariant
function from~\(X\) to the unit sphere of~\(\ell^2(SX)\).  For this we need
the barycentric subdivision~\(X^{(1)}\) of~\(X\).  Recall that the vertices of
this subdivision are in bijection with \(SX\).  Let \(x\in X\) and
let~\(\sigma^{(1)}\) be some simplex of the barycentric subdivision that
contains~\(x\).  The vertices of~\(\sigma^{(1)}\) form a strictly increasing
chain \(\sigma_0\subset\dotsb\subset\sigma_k\) in \(SX\); we view these as
basis vectors of \(\ell^2(SX)\).  Any point of~\(\abs{\sigma^{(1)}}\) can be
written uniquely as a convex combination of the vertices~\(\sigma_j\);
formally, \(x=\sum_{j=0}^k t_j\sigma_j\), where \(t_j\ge0\) for all
\(j\in\nok\) and \(\sum t_j=1\).  The barycentric subdivision of a
\(2\)\nbd{}simplex is illustrated in Figure~\ref{fig:barycentric}; we have
represented the vertices by their colours in \(\{0,1,2\}\).  The shaded
maximal simplex of the barycentric subdivision is labelled by the chain \(0,
01, 012\).

Let
\begin{equation}
  \label{eq:def_vprime}
  v'(x) = v'\Bigl(\sum_{j=0}^k t_j\sigma_j\Bigr)
  \defeq \sum_{j=0}^k \sqrt{t_j} \sigma_j \in \ell^2(SX).
\end{equation}
This defines a continuous map from~\(\abs{\sigma^{(1)}}\) to the unit sphere
of \(\ell^2(SX)\).  If some of the coordinates~\(t_j\) of~\(x\) vanish, then
we may replace~\(\sigma^{(1)}\) by the face that is spanned by
the~\(\sigma_j\) with \(t_j\neq0\).  Since this does not change \(v'(x)\), the
maps~\(v'\) may be glued together to a continuous map \(v'\colon
X\to\ell^2(SX)\), whose range is contained in the unit sphere.  Now let
\(P'\colon X\to\Comp(\ell^2(SX))\) be the function whose value at \(x\in X\)
is the rank-1-projection onto \(\C\cdot v'(x)\).  The maps \(v'\) and~\(P'\)
are evidently \(G\)\nbd{}equivariant.

An important point about this definition is that the basis vectors involved
in~\(v'(x)\) form a chain in~\(SX\); hence there is some region in~\(E\) where
\(P'(x)\) is a possible value for an element of~\(\sour\).  Observe
that~\(P'(x)\) is a diagonal operator in the basis \(SX\) if and only
if~\(v'(x)\) is a basis vector, if and only if~\(x\) is a vertex of the
barycentric subdivision; equivalently, \(x=\xi_\sigma\) for some \(\sigma\in
SX\).  In order to proceed with the construction, we need a projection valued
function that is diagonal not merely \emph{at} these points but \emph{near}
them.  Therefore, we replace \(v'\) and~\(P'\) by \(v\defeq v'\circ\collapse\)
and \(P\defeq P'\circ\collapse\) with a certain \emph{collapsing map}
\(\collapse\colon X\to X\).

Choose \(L\in (0,1/(n+1))\).  We first define a map \(\collapse_\Sigma\colon
\Sigma\to\Sigma\) on the standard simplex by
\begin{equation}
  \label{eq:def_collapse}
  \collapse_\Sigma\bigl((t_0,\dotsc,t_n)\bigr)
  \defeq [\min\{t_0,L\},\dotsc,\min\{t_n,L\}],
\end{equation}
where~\([\dotso]\) denotes homogeneous coordinates as
in~\eqref{eq:homogeneous_coordinates}.  If \(t_j=0\), then \(\min\{t_j,L\}=0\)
as well.  This means that \(\collapse_\Sigma(\abs{f})\subseteq\abs{f}\) for
each face~\(f\) of~\(\Sigma\) (these faces are defined
in~\eqref{eq:def_face}).  Therefore, if \(\sigma\in SX\) we may define
\[
\collapse\colon \abs{\sigma}\to\abs{\sigma},
\qquad x\mapsto
\abs{\nu}_\sigma^{-1}\circ\collapse_\Sigma\circ\abs{\nu}_\sigma(x),
\]
using the homeomorphisms~\(\abs{\nu}_\sigma\) defined
after~\eqref{eq:colouring_map}.  These maps on simplices match on
\(\abs{\sigma\cap\sigma'}\), so that we get a global map \(\collapse\colon
X\to X\).  Now we put \(v\defeq v'\circ\collapse\) and \(P\defeq
P'\circ\collapse\); thus~\(P(x)\) is the rank-1-projection onto \(\C\cdot
v(x)\) for all \(x\in X\).

In order to formulate some properties of the collapsing map, we define
\begin{equation}
  \label{eq:def_CRf}
  CR_f \defeq \{(t_0,\dotsc,t_n)\in \Sigma \mid
  \text{\(t_i\ge L\) for \(i\in f\), \(t_i\le L\) for \(i\in \no\setminus
    f\)}\}.
\end{equation}%
\begin{figure}
  \includegraphics{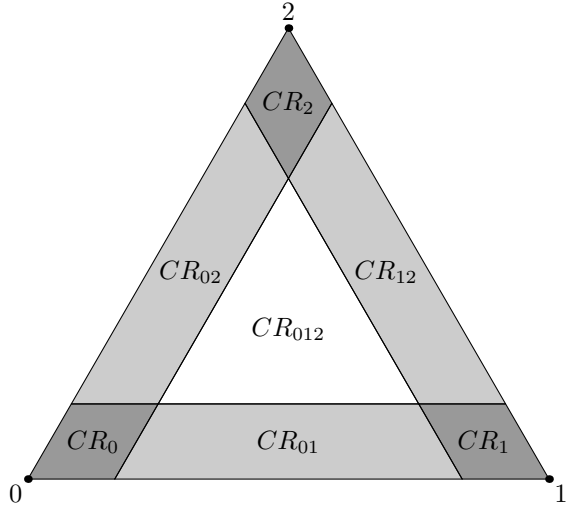}
  \centering
  \caption{The regions \(CR_f\)}
  \label{fig:CR_f}
\end{figure}%
for \(f\in\poset\).  These regions cover~\(\Sigma\) because \(L<1/(n+1)\).
Figure~\ref{fig:CR_f} illustrates them for the \(2\)\nbd{}simplex.  Combining
\eqref{eq:def_face} and~\eqref{eq:def_CRf}, we get
\begin{equation}
  \label{eq:fCRf}
  \abs{f}\cap CR_f = \{(t_0,\dotsc,t_n)\in \Sigma \mid
  \text{\(t_i\ge L\) for \(i\in f\), \(t_i=0\) for \(i\in \no\setminus f\)}\}.
\end{equation}
Hence \(\collapse_\Sigma(\abs{f}\cap CR_f)\) consists of a single point
\((t_i')\), with homogeneous coordinates \(t_i'=L\) for \(i\in f\) and
\(t_i'=0\) for \(i\in\no\setminus f\).  The rescaling replaces~\(L\) by \(1/\#
f\) and thus produces the barycentre~\(\xi_f\) of the face~\(f\); that is,
\begin{equation}
  \label{eq:collapse_fCRf}
  \collapse_\Sigma(\abs{f}\cap CR_f) = \{\xi_f\}.
\end{equation}

\begin{lemma}
  \label{lem:SRf}
  If \(x\in X\) satisfies \(\abs{\nu}(x)\in CR_f\), then
  \(P(x)_{\sigma\sigma'}=0\) or \(f\subseteq\nu(\sigma\cap\sigma')\).
\end{lemma}

\begin{proof}
  Let~\(\tau\) be some simplex of~\(X\) that contains~\(x\).  Since \(v(x)\)
  only has non-zero coefficients at the faces of~\(\tau\) and since the
  restriction of~\(\abs{\nu}\) to~\(\abs{\tau}\) is injective, we may assume
  without loss of generality that \(X=\Sigma\) and \(\abs{\nu}=\ID_\Sigma\).
  Moreover, the assertion is invariant under simplicial automorphisms
  of~\(\Sigma\), that is, permutations of coordinates.  (We transform
  both~\(x\) and~\(f\), of course.)  We can achieve that \(x_0\ge
  x_1\ge\dotsb\ge x_n\) by a coordinate permutation.  Hence the only
  possibilities for~\(f\) are \(f=\{0,\dotsc,k\}\) for some \(k\in\no\).  Let
  \(x'\defeq \collapse_\Sigma(x)\), then we get \(P(x)=P'(x')\) and
  \[
  x'_0=x'_1=\dotsb=x'_k\ge x'_{k+1} \ge\dotsb\ge x'_n.
  \]

  Let \(\alpha_j\in\Sigma\) be the vertex of the barycentric subdivision
  of~\(\Sigma\) that is labelled by the simplex \(\{0,\dotsc,j\}\);
  equivalently, \(\alpha_j\) is the barycentre of that simplex; explicitly,
  the first \(j+1\) coordinates of~\(\alpha_j\) are \(1/(j+1)\), the remaining
  ones vanish.  It is straightforward to see that~\(x'\) is a convex
  combination of~\(\alpha_j\) with \(j\ge k\).  Such convex combinations form
  a single simplex in the barycentric subdivision.  Hence the vector
  \(v(x)=v'(x')\in\ell^2(SX)\) only contains the basis vectors~\(\alpha_j\)
  with \(j\ge k\).  Therefore, if \(P(x)_{\sigma,\sigma'}\neq0\), then
  \(\sigma\) and~\(\sigma'\) are among the~\(\alpha_j\) with \(j\ge k\).  This
  implies \(f\subseteq \nu(\sigma\cap\sigma')\) as asserted.
\end{proof}

For \(\lambda>1\), let \(r_\lambda\colon E \to E\) be the radial expansion map
around the barycentre of~\(\Sigma\).  Explicitly,
\begin{equation}
  \label{eq:radial_expansion}
  r_\lambda(t_0,\dotsc,t_n) = \biggl(\lambda t_0 - \frac{\lambda-1}{n+1}, 
  \dotsc, \lambda t_n - \frac{\lambda-1}{n+1}\biggr).
\end{equation}
If \((\lambda-1)/(n+1)=\lambda L\), that is, \(\lambda = (1-(n+1)L)^{-1}\),
then we get
\begin{equation}
  \label{eq:rlambdainv_Rf}
  CR_f = r_\lambda^{-1}(R_f)
\end{equation}
and hence \(r_\lambda(CR_f)\subseteq R_f\subseteq R_{\le f}\); this follows
immediately from the definitions \eqref{eq:def_Rf} and~\eqref{eq:def_CRf}, see
also Figures \ref{fig:regions} and~\ref{fig:CR_f}.  We shall need a slightly
different result, as follows.  For \(\delta>0\), let
\[
B(\delta) \defeq
\bigl\{ (t_0,\dotsc,t_n)\in\R^{n+1} \bigm|
\text{\(\sum t_j=0\) and \(\abs{t_j}<\delta\) for all \(j=0,\dotsc,n\)} \bigr\}.
\]

\begin{lemma}
  \label{lem:expand_CRf}
  If \(\lambda> (1-(n+1)L)^{-1}\), then there exists \(\delta>0\) such that
  \(r_\lambda(s)+B(\delta) \subseteq R_{\le f}\) for all \(f\in\poset\),
  \(s\in CR_f\).
\end{lemma}

\begin{proof}
  This follows immediately from the definition of~\(r_\lambda\) and the
  definitions \eqref{eq:def_Rlef} and~\eqref{eq:def_CRf} of the regions
  \(R_{\le f}\) and \(CR_f\).
\end{proof}

Choose \(\lambda\) and \(\delta>0\) as in Lemma~\ref{lem:expand_CRf} and
choose an orientation-preserving diffeomorphism \(h\colon E\congto
B(\delta)\).  Let~\(E_+\) be the one-point compactification of~\(E\).
Extend~\(h^{-1}\) to a map \(h^{-1}\colon E_+\to E_+\) by
\(h^{-1}(t)\defeq\infty\) for \(t\notin B(\delta)\) and extend \(\varphi\in
C_0(E)\) to~\(E_+\) by \(\varphi(\infty)\defeq0\).  We get a continuous family
of \(*\)\nbd{}homomorphisms
\begin{equation}
  \label{eq:def_hs}
  h_s!\colon C_0(E)\to C_0(E),
  \qquad h_s! \varphi(t) \defeq \varphi\circ
  h^{-1}\bigl(t-r_\lambda(s)\bigr)
\end{equation}
for \(s\in\Sigma\), where \(r_\lambda\colon E\to E\) is defined
in~\eqref{eq:radial_expansion}.  Our notation stems from the fact
that~\(h_s!\) is the wrong-way map associated to the open embedding
\[
h_s\colon E\to E,\qquad t\mapsto h(t)+r_\lambda(s).
\]
By construction, \(h_s!(\varphi)\) vanishes outside \(r_\lambda(s)+B(\delta)\)
for all \(\varphi\in C_0(E)\).  Using the map \(\abs{\nu}\colon X\to\Sigma\)
defined in \eqref{eq:colouring_map}, we get a \(G\)\nbd{}invariant continuous
family of \(*\)\nbd{}homomorphisms \(h_{\abs{\nu}(x)}\colon C_0(E)\to C_0(E)\)
parametrised by \(x\in X\).

\begin{lemma}
  \label{lem:def_rho}
  The formula \(\vartheta_x(\varphi)\defeq h_{\abs{\nu}(x)}!(\varphi) \otimes
  P(x)\) for \(x\in X\) defines a \(G\)\nbd{}equivariant continuous family of
  \(*\)\nbd{}homomorphisms \(C_0(E)\to\sour\) and hence a class
  \([\vartheta]\in\RKK^G_0(X; C_0(E),\sour)\).  We define
  \[
  \Theta \defeq \hat{\beta}_E \otimes_{C_0(E)} [\vartheta] \in
  \RKK^G_{-n}(X;\C,\sour).
  \]
\end{lemma}

\begin{proof}
  It is clear that~\(\vartheta_x\) is a \(G\)\nbd{}equivariant continuous
  family of \(*\)\nbd{}homomorphisms into
  \(C_0(E,\Comp(\ell^2(SX)))\supseteq\sour\).  We must check that its range is
  contained in~\(\sour\).  Fix \(\sigma,\sigma'\in SX\) and \(x\in X\) such
  that \(P(x)_{\sigma\sigma'}\neq0\).  We have \(\abs{\nu}(x)\in CR_f\) for
  some \(f\in\poset\) because these regions cover~\(\Sigma\).
  Lemma~\ref{lem:SRf} yields \(f\subseteq \nu(\sigma\cap\sigma')\).  Let
  \(V\defeq r_\lambda\abs{\nu}(x)+B(\delta)\), then \(V\subseteq R_{\le f}
  \subseteq R_{\le\nu(\sigma\cap\sigma')}\) by Lemma~\ref{lem:expand_CRf}.
  Since \(h_{\abs{\nu}(x)}!\varphi\) is supported in~\(V\) for all
  \(\varphi\in C_0(E)\), we get \(\supp \vartheta_x(\varphi)_{\sigma\sigma'}
  \subseteq R_{\le\nu(\sigma\cap\sigma')}\).  This means that
  \(\vartheta_x(\varphi)\in\sour\) (see \eqref{eq:def_sour_combinatorial}).
\end{proof}

\begin{theorem}
  \label{the:Kd_combinatorial}
  The triple \((\sour,D,\Theta)\) defined above is Kasparov dual for~\(X\) of
  dimension~\(-n\).
\end{theorem}

\begin{proof}
  First we check condition \ref{def:Kasparov_dual}.1, that is, \(\Theta
  \otimes_\sour D=1_\C\).  Let \(\vartheta= (\vartheta_x)_{x\in X}\) be as in
  Lemma \ref{lem:def_rho}.  Let~\(i\) be the embedding \(\sour\to
  C_0(E)\otimes\Comp(\ell^2(SX))\); it defines a class
  \([i]\in\KK^G_0(\sour,C_0(E))\).  Then
  \[
  \Theta\otimes_\sour D
  = \hat\beta_E \otimes_{C_0(E)} [\vartheta]\otimes_\sour
  [i]\otimes_{C_0(E)}\beta_E.
  \]
  We are done if we show \([i\circ\vartheta]=1_{C_0(E)}\) in
  \(\RKK^G_0(X;C_0(E),C_0(E))\) because \(\hat\beta_E\otimes_{C_0(E)}
  \beta_E=1_\C\).  Since we no longer impose any support restrictions on the
  range of \(i\circ\vartheta_x\), the family of maps
  \[
  \vartheta_x^s(\varphi)(t) \defeq
  \varphi\circ h^{-1}\bigl(t-r_{s\lambda}\bigl(\abs{\nu}(x)\bigr)\bigr) P(x)
  \]
  for \(s\in[0,1]\) provides a natural homotopy between \(i\circ
  \vartheta_x=\vartheta_x^1\) and the map
  \[
  \vartheta_x^0\varphi(t)\defeq \varphi\circ h^{-1}(t-\xi) P(x),
  \]
  where~\(\xi\) is the barycentre of~\(\Sigma\).  Since~\(h\) is an
  orientation-preserving homeomorphism \(E\to B(\delta)\), the endomorphism
  \(\varphi\mapsto \varphi\circ h^{-1}(t-\xi)\) of \(C_0(E)\) is homotopic to
  the identity map.  Thus~\([\vartheta]\) is the exterior product of
  \(1_{C_0(E)}\) and the class \([P]\in\RKK^G_0(X;\C,\C)\) determined by the
  continuous family of projections \(P(x)\), \(x\in X\).  The continuous
  family of unit vectors~\(v(x)\) may be viewed as a \(G\)\nbd{}equivariant
  continuous family of isometries \(\hat{v}(x)\colon \C\to\ell^2(SX)\) with
  \(\hat{v}(x)\hat{v}^*(x)=P(x)\).  This means that \([P]=[1_\C]\).  This
  finishes the proof that \(\Theta \otimes_\sour D=1_\C\).

  Next we verify \ref{def:Kasparov_dual}.2, which asserts that \(\Theta\hot_X
  f = \Theta \hot_\sour \sigma_{X,\sour}(f)\) in \(\RKK^G_*(X\times
  Y;A,B\hot\sour)\) for all \(f\in\RKK^G_{*+n}(X\times Y;A,B)\) and all
  \(Y,A,B\).  Since the classes \(\beta_E\) and~\(\hat\beta_E\) are inverse to
  each other, this is equivalent to \([\vartheta]\hot_X f = [\vartheta]
  \hot_\sour \sigma_{X,\sour}(f)\) in \(\RKK^G_*(X\times Y;A\hot
  C_0(E),B\hot\sour)\).  By Lemma \ref{lem:weaken_def_Kd}, it suffices to
  prove this in the special case where~\(f\) is an \((X\times Y)\cross
  G\)-equivariant \(*\)\nbd{}homomorphism.  Thus both factors in our product
  are now classes of equivariant \(*\)\nbd{}homomorphisms.

  We view \(f\) as a \(G\)\nbd{}equivariant continuous family of
  \(C_0(Y)\)\nbd{}linear \(*\)\nbd{}homomorphisms \(f_x\colon C_0(Y,A)\to
  C_0(Y,B)\) for \(x\in X\).  Then \(\vartheta\hot_X f\) corresponds to the
  continuous family of maps \(\vartheta_x\hot f_x\colon C_0(E)\hot C_0(Y,A)\to
  \sour\hot C_0(Y,B)\), \(x\in X\).  Explicitly,
  \begin{multline}
    \label{eq:f_X_rho}
    (\vartheta\hot_X f)_x(\varphi\hot a)_{\sigma\sigma'}(t)
    \defeq \vartheta_x(\varphi)_{\sigma\sigma'}(t)  f_x(a)
    \\ = (h_{\abs{\nu}(x)}!\varphi)(t) P(x)_{\sigma\sigma'} f_x(a)
    \qquad\text{in \(C_0(Y,B)\)}
  \end{multline}
  for all \(\varphi\in C_0(E)\), \(a\in C_0(Y,A)\), \(t\in E\),
  \(\sigma,\sigma'\in SX\).

  By definition, we have \(\sigma_{X,\sour}(f)(\varphi_1\cdot \varphi_2\hot a)
  = \varphi_1\cdot f(\varphi_2\hot a)\) for all \(\varphi_1\in\sour\),
  \(\varphi_2\in C_0(X)\), \(a\in C_0(Y,A)\).  Using \eqref{eq:def_m}, we
  rewrite this as
  \begin{equation}
    \label{eq:sigma_Xsour_explicit}
    \sigma_{X,\sour}(f)(\varphi\hot a)_{\sigma\sigma'}(t)
    = \varphi_{\sigma\sigma'}(t) f_{\bar{q}(\sigma,t)}(a)
    \qquad\text{in \(C_0(Y,B)\)}
  \end{equation}
  for all \(\varphi\in\sour\), \(a\in C_0(Y,A)\), \(\sigma,\sigma'\in SX\),
  \(t\in E\).  Composition with~\(\vartheta\) yields the continuous family of
  maps \((\vartheta\hot_\sour \sigma_{X,\sour}(f))_x =
  \sigma_{X,\sour}(f)\circ(1_A\hot\vartheta_x)\) from \(C_0(E)\hot C_0(Y,A)\)
  to \(\sour\hot C_0(Y,B)\) for \(x\in X\).  Thus
  \begin{equation}
    \label{eq:rho_sour_sigmaf}
    (\vartheta\hot_\sour \sigma_{X,\sour}(f))_{x}(\varphi\hot a)_{\sigma\sigma'}(t)
    = (h_{\abs{\nu}(x)}!\varphi)(t) P(x)_{\sigma\sigma'} f_{\bar{q}(\sigma,t)}(a).
  \end{equation}
  The only difference between the two families in \eqref{eq:f_X_rho}
  and~\eqref{eq:rho_sour_sigmaf} is that we use \(f_x\) and
  \(f_{\bar{q}(\sigma,t)}\), respectively.  Whenever
  \(P(x)_{\sigma\sigma'}\neq0\), \(x\) and \(\bar{q}(\sigma,t)\) lie in the
  same simplex \(\sigma\cap\sigma'\) of~\(X\).  Arguments as in the proofs of
  Lemmas \ref{lem:mprime_central} and \ref{lem:def_rho} show that we still get
  a homomorphism from \(C_0(E)\hot C_0(Y,A)\) to \(\sour\hot C_0(Y,B)\) if we
  replace~\(f_{\bar{q}(\sigma,t)}\) in \eqref{eq:rho_sour_sigmaf} with
  \(f_{(1-s)x+s\bar{q}(\sigma,t)}\) for \(s\in[0,1]\).  Thus \(\vartheta\hot_X
  f\) and \(\vartheta\hot_\sour \sigma_{X,\sour}(f)\) are homotopic.  This
  finishes the proof of \ref{def:Kasparov_dual}.2.

  It remains to verify \ref{def:Kasparov_dual}.3.  Since~\(\beta_E\) is
  invertible, we may replace \(\sigma_{X,\sour}(\Theta)\) with
  \(\sigma_{X,\sour}([\vartheta])\) in this statement.  This is the class of a
  \(G\)\nbd{}equivariant \(*\)\nbd{}homomorphism
  \(\sigma_{X,\sour}(\vartheta)\colon \sour\otimes C_0(E)\to
  \sour\otimes\sour\).  We must check
  \[
  [\Phi_\sour\circ\sigma_{X,\sour}(\vartheta)] =
  (-1)^n[\sigma_{X,\sour}(\vartheta)]
  \qquad\text{in \(\KK^G_0(\sour\otimes C_0(E), \sour\otimes\sour)\);}
  \]
  here~\(\Phi_\sour\) denotes the flip automorphism on \(\sour\otimes\sour\).
  We describe \(\sigma_{X,\sour}(\vartheta)\) by specifying its matrix
  coefficients with respect to the basis \(SX\times SX\) of
  \(\ell^2(SX)\hot\ell^2(SX)\cong\ell^2(SX\times SX)\).  Equation
  \eqref{eq:sigma_Xsour_explicit} yields
  \begin{multline}
    \label{eq:technical_comul}
    \sigma_{X,\sour}(\vartheta)(\varphi_1\hot\varphi_2)
    (t_1,t_2)_{(\sigma_1,\sigma_2),(\sigma_1',\sigma_2')}
    \\ = \varphi_1(t_1)_{\sigma_1\sigma_1'} \cdot \varphi_2\circ
    h^{-1}(t_2-r_\lambda\circ \abs{\nu}\circ\bar{q}(t_1,\sigma_1)) \cdot
    P(\bar{q}(\sigma_1,t_1))_{\sigma_2\sigma_2'}
  \end{multline}
  for all \(\varphi_1\in\sour\), \(\varphi_2\in C_0(E)\), \(t_1,t_2\in E\),
  \(\sigma_1,\sigma_1',\sigma_2,\sigma_2'\in SX\).  Fix \(t_1,t_2\in E\) and
  choose \(f\subseteq\no\) minimal such that~\(t_1\) belongs to the interior
  of~\(R_{\le f}\).  Thus \(q(t_1)\in\abs{f}\) by \eqref{eq:qinv_f}.  If
  \(\varphi_1(t_1)_{\sigma_1\sigma_1'}\neq0\), then
  \(f\subseteq\nu(\sigma_1\cap\sigma_1')\) by the definition of~\(\sour\),
  see~\eqref{eq:def_sour_combinatorial}.  Hence \eqref{eq:def_barq} yields
  \(\bar{q}(\sigma_1,t)=\abs{\nu}^{-1}_{\sigma_1} q(t_1)\).  Thus
  \(\abs{\nu}\bar{q}(\sigma_1,t_1)=q(t_1)\) and we can rewrite the right hand
  side of \eqref{eq:technical_comul} as
  \[
  \varphi_1(t_1)_{\sigma_1\sigma_1'} \cdot \varphi_2\circ
  h^{-1}(t_2-r_\lambda q(t_1)) \cdot
  P(\bar{q}(\sigma_1,t_1))_{\sigma_2\sigma_2'}.
  \]

  For \(\sigma\in SX\), \(f\subseteq\no\), we let
  \begin{equation}
    \label{eq:colour_partition}
    \begin{aligned}
      SX_{\ge\sigma} &\defeq \{\sigma'\in SX\mid \sigma'\ge\sigma\},
      \\
      SX_{\ge f} &\defeq \{\sigma'\in SX\mid \nu(\sigma')\supseteq f\},
      \\
      SX_{=f} &\defeq \{\sigma'\in SX\mid \nu(\sigma')=f\},
    \end{aligned}
  \end{equation}
  Since~\(\nu\) is a colouring, any simplex in \(SX_{\ge f}\) contains a
  unique face~\(\sigma\) with \(\nu(\sigma)=f\).  This means that~\(SX_{\ge
    f}\) is the disjoint union of the subsets \(SX_{\ge\sigma}\), where
  \(\sigma\in SX_{=f}\).  We write \(\sigma_{1,f}\defeq\tau\) if
  \(\nu(\tau)=f\) and \(\sigma_1\in SX_{\ge\tau}\).  With~\(f\) as defined
  above, we have \(\bar{q}(\sigma_1,t_1)\in \abs{\sigma_{1,f}}\).  Thus
  \(P(\bar{q}(\sigma_1,t_1))_{\sigma_2\sigma_2'}=0\) unless \(\sigma_2\)
  and~\(\sigma_2'\) are faces of~\(\sigma_{1,f}\).  We also choose~\(f_2\)
  such that \(q(t_1)\in CR_{f_2}\).  Then it follows from Lemma \ref{lem:SRf}
  that \(f_2\subseteq \nu(\sigma_2\cap\sigma_2')\) for all
  \(\sigma_2,\sigma_2'\in SX\) with
  \(P(\bar{q}(\sigma_1,t_1))_{\sigma_2\sigma_2'}\neq0\).  Since both
  \(\sigma_2\) and~\(\sigma_2'\) are faces of~\(\sigma_1\), this is equivalent
  to \(\sigma_2,\sigma_2'\in SX_{\ge \sigma_{1,f_2}}\).  Moreover, Lemma
  \ref{lem:expand_CRf} yields that \(t_2\) belongs to the interior of~\(R_{\le
    f_2}\).

  It follows from the definition of~\(\sour\) that the possible values of
  \(\varphi_1(t_1)\) for \(\varphi_1\in\sour\) are exactly the elements of
  \[
  \bigoplus_{\sigma\in SX_{=f}} \Comp(\ell^2(SX_{\ge\sigma}))
  \subseteq \Comp(\ell^2(SX)).
  \]
  A similar description is available for the possible values at~\(t_2\), of
  course; the relevant face \(f'\subseteq\no\) is the minimal subset for
  which~\(t_2\) is an interior point of \(R_{\le f'}\).  We have \(f'\le f_2\)
  by Lemma \ref{lem:expand_CRf}.  The map \(\sigma_{X,\sour}(\vartheta)\)
  gives rise to an embedding
  \[
  \bigoplus_{\tau\in SX_{=f}} \Comp(\ell^2(SX_{\ge\tau}))
  \to
  \bigoplus_{\tau\in SX_{=f}}
  \Comp(\ell^2(SX_{\ge\tau}\times SX_{\ge\tau_{f'}})),
  \]
  which is induced by the family of isometries
  \[
  J(\tau,t_1)\colon \ell^2(SX_{\ge\tau})\to
  \ell^2(SX_{\ge\tau}\times SX_{\ge\tau_{f'}}),
  \qquad \eta\mapsto \eta\otimes v(\bar{q}(\tau,t_1)),
  \]
  for \(\tau\in SX_{=f}\).  Here we use the definition of~\(P(x)\) as the
  rank-1-projection onto the span of the unit vector \(v(x)\).

  Since the coefficients of \(v(x)\) are non-negative for all \(x\in X\), we
  have \((1-s)v(x)+s\sigma\neq0\) for any \(s\in [0,1]\), \(\sigma\in SX\) (we
  view~\(\sigma\) as a basis vector of \(\ell^2(SX)\)).  Therefore, we may
  deform the isometry \(J(\tau,t_1)\) by a continuous path of isometries
  \(J^s(\tau,t_1)(\sigma)\colon \ell^2(SX_{\ge\tau})\to
  \ell^2(SX_{\ge\tau}\times SX_{\ge\tau_{f'}})\), defined by
  \[
  J^s(\tau,t_1)(\sigma) \defeq \sigma\otimes
  \frac{(1-s) v(\bar{q}(\tau,t_1)) + s\sigma}{\norm{(1-s)
      v(\bar{q}(\tau,t_1)) + s\sigma}}
  \]
  for \(\sigma\in SX_{\ge\tau}\).  These isometries yield a homotopy of
  \(*\)\nbd{}homomorphisms
  \[
  \Ad J^s(t_1)\colon \bigoplus_{\tau\in SX_{=f}} \Comp(\ell^2(SX_{\ge\tau}))
  \to
  \bigoplus_{\tau\in SX_{=f}} \Comp(\ell^2(SX_{\ge\tau}\times
  SX_{\ge\tau_{f'}})).
  \]
  Letting \(t_1,t_2\) vary again, we get a homotopy of \(G\)\nbd{}equivariant
  \(*\)\nbd{}homomorphisms \(\sour\otimes C_0(E)\to \sour\otimes\sour\) by
  sending \(\varphi_1\otimes \varphi_2\) to the function
  \[
  (t_1,t_2)\mapsto \Ad J^s(t_1) \varphi_1(t_1) \cdot
  \varphi_2\circ h^{-1}(t_2 - r_\lambda q(t_1)).
  \]
  Thus \(\sigma_{X,\sour}(\vartheta)\) is homotopic to the
  \(G\)\nbd{}equivariant \(*\)\nbd{}homomorphism \(\alpha\colon \sour\otimes
  C_0(E)\to \sour\otimes\sour\) defined by
  \[
  \alpha(\varphi_1\otimes\varphi_2)(t_1,t_2)
  \defeq \Ad J^1 \varphi_1(t_1) \cdot
  \varphi_2\circ h^{-1}(t_2 - r_\lambda q(t_1)),
  \]
  where \(J^1\colon \ell^2(SX)\to\ell^2(SX\times SX)\) is the diagonal
  embedding that sends the basis vector \(\sigma\in SX\) to
  \(\sigma\otimes\sigma\).

  Fix~\(t_1\) once again and let \(f\subseteq\no\) be as above.  Then \(\Ad
  J^1 \varphi(t_1)\) is an allowed value for functions in
  \(\sour\otimes\sour\) if~\(t_2\) belongs to the interior of \(R_{\ge f}\).
  Since this holds for the points in \(r_\lambda((1-s)q(t_1)+st_1)+B(\delta)\)
  for \(s\in[0,1]\), the map~\(\alpha\) is homotopic to
  \[
  \alpha'(\varphi_1\otimes\varphi_2)(t_1,t_2)
  \defeq \Ad J^1 \varphi_1(t_1) \cdot \varphi_2\circ h^{-1}(t_2 - r_\lambda t_1).
  \]
  Observe that \(\varphi\mapsto \varphi\circ r_{\lambda^{-s}}\) for \(s\ge0\)
  defines an endomorphism of~\(\sour\).  Hence the homotopy
  \(\alpha'_s(\varphi)(t_1,t_2) \defeq
  \alpha'(\varphi)(r_{\lambda^{-s}}t_1,t_2)\) connects~\(\alpha'\) with
  \[
  \alpha''(\varphi)(t_1,t_2)
  \defeq \Ad J^1 \varphi\bigl(r_{\lambda^{-1}}t_1, h^{-1}(t_2 - t_1)\bigr).
  \]
  Since \([\sigma_{X,\sour}(\vartheta)]=[\alpha'']\), condition
  \ref{def:Kasparov_dual}.3 is equivalent to
  \([\Phi_\sour\circ\alpha'']=(-1)^n[\alpha'']\).  We have
  \[
  \Phi_\sour\circ \alpha''(\varphi)(t_1,t_2)
  = \Ad J^1 \varphi\bigl(r_{\lambda^{-1}}t_2, h^{-1}(t_1 - t_2)\bigr)
  \]
  because the range of~\(J^1\) is invariant under
  \(\Phi_{\Comp(\ell^2(SX))}\).  We define yet another homotopy of
  homomorphisms \(\sour\otimes C_0(E)\to \sour\otimes\sour\) by
  \[
  \alpha''_s(\varphi)(t_1,t_2)
  \defeq \Ad J^1 \varphi\bigl(r_{\lambda^{-1}}(st_1+(1-s)t_2),
  h^{-1}(t_1 - t_2)\bigr).
  \]
  It connects \(\Phi_\sour\circ\alpha''\) and \(\alpha''\circ
  (\ID_\sour\otimes f)\), where \(f\colon C_0(E)\to C_0(E)\) is induced by the
  map \(t\mapsto -t\) on~\(E\); here we assume that~\(h\) is an even function,
  as we may.  Of course, \([f]=(-1)^n\) in \(\KK_0(C_0(E),C_0(E))\).  This
  finishes the proof of \ref{def:Kasparov_dual}.3.
\end{proof}

\begin{theorem}
  \label{the:Euler_combinatorial}
  Let~\(X\) be a simplicial complex equipped with a simplicial action
  of~\(G\).  Then \(\Eul_X=\Eul_X^\comb\) in \(\KK^G_0(C_0(X),\C)\).
\end{theorem}

\begin{proof}
  Lemma \ref{lem:abstract_Euler_Kd} asserts that
  \[
  \Eul_X = \bar\Theta \hot_{C_0(X)\hot\sour} [m] \hot_{\sour} D
  \in \KK^G_0(C_0(X),\C);
  \]
  here \(\bar\Theta\in\KK^G_{-n}(C_0(X),C_0(X)\hot\sour)\) is obtained from
  the class~\(\Theta\) defined in Lemma \ref{lem:def_rho} by forgetting the
  \(X\)\nbd{}structure; \(m\colon C_0(X)\hot\sour\to\sour\) is the
  multiplication homomorphism, which is described in~\eqref{eq:def_m};
  and~\(D\) is defined in~\eqref{eq:def_D_comb}.  Since the Bott periodicity
  classes \(\beta_E\) and~\(\hat\beta_E\), which appear in \(\Theta\)
  and~\(D\), are inverse to each other, our assertion is equivalent to
  \[
  \ID_{C_0(E)}\hot \Eul_X^\comb
  = \beta_E \hot \Eul_X \hot \hat\beta_E
  = [\vartheta] \hot_{C_0(X)\hot\sour} [m] \hot_{\sour} [i],
  \]
  where \(\vartheta\colon C_0(X)\hot C_0(E)\to C_0(X)\hot\sour\) is the
  continuous family of \(*\)\nbd{}homomorphisms defined in Lemma
  \ref{lem:def_rho} and \([i]\in\KK^G_0(\sour, C_0(E))\) is the class of the
  inclusion map \(i\colon \sour\to C_0(E)\hot\Comp(\ell^2(SX))\).
  The above Kasparov product is the class of the composite homomorphism
  \[
  i\circ m\circ\vartheta\colon C_0(X\times E)\to C_0(E,\Comp(\ell^2(SX))).
  \]

  Plugging in the definition of~\(\vartheta\) and~\eqref{eq:def_m}, we get
  \[
  i\circ m\circ \vartheta(\varphi)_{\sigma\sigma'}(t)
  = \varphi(\bar{q}(t,\sigma),h^{-1}(t-r_\lambda q(t))) \cdot
  P(\bar{q}(t,\sigma))_{\sigma\sigma'}.
  \]
  for all \(\sigma,\sigma'\in SX\), \(t\in E\), \(\varphi\in C_0(X\times E)\).
  We want to simplify this expression.  Assume that it is non-zero.  Let
  \(t'\defeq q(t)\in\Sigma\) and \(x'\defeq \bar{q}(t,\sigma)\).  Since the
  regions \(CR_f\) cover~\(\Sigma\), we have \(t'\in CR_f\) for some
  \(f\in\poset\).  Lemma \ref{lem:expand_CRf} yields \(r_\lambda
  t'+B(\delta)\subseteq R_{\le f}\); hence we must have \(t\in R_{\le f}\) in
  order for \(\varphi(\bar{q}(t,\sigma),h^{-1}(t-r_\lambda q(t)))\) to be
  non-zero.  Since \(t\in R_{\le f}\), we get \(q(t)\in\abs{f}\cap CR_f\) by
  \eqref{eq:qinv_f} and \(\collapse_\Sigma(q(t))=\xi_f\)
  by~\eqref{eq:collapse_fCRf}.  Therefore,
  \(P(\bar{q}(t,\sigma))=P'\circ\collapse(\bar{q}(t,\sigma))\) is the
  projection onto a basis vector of \(\ell^2(SX)\).  Hence
  \(P(\bar{q}(t,\sigma))_{\sigma\sigma'}=0\) unless \(\sigma=\sigma'\) and
  \(\nu(\sigma)=f\).

  Summing up, \(i\circ m\circ \vartheta(\varphi)(t)\in\Comp(\ell^2(SX))\) is
  diagonal with respect to the basis \(SX\); the diagonal entry for the basis
  vector~\(\sigma\) is supported in
  \[
  D_{\nu(\sigma)} \defeq q^{-1}\bigl(CR_{\nu(\sigma)}\cap\abs{\nu(\sigma)}\bigr)
  \]
  and given there by the formula
  \[
  \Lambda_\sigma(\varphi)(t)
  \defeq \varphi(\bar{q}(t,\sigma),h^{-1}(t-r_\lambda q(t)))
  = \varphi(\abs{\nu}_\sigma^{-1} q(t),h^{-1}(t-r_\lambda q(t))).
  \]
  The last term is defined for \(t\in D_{\nu(\sigma)}\) because
  \(q(t)\in\abs{\nu(\sigma)}\) (see the definition of~\(\bar{q}\)
  in~\eqref{eq:def_barq}).  Let \(\Lambda\defeq (\Lambda_\sigma)\colon
  C_0(X\times E)\to C_0(SX\times E)\).

  In the definition of \(\Eul_X^\comb\) we use the map \(\xi^*\colon C_0(X)\to
  C_0(SX)\) induced by \eqref{eq:def_xi}.  Moving
  \(\abs{\nu}_\sigma^{-1}q(t)\in\abs{\sigma}\) linearly
  towards~\(\xi_\sigma\), we get a homotopy between~\(\Lambda\) and the
  \(*\)\nbd{}homomorphism \(\Lambda'\circ (\xi^*\otimes \ID_{C_0(E)})\), where
  we define \(\Lambda'\colon C_0(SX\times E)\to C_0(SX\times E)\) by
  \begin{equation}
    \label{eq:Lambdaprime}
    \Lambda'\varphi(\sigma,t) \defeq
    \begin{cases}
      \varphi(\sigma,h^{-1}(t-r_\lambda q(t))) & \text{for \(t\in
        D_{\nu(\sigma)}\),}
      \\
      0 & \text{otherwise.}
    \end{cases}
  \end{equation}

  We may describe \(\Lambda'\) by a family of maps \(\Lambda'_\sigma\colon
  C_0(E)\to C_0(E)\) for \(\sigma\in SX\).  Equation~\eqref{eq:Lambdaprime}
  shows that \(\Lambda'_\sigma\) only depends on \(\nu(\sigma)\), so that we
  also denote it by \(\Lambda'_{\nu(\sigma)}\).  Partitioning \(SX\) into the
  subsets \(SX_{=f}\) defined in~\eqref{eq:colour_partition}, we obtain
  \[
  [i\circ m\circ\vartheta] = \sum_{f\in\poset}
  [\xi^* \hot \Lambda'_f]\hot_{C_0(SX)} [C_0(SX)\to\Comp(\ell^2(SX_{=f}))]
  \]
  in \(\KK^G_0(C_0(X)\hot C_0(E),C_0(E))\).  The combinatorial Euler
  characteristic is defined by
  \[
  \Eul_X^\comb \defeq \sum_{f\in\poset} (-1)^{\dim f} [\xi^*]\hot_{C_0(SX)}
  [C_0(SX)\to\Comp(\ell^2(SX_{=f}))].
  \]
  Therefore, \([i\circ m\circ\vartheta]=\ID_{C_0(E)}\hot\Eul_X^\comb\) follows
  if \([\Lambda'_f]=(-1)^{\dim f}\) in \(\KK_0(C_0(E),C_0(E))\).  It remains
  to verify this assertion.

  Since all our constructions are invariant under coordinate permutations, we
  may assume \(f=\{0,\dotsc,k\}\) with \(k=\dim f\).  If \(t\in D_f\), then
  \(t_i>0\) for \(i\in f\) and \(t_i<0\) for \(i\in \no\setminus f\).
  Hence~\(q\) is given by
  \[
  q(t)_i =
  \begin{cases}
    \bigl(1-\sum_{j\in\no\setminus f} t_j\bigr)^{-1}t_i
    & \text{for \(i\in f\),} \\
    0 & \text{for \(i\in \no\setminus f\).}
  \end{cases}
  \]
  The point \(r_\lambda\xi_f\) belongs to~\(D_f\) and satisfies \(q
  r_\lambda(\xi_f)=\xi_f\).  Hence it is a fixed-point of the map
  \(r_\lambda\circ q\).  We reparametrise our maps and consider
  \(\psi_f(t)\defeq r_\lambda \xi_f + t - r_\lambda q(r_\lambda \xi_f+t)\),
  where \(t\in D_f-r_\lambda\xi_f\).  Thus \(\psi_f(0)=0\).  Since points in
  the range of~\(\psi_f\) satisfy \(\sum \psi_f(t)_i=0\), we may drop one
  coordinate; we choose the \(0\)th coordinate, which belongs to~\(f\).

  It is easy to see that \(\psi_f(t)_i=t_i\) for \(i\in\no\setminus f\).
  Moreover, if we fix the coordinates \(t_j\) with \(j\in\no\setminus f\),
  then \(\psi_f(t)_i = -a\bigl(\sum_{j\in\no\setminus f} t_j\bigr)\cdot t_i +
  b\bigl(\sum_{j\in\no\setminus f} t_j\bigr)\) with certain rational functions
  \(a,b\) of one variable.  Explicitly,
  \[
  a(s) = \frac{(k+1)(\lambda-1) + (n+1) s}{(n-k)\lambda + k+1- (n+1) s}.
  \]
  The important point here is that \(a(0)>0\).

  Since \(h^{-1}\circ\psi_f(t)=\infty\) unless \(\abs{\psi_f(t)_i}<\delta\)
  for all \(i\in\no\), we may restrict attention to~\(t\) with
  \(\abs{t_j}<\delta\) for \(j\in\no\setminus f\), so that
  \(\abs{s}<(n-k)\delta\).  We may choose~\(\delta\) as small and~\(\lambda\)
  as great as we like.  Therefore, the difference between~\(a\) and the
  constant function \(a(0)\) is negligible.  Hence the maps \(a(rs)\cdot
  t_i+b(rs)\) for \(r\in[0,1]\) give rise to an isotopy between \(\psi_f\) and
  the invertible linear map
  \[
  \psi_f'(t)_i =
  \begin{cases}
    t_i & \text{for \(i\in\no\setminus f\),}
    \\
    - a(0)\cdot t_i & \text{for \(i\in f\), \(i\neq0\).}
  \end{cases}
  \]
  Recall that we have dropped one coordinate, so that we do not have to worry
  about the condition \(\sum t_i=0\) any more.  This also means that there
  only remain \(\dim f\) relevant coordinates in~\(f\), which are multiplied
  by a negative number.  Hence \([\psi_f']=(-1)^{\dim f}\) in
  \(\KK_0(C_0(\R^n),C_0(\R^n))\).  Since~\(h\) is orientation-preserving, we
  have \([h^{-1}]=1\).  Therefore, \([\Lambda'_f]=
  [h^{-1}\circ\psi_f]=[h^{-1}\circ\psi'_f]=(-1)^{\dim f}\).
\end{proof}

\section{Gysin sequence in the simplicial case}
\label{sec:Gysin_combinatorial}

\begin{theorem}
  \label{the:Gysin_combinatorial}
  Let~\(G\) be a locally compact group and let \(\bd{X}=\cl{X}\setminus X\) be
  a boundary action as in Definition \ref{def:boundary_action}.  Suppose
  that~\(X\) is a finite-dimensional, locally finite simplicial complex with a
  simplicial action of~\(G\).  Suppose also that~\(G\) satisfies the
  Baum-Connes conjecture with coefficients \(\C\) and \(C(\bd{X})\).  Then
  there is an exact sequence
  \[\xymatrix{
    0 \ar[r] &
    \K_1(\Cred G) \ar[r]^-{u_*} &
    \K_1(C(\bd{X})\rcross G) \ar[r]^{\delta} &
    \K_0(C_0(X)\rcross G) \ar[d]^-{\Eul_X^\comb} \\
    0 &
    \K_1(C_0(X)\rcross G) \ar[l] &
    \K_0(C(\bd{X})\rcross G) \ar[l]_-{\delta} &
    \K_0(\Cred G) \ar[l]_-{u_*}
  }
  \]
  where \(\Eul_X^\comb\) denotes the Kasparov product with the combinatorial
  equivariant Euler characteristic \(\Eul_X^\comb \in \KK^G_0(C_0(X),\C)\).
  More explicitly,
  \[
  \Eul_X^\comb(x) = \sum_{(H)} \chi(X,H) \dim_H(x),
  \]
  where the summation runs over the conjugacy classes in~\(G\) of stabilisers
  of simplices in~\(X\), and \(\dim_H\in\KK^G_0(C_0(X),\C)\) and
  \(\chi(X,H)\in\Z\) are defined as on
  page~\pageref{eq:Eul_decompose_simplify}.
\end{theorem}

\begin{proof}
  We plug the formula \(\Eul_X=\Eul_X^\comb\) of
  Theorem~\ref{the:Euler_combinatorial} into the abstract Gysin sequence of
  Proposition \ref{pro:abstract_Gysin}.  By definition, \(\Eul_X^\comb\)
  factors through the homomorphism \(\xi^*\colon C_0(X)\to C_0(SX)\) induced
  by the barycentre embedding~\eqref{eq:def_xi}.  Writing~\(SX\) as a disjoint
  union of \(G\)\nbd{}orbits, we get
  \begin{multline*}
    \Ktop_*(G,C_0(SX))
    \cong \K_*(C_0(SX)\rcross G)
    \\ \cong \bigoplus_{\sigma\in G\backslash SX} \K_*(C_0(G/G_\sigma)\rcross G)
    \cong \bigoplus_{\sigma\in G\backslash SX} \K_*(\Cred(G_\sigma)).
  \end{multline*}
  Here~\(G_\sigma\) denotes the stabiliser of the simplex~\(\sigma\), which is
  a compact-open subgroup of~\(G\).  The map \(\Eul_X^\comb\colon
  \K_1(C_0(X)\rcross G)\to \K_1(\Cred G)\) vanishes because it factors through
  \(\K_1(C_0(SX)\rcross G)=0\).  This yields the asserted long exact sequence.
  Equation~\eqref{eq:Eul_decompose} yields the formula for
  \(\Eul_X^\comb(x)\).
\end{proof}

We now describe the map \(\dim_H\colon \K_0(C_0(X)\rcross G)\to \K_0(\Cred
G)\) for a compact-open subgroup \(H\subseteq G\), which occurs in
Theorem~\ref{the:Gysin_combinatorial}.  It factors through the map
\(\K_0(C_0(X)\rcross G) \to \K_0(C_0(G/H)\rcross G)\cong \Rep(H)\) that is
induced by an orbit restriction map \(X\to G/H\).  The composite map
\[
\Rep(H) \cong \K_0(\Cred H) \cong \KK^G_0(C_0(G/H),\C) \to \Ktop_0(G) \to
\K_0(\Cred G)
\]
is equal to the induction map \(i_H^G\colon \Rep(H)\to\K_0(\Cred G)\), which
is induced by the embedding \(\Cred H\subseteq\Cred G\).  Thus \(\dim_H(x)\)
is the composite of an orbit restriction map and the induction map.  Since
there are relations between the orbit restriction and induction maps for
different~\(H\), it is hard to describe the range and kernel of
\(\Eul_X^\comb\) in general.

The following corollary is equivalent to Theorem \ref{the:bigtheorem} and
Corollary \ref{cor:torsion_unit}.  Let \(1_{\Cred G}\) be the unit projection
in~\(\Cred G\).

\begin{corollary}
  \label{cor:torsion_free}
  In the situation of Theorem \ref{the:Gysin_combinatorial}, suppose in
  addition that~\(G\) is discrete and torsion-free.  If \(G\backslash X\) is
  compact and \(\chi(G\backslash X)\neq0\), then there are exact sequences
  \[\xymatrix@C-1.1em{
    0 \ar[r] &
    \gen{\chi(G\backslash X)[1_{\Cred G}]} \ar[r]^-{\subseteq} &
    \K_0(\Cred G) \ar[r]^-{u_*} &
    \K_0(C(\bd{X}) \rcross G) \ar[r] &
    \K^1(G\backslash X) \ar[r] & 0,
  }
  \]
  \[\xymatrix@C-.3em{
    0 \ar[r] &
    \K_1(\Cred G) \ar[r]^-{u_*} &
    \K_1(C(\bd{X}) \rcross G) \ar[r] &
    \K^0(G\backslash X) \ar[r]^-{\dim} &
    \Z \ar[r] & 0,
  }
  \]
  and the class of the unit element in \(\K_0(C(\bd{X})\rcross G)\) has
  torsion of order \(\abs{\chi(G\backslash X)}\).  Otherwise, there are exact
  sequences
  \[\xymatrix{
    0 \ar[r] &
    \K_0(\Cred G) \ar[r]^-{u_*} &
    \K_0(C(\bd{X})\rcross G) \ar[r] &
    \K^1(G\backslash X) \ar[r] & 0,
  }
  \]
  \[\xymatrix{
    0 \ar[r] &
    \K_1(\Cred G) \ar[r]^-{u_*} &
    \K_1(C(\bd{X})\rcross G) \ar[r] &
    \K^0(G\backslash X) \ar[r] & 0,
  }
  \]
  and the class of the unit element in \(\K_0(C(\bd{X})\rcross G)\) has no
  torsion.
\end{corollary}

\begin{proof}
  Since the action on~\(X\) is free and proper, \(C_0(X)\rcross G\) is
  strongly Morita equivalent to \(C_0(G\backslash X)\).  Hence
  \(\K_*(C_0(X)\rcross G) \cong \K^*(G\backslash X)\).  Furthermore, we have
  an isomorphism \(\KK^G_0(C_0(X),\C)\cong\KK_0(C_0(G\backslash X),\C)\).  It
  maps \(\Eul_X^\comb\) to \(\Eul_{G\backslash X}^\comb = \chi(G\backslash X)
  \cdot\dim\).  The Kasparov product with the class
  \(\dim\in\KK^G_0(C_0(X),\C)\) factors through \(\KK^G_0(C_0(G),\C)\cong\Z\);
  one checks easily that it corresponds to the map
  \[
  \K^0(G\backslash X) \to \Z \to \K_0(\Cred G),
  \qquad x\mapsto \dim(x) \cdot [1_{\Cred G}];
  \]
  the reason for this is that \(\dim\in\KK_0^G(C_0(G),\C)\to \Ktop_0(G)\) is a
  pre-image for \([1_{\Cred G}]\) under the Baum-Connes assembly map.  Hence
  the kernel and range of the map \(\K^0(G\backslash X)\to\K_0(\Cred G)\) are
  the kernel of \(\chi(G\backslash X)\dim\) and \(\gen{\chi(G\backslash
    X)[1_{\Cred G}]}\), respectively.  Now the exact sequences follow from
  Theorem~\ref{the:Gysin_combinatorial}.  The assertions about the torsion of
  the unit element follow because \(u_*[1_{\Cred G}] = [1_{C(\bd{X})\rcross
    G}]\).
\end{proof}

\begin{example}
  \label{exa:free_group}
  Let~\(\F_n\) be the non-Abelian free group on \(n\)~generators for
  \(n\ge2\).  Let~\(X\) be the Cayley graph of~\(\F_n\), which is a regular
  tree, and let~\(\overline{X}\) be its ends compactification.  Let
  \(\bd{X}\defeq \overline{X}\setminus X\) be the set of ends of~\(X\), which
  is a Cantor set.  This compactification is the Gromov compactification of
  the hyperbolic group~\(\F_n\) and the visibility compactification of the
  \(\CAT(0)\) space~\(X\).  Of course, the group~\(\F_n\) is torsion-free, so
  that we are in the situation of Corollary \ref{cor:torsion_free}.  The
  group~\(\F_n\) satisfies the Baum-Connes conjecture with arbitrary
  coefficients by~\cite{Higson-Kasparov:Amenable}.

  The orbit space \(\F_n\backslash X\) is a wedge of \(n\)~circles, hence
  compact.  Therefore, \(\K_*(\Cred G)\cong \Ktop_*(G)\cong
  \K_*(\F_n\backslash X)\) and \(\Eul_X\in\Ktop_0(G)\) is the Euler
  characteristic of~\(G\).  The \(\K\)\nbd{}homology and \(\K\)\nbd{}theory of
  \(\F_n\backslash X\) are isomorphic to~\(\Z\) in degree~\(0\) and~\(\Z^n\)
  in degree~\(1\), and \(\chi(\F_n\backslash X)=1-n\).  Corollary
  \ref{cor:torsion_free} yields
  \[
  \K_0(C(\bd{X})\rcross \F_n) \cong \Z/\gen{n-1} \oplus \Z^n,
  \qquad
  \K_1(C(\bd{X})\rcross \F_n) \cong \Z^n.
  \]
  Therefore, the class of the unit element in \(\K_0(C(\bd{X})\rcross\F_n)\)
  is a torsion element of order \(n-1\).  This example is also studied
  in~\cite{Spielberg:Free-product}.

  If \(n=1\), we get \(\F_1=\Z\), \(X=\R\), and \(\bd{X}=\{\pm\infty\}\)
  with~\(\Z\) acting trivially.  In this case, the Euler characteristic
  vanishes; this already follows from Example~\ref{exa:fixed_boundary_point}.
\end{example}

\begin{example}
  \label{exa:surface_groups}
  Let~\(\Sigma_g\) be a closed surface of genus \(g\ge 2\) and let \(\Gamma_g
  \defeq\pi_1(\Sigma_g)\) be its fundamental group.  The universal cover
  of~\(\Sigma_g\) is homeomorphic to~\(\Hyp^2\), so that \(\Gamma_g\subseteq
  \Isom(\Hyp^2)\).  This implies that~\(\Gamma_g\) satisfies the Baum-Connes
  conjecture with coefficients (\cite{Kasparov:Warsaw}).  The usual
  compactification of~\(\Hyp^2\) by a circle at infinity \(\bd{\Hyp^2} \cong
  S^1\) is both the visibility compactification and the Gromov
  compactification of~\(\Hyp^2\) and therefore produces a boundary action.

  As in the previous example, \(\Gamma_g\) is torsion-free and
  \(\Gamma_g\backslash\Hyp^2 \cong \Sigma_g\) is compact, so that
  \(\Eul_{\Hyp^2}\in\Ktop_0(G)\) is the Euler characteristic of~\(G\).  The
  \(\K\)\nbd{}theory and \(\K\)\nbd{}homology of~\(\Sigma_g\) are isomorphic
  to~\(\Z^2\) in degree~\(0\) and~\(\Z^{2g}\) in degree~\(1\), and
  \(\chi(\Sigma_g)=2-2g\).  Therefore, Corollary \ref{cor:torsion_free} yields
  \[
  \K_0(C(\bd{\Hyp^2})\rcross\Gamma_g) \cong \Z/\gen{2g-2} \oplus \Z^{2g+1},
  \qquad
  \K_1(C(\bd{\Hyp^2})\rcross \Gamma_g) \cong \Z^{2g+1}.
  \]
  Explicit generators, as well as a dynamical proof of these assertions, can
  be found in~\cite{Emerson}.  This example is also studied in
  \cites{Connes:Cyclic_transverse, Connes:NCG, Delaroche, Natsume:Torsion}.
\end{example}

\begin{example}
  \label{exa:PSL2Z}
  Consider \(G\defeq \PSL(2,\Z)\), acting properly on the tree~\(X\) discussed
  in Example \ref{exa:Eul_comb_PSL2Z}.  Let~\(\cl{X}\) be the ends
  compactification of~\(X\) and let \(\bd{X}\defeq \cl{X}\setminus X\); this
  is the same as the Gromov or the visibility boundary of the tree~\(X\).  The
  group~\(G\) satisfies the Baum-Connes conjecture with arbitrary coefficients
  because it is a closed subgroup of \(\Isom(\Hyp^2)\).

  Since \(G\backslash X\) is compact, we have
  \(\KK^G_0(C_0(X),\C)\cong\Ktop_0(G)\), and the Euler characteristic
  \(\Eul_X\in\KK^G_0(C_0(X),\C)\cong\Ktop_0(G)\) is the Euler characteristic
  of~\(G\).  We have already computed \(\Eul_X^\comb=\Eul_X\) in Example
  \ref{exa:Eul_comb_PSL2Z}.  Hence
  \[
  \Eul_\EG = \dim_{\Z/2} + \dim_{\Z/3} - \dim_{\{1\}} \in \Ktop_0(G).
  \]

  Functions vanishing on the vertices form a \(G\)\nbd{}invariant ideal in
  \(C_0(X)\) that is isomorphic to \(C_0(\R\times G)\) with the free action
  of~\(G\); the quotient \(C^*\)\nbd{}algebra is isomorphic to
  \(C_0(G/\Ztwo)\oplus C_0(G/\Z/3)\).  The corresponding long exact sequences
  for \(\K_*(\blank\rcross G)\) and \(\KK^G_*(\blank,\C)\) are
  \[\xymatrix{
    0 \ar[r] &
    \K_0(C_0(X)\rcross G) \ar[r] &
    \Rep(\Z/2)\oplus\Rep(\Z/3) \ar[d]^{(\dim,\dim)} \\
    0 \ar[u] &
    \K_1(C_0(X)\rcross G) \ar[l] &
    \Z, \ar[l]
  }
  \]
  \[\xymatrix{
    \Rep(\Z/2)\oplus\Rep(\Z/3) \ar[r] &
    \KK_0^G(C_0(X),\C) \ar[r] &
    0 \ar[d] \\
    \Z \ar[u] &
    \KK_1^G(C_0(X),\C) \ar[l] &
    0. \ar[l]
  }
  \]
  The vertical map in the second exact sequence sends \(1\in\Z\) to
  \((\varrho,-\varrho)\), where~\(\varrho\) denotes the regular
  representation.  Thus \(\Ktop_1(G)\cong0\), \(\K_1(C_0(X)\rcross G)\cong0\),
  \(\Ktop_0(G)\cong\Z^4\), \(\K_0(C_0(X)\rcross G)\cong\Z^4\).  One can check
  that multiplication by~\(\Eul_\EG\) is a bijective map \(\K_0(C_0(X)\rcross
  G)\to \Ktop_0(G)\).  Hence we get \(\K_*(C(\bd{X})\rcross G)\cong0\).

  If we consider the boundary action on \(\bd{\Hyp^2}\), then we replace
  \(C_0(X)\) by \(C_0(\Hyp^2)\) in the Gysin sequence.  Equivariant Bott
  periodicity applies here and yields \(\K_*(C_0(\Hyp^2)\rcross G)\cong
  \Ktop_*(G)\).  Hence this group is concentrated in degree~\(0\).  One can
  check that the map \(\K_0(C_0(\Hyp^2)\rcross G)\to \K_0(C_0(X)\rcross G)\)
  that is induced by the embedding \(X\to\Hyp^2\) has kernel and cokernel
  isomorphic to~\(\Z\).  Since the map \(K_0(C_0(X)\rcross G)\to \Ktop_0(G)\)
  in the Gysin sequence for \(\bd{X}\) is invertible, the map
  \(K_0(C_0(\Hyp^2)\rcross G)\to \Ktop_0(G)\) in the Gysin sequence for
  \(\bd{\Hyp^2}\) has kernel and cokernel~\(\Z\) as well.  Thus
  \(\K_*(C(\bd{\Hyp^2})\rcross G)\cong\Z\) for \(*=0,1\).
\end{example}

\begin{example}
  \label{exa:reductive_on_Furstenberg}
  Let~\(G\) be a reductive \(p\)\nbd{}adic group and let \(\Gamma\subseteq G\)
  be a torsion-free discrete subgroup.  Let~\(X\) be the affine Bruhat-Tits
  building of~\(G\) and let \(\bd{X}_\vis\) be its visibility boundary.
  Recall that this is a boundary action of~\(G\).  The Fürstenberg boundary is
  \(G/P\), where \(P\subseteq G\) is a minimal parabolic subgroup.  Since
  there exist points in \(\bd{X}_\vis\) that are fixed by~\(P\), we get an
  embedding \(G/P\subseteq \bd{X}_\vis\), which induces a map
  \(\varphi\colon C(\bd{X}_\vis)\to C(G/P)\).

  We assume that \(\Gamma\subseteq G\) is cocompact or, equivalently, that
  \(\Gamma\backslash X\) is compact, and that \(\chi(\Gamma\backslash
  X)\neq0\).  We want to show that the class of the unit element in
  \(\K_0(C(G/P)\rcross\Gamma)\) is a torsion element whose order divides
  \(\chi(\Gamma\backslash X)\).  For certain buildings, this result has been
  obtained previously by Guyan Robertson in \cite{Robertson:Torsion}.  We
  remark that we get no information about the torsion of the unit element
  if~\(\Gamma\) fails to be cocompact or if \(\chi(\Gamma\backslash X)=0\).

  Observe first that \(\dim\in\Ktop_0(\Gamma)\) is a canonical choice of a
  pre-image for the class of the unit element in \(\Cred\Gamma\).  As in
  Corollary~\ref{cor:torsion_free}, we find that the image of \(\dim\) in
  \(\Ktop_0(\Gamma,C(\bd{X}_\vis))\) is a torsion element of order exactly
  equal to \(\abs{\chi(\Gamma\backslash X)}\).  Mapping further via
  \(\varphi\colon C(\bd{X}_\vis)\to C(G/P)\), we find that the image \(\dim'\)
  of \(\dim\) in \(\Ktop_0(\Gamma,C(G/P))\) is a torsion element whose order
  divides \(\chi(\Gamma\backslash X)\).  It is easy to see that the
  Baum-Connes assembly map sends \(\dim'\) to the class \([1]\) of the unit
  element in \(C(G/P)\rcross\Gamma\).  Hence \(\chi(\Gamma\backslash X)[1]=0\)
  as asserted.

  Similarly, let~\(G\) be an almost connected Lie group whose connected
  component is reductive and linear, and let \(\Gamma\subseteq G\) be a
  torsion-free, cocompact discrete subgroup.  Let \(X\defeq G/K\), where~\(K\)
  is the maximal compact subgroup of~\(G\), and let \(\bd{X}_\vis\) be its
  visibility boundary.  Again, this is a boundary action of~\(G\), and there
  exists an embedding \(G/P\subseteq \bd{X}_\vis\) of the Fürstenberg boundary
  because there exist points in \(\bd{X}_\vis\) that are fixed by~\(P\).
  Corollary \ref{cor:torsion_free} applies because~\(X\) has a
  \(G\)\nbd{}invariant triangulation.  Arguing as above, we see that the class
  of the unit element in \(\K_0(C(G/P)\rcross\Gamma)\) is a torsion element
  whose order divides \(\chi(\Gamma\backslash X)\) if the latter is non-zero.
\end{example}

\section{Equivariant Euler characteristics for smooth manifolds}
\label{sec:Euler_smooth}

If a locally compact group~\(G\) acts properly by diffeomorphisms on a smooth
manifold~\(M\), then there exists a complete \(G\)\nbd{}invariant Riemannian
metric on~\(M\).  Hence the action of~\(G\) factors through \(\Isom(M)\),
which is a Lie group unless~\(M\) has infinitely many connected components.
Throughout this section, we consider the situation where~\(M\) is a complete
Riemannian manifold and~\(G\) is a locally compact group that acts
isometrically on~\(M\).  It does not matter whether or not this action is
proper because \(\Isom(M)\) acts properly on~\(M\) in any case.

We first recall the construction of a Kasparov dual for~\(M\) in
\cite{Kasparov:Novikov}*{Section 4}.  Then we identify \(\Eul_M\) with the
\(\K\)\nbd{}homology class of the de Rham differential operator.  If~\(G\) is
discrete and acts properly, then we can also compute the Euler characteristic
combinatorially, using a \(G\)\nbd{}invariant triangulation of~\(M\).  Since
the Euler characteristic is independent of the choice of Kasparov dual, the
combinatorial and the de-Rham-Euler characteristics agree; this result is due
to Jonathan Rosenberg and Wolfgang Lück (\cite{Luck-Rosenberg:Euler}).

Let \(\Cliff\) be the bundle whose fibre at \(x\in M\) is the Clifford algebra
for the vector space~\(T^*_xM\) equipped with the inner product defined by the
Riemannian metric.  This is a bundle of \(\Ztwo\)\brd{}graded
finite-dimensional \(C^*\)\nbd{}algebras, on which~\(G\) acts in a canonical
way.  Let \(\sour \defeq C_0(M,\Cliff)\) be the \(\Ztwo\)\brd{}graded
\(C^*\)-algebra of \(C_0\)\nbd{}sections of this bundle, equipped with its
canonical action of~\(G\).  (It is denoted \(C_\tau(M)\)
in~\cite{Kasparov:Novikov}*{4.1}.)  We have \(C_0(M)\subseteq\sour\) as scalar
valued functions.  This embedding is central, so that~\(\sour\) becomes a
\(\Ztwo\)-graded \(M\cross G\)-\(C^*\)-algebra.

Now we describe \(D\in\KK^G_0(\sour,\C)\) (see \cite{Kasparov:Novikov}*{4.2}).
Let \(\Lambda^*M = \bigoplus_n \Lambda^n M\) be the bundle of differential
forms on~\(M\), graded by parity.  Let \(C_c^\infty(\Lambda^*M)\) be the space
of smooth, compactly supported sections of \(\Lambda^*M\), and let
\(L^2(\Lambda^*M)\) be the Hilbert space completion of
\(C_c^\infty(\Lambda^*M)\) with respect to the standard inner product given by
the Riemannian metric.  We let~\(\sour\) act on \(L^2(\Lambda^*M)\) by
Clifford multiplication.  Let \(d \colon C_c^\infty(\Lambda^*M)\to
C_c^\infty(\Lambda^*M)\) be the de Rham differential.  The operator \(d+d^*\)
is an essentially self-adjoint, \(G\)\nbd{}invariant unbounded operator on
\(L^2(\Lambda^*M)\).  Together with the representation of~\(\sour\) it defines
a class \(D\in\KK^G_0(\sour,\C)\).  Here we use the unbounded picture of
Kasparov theory by Saad Baaj and Pierre Julg (\cite{Baaj-Julg:KK_unbounded}).

Next we define \(\Theta\in\RKK^G_0(M;\C,\sour)\) as in
\cite{Kasparov:Novikov}*{4.3--4.4}.  The basic ingredients are the geodesic
distance function \(\varrho\colon M\times M\to M\) and a \(G\)\nbd{}invariant
function \(r\colon M\times M\to \R_{>0}\) such that any \(x,y\in M\) with
\(\varrho(x,y)<r(x)\) are joined by a unique geodesic.  Let
\[
U\defeq \{(x,y)\in M\times M\mid \varrho(x,y)<r(x)\}
\]
and pull back \(T^*M\) to a bundle \(\pi_2^*T^*M\) on~\(U\) via the coordinate
projection \(\pi_2\colon (x,y)\mapsto y\).  Let \(J_U\subseteq C_0(M)\hot
\sour\) be the ideal of sections that vanish outside~\(U\).  We view~\(J_U\)
as a \(G\)\nbd{}equivariant \(\Ztwo\)\nbd{}graded Hilbert module over
\(C_0(M)\hot \sour\).

Define a covector field~\(F\) on~\(U\) by
\[
F(x,y)\defeq \frac{\varrho(x,y)}{r(x)} \cdot d_2 \varrho(x,y)
\in C_0(U,\pi_2^* T^*M).
\]
where~\(d_2\) is the exterior derivative in the second variable~\(y\).  The
covector field~\(F\) defines a \(G\)\nbd{}invariant self-adjoint, odd
multiplier of~\(J_U\); it satisfies \((f\hot 1)\cdot (1-F^2)\in J_U\) for all
\(f\in C_0(M)\).  Thus \(\Theta = (J_U,F)\) is a cycle for
\(\RKK^G_0(M;\C,\sour)\).  It is asserted in \cite{Kasparov:Novikov}*{4.5,
  4.6, 4.8} that \((\sour, D, \Theta)\) is a Kasparov dual for~\(M\) in the
sense of Definition \ref{def:Kasparov_dual}.

\begin{definition}
  \label{def:deRham_Euler}
  Let \(\Eul_M^\dR\in \KK^G_0(C_0(M),\C)\) be the class determined by the
  representation of \(C_0(M)\) on \(L^2(\Lambda^*M)\) and the operator
  \(d+d^*\) described above.  We call \(\Eul_M^\dR\) the
  \emph{\(G\)\nbd{}equivariant de-Rham-Euler characteristic of~\(M\)}.
\end{definition}

Equivalently, we get \(\Eul_M^\dR\) from \(D\in\KK^G_0(\sour,\C)\) by
restricting the representation of~\(\sour\) to \(C_0(M)\subseteq\sour\).

\begin{theorem}
  \label{the:Euler_smooth}
  Let~\(M\) be a complete Riemannian manifold and let~\(G\) be a locally
  compact group acting isometrically on~\(M\).  Then the abstract
  \(G\)\nbd{}equivariant Euler characteristic \(\Eul_M\) is equal to
  \(\Eul_M^\dR\).
\end{theorem}

\begin{proof}
  Lemma~\ref{lem:abstract_Euler_Kd} asserts that
  \[
  \Eul_M = \bar\Theta \hot_{C_0(M)\hot\sour} [m] \hot_{\sour} D
  \in\KK^G_0(C_0(M),\C),
  \]
  where \(\bar\Theta\in\KK^G_0(C_0(M),C_0(M)\hot\sour)\) is obtained
  from~\(\Theta\) by forgetting the \(M\)\nbd{}structure and where \(m\colon
  C_0(M)\hot\sour \to \sour\) is the \(M\)\nbd{}structure homomorphism
  of~\(\sour\).  We first compute
  \[
  \bar\Theta \hot_{C_0(M)\hot\sour} [m]
  = m_*(\bar\Theta)\in \KK^G_0(C_0(M),\sour).
  \]
  Its underlying Hilbert module is \(J_U\hot_{C_0(M)\hot\sour} \sour\); this
  is isomorphic to~\(\sour\) because \(J_U\) is an ideal in
  \(C_0(M)\hot\sour\) that contains all functions supported in some
  neighbourhood of the diagonal, and the multiplication homomorphism~\(m\)
  restricts to the diagonal.  The action of \(C_0(M)\) is by multiplication
  on~\(J_U\); this corresponds to the embedding \(m'\colon C_0(M)\to\sour\) by
  scalar valued functions.  Since this homomorphism maps into~\(\sour\) and
  not just into \(\Mult(\sour)\), the Fredholm operator is irrelevant.  Thus
  \(m^*(\bar\Theta)=m'\).

  Taking the Kasparov product with \(D\in\KK^G_0(\sour,\C)\), we get
  \[
  \Eul_M
  = \bar\Theta \hot_{C_0(M)\hot\sour} [m] \hot_{\sour} D
  = [m'] \hot_{\sour} D
  = (m')^*(D).
  \]
  This is equal to \(\Eul_M^\dR \in \KK^G_0(C_0(M),\C)\) by definition.
\end{proof}

\begin{theorem}[Lück and Rosenberg \cite{Luck-Rosenberg:Euler}]
  \label{the:luckrosenberg}
  Let~\(M\) be a smooth manifold and let~\(G\) be a discrete group acting
  on~\(M\) properly by diffeomorphisms.  Then the de-Rham-Euler characteristic
  and the combinatorial Euler characteristic agree:
  \[
  \Eul_M^\dR = \Eul_M = \Eul_M^\comb \in \KK^G_0(C_0(M),\C).
  \]
\end{theorem}

\begin{proof}
  We have seen in Section~\ref{sec:Kasparov_duality} that the (abstract)
  equivariant Euler characteristic is independent of the Kasparov dual.  Hence
  the assertion follows by combining Theorems \ref{the:Euler_combinatorial}
  and~\ref{the:Euler_smooth}.  Recall that a smooth manifold equipped with a
  proper action of a discrete group~\(G\) always admits a Riemannian metric
  for which the group acts isometrically and a triangulation for which the
  group acts simplicially (\cite{Illman:Equivariant_Triangulations}).
\end{proof}

\begin{remark}
  The analogues of Theorems \ref{the:Euler_smooth} and \ref{the:luckrosenberg}
  in real \(\K\)\nbd{}homology also hold, by exactly the same arguments.
\end{remark}

As before, let~\(M\) be an \(n\)\nbd{}dimensional Riemannian manifold equipped
with an isometric action of~\(G\).  Assume, in addition, that~\(M\) is
\(G\)\nbd{}equivariantly \(K\)\nbd{}oriented.  This means that the tangent
bundle of~\(M\) has a \(G\)\nbd{}equivariant spinor bundle \(\Spinor\) (see
\cite{Connes-Skandalis}).  This bundle gives rise to a \(G\)\nbd{}equivariant
Morita equivalence between \(\sour=C_0(M,\Cliff)\) and the trivial Clifford
algebra bundle \(C_0(M)\hot \Cliff(\R^n)\) (see \cite{Plymen:Spin_Morita}).
Therefore, \(C_0(M)\) and~\(\sour\) are \(\KK^{M\cross G}\)\nbd{}equivalent
with a dimension shift of~\(n\).  We may transport the structure of a Kasparov
dual from~\(\sour\) to \(C_0(M)\).  It is easy to see that
\(D\in\KK^G_0(\sour,\C)\) corresponds to the Dirac operator
\(\Dslash_M\in\KK^G_n(C_0(M),\C)\), which acts on sections of \(\Spinor\).
The map \(\sigma_{M,C_0(M)}\) is simply the forgetful map as
in~\eqref{eq:forget_Y}.  Hence the inverse of the Poincaré duality map is
given by
\[
\PD^{-1}\colon \RKK^G_*(M;\C,\C)\to\KK^G_{*+n}(C_0(M),\C),
\qquad f\mapsto (-1)^n\bar{f}\hot_{C_0(M)} \Dslash_M.
\]
We also get a class \([\Spinor]\in\RKK^G_{-n}(M;\C,\C)\) by taking
\(C_0\)\nbd{}sections of \(\Spinor\) with \(F=0\).  Using
\(\Spinor\hot\Spinor\cong \Lambda^*M\), one shows easily that
\begin{equation}
  \label{eq:PD_Spinor}
  \PD^{-1}[\Spinor] = (-1)^n \overline{[\Spinor]}\hot_{C_0(M)} \Dslash_M
  = (-1)^n \Eul_M^\dR = (-1)^n \Eul_M.
\end{equation}
Together with Theorem \ref{the:luckrosenberg}, this allows us to compute the
classical, commutative Gysin sequence in \(\K\)\nbd{}theory.  We have seen in
the introduction how to get a long exact Gysin sequence of the form
\[\xymatrix@-2mm{
  \dotsb \ar[r] &
  \K^{*-n}(M) \ar[r]^-{\epsilon^*} &
  \K^*(M) \ar[r]^-{\pi^*} &
  \K^*(SM) \ar[r]^-{\delta} &
  \K^{*-n+1}(M) \ar[r] &
  \dotsb,
}
\]
where \(\epsilon^*(x)= x\hot\Spinor\) and \(\pi\colon SM\to M\) is the
bundle projection.  This is proved in \cite{Karoubi:K-theory}*{IV.1.13} in a
more general context, where the tangent bundle is replaced by any
\(\K\)\nbd{}oriented vector bundle on~\(M\).  If we specialise to the tangent
bundle, then \eqref{eq:PD_Spinor} yields:

\begin{theorem}
  \label{the:Gysin_classical}
  Let~\(M\) be a \(\K\)\nbd{}oriented, connected \(n\)\nbd{}dimensional
  manifold.  Let \(SM\) be its sphere bundle and let \(\pi\colon SM\to M\) be
  the bundle projection.  If~\(M\) is compact and \(\chi(M)\neq0\), then there
  are exact sequences
  \[\xymatrix@C-1.1em{
    0 \ar[r] &
    \gen{\chi(M) \pnt!} \ar[r]^-{\subseteq} &
    \K^n(M) \ar[r]^-{\pi^*} &
    \K^n(SM) \ar[r] &
    \K^1(M) \ar[r] & 0,
  }
  \]
  \[\xymatrix@C-.3em{
    0 \ar[r] &
    \K^{n+1}(M) \ar[r]^-{\pi^*} &
    \K^{n+1}(SM) \ar[r] &
    \K^0(M) \ar[r]^-{\dim} &
    \Z \ar[r] & 0,
  }
  \]
  where \(\pnt!\in\KK_{-n}(C(*),C(M))\cong \K^n(M)\) is the wrong way element
  associated to the inclusion of a point in~\(M\).

  If~\(M\) is not compact or if \(\chi(M)=0\), then there are exact sequences
  \[\xymatrix{
    0 \ar[r] &
    \K^n(M) \ar[r]^-{\pi^*} &
    \K^n(SM) \ar[r] &
    \K^1(M) \ar[r] & 0,
  }
  \]
  \[\xymatrix{
    0 \ar[r] &
    \K^{n+1}(M) \ar[r]^-{\pi^*} &
    \K^{n+1}(SM) \ar[r] &
    \K^0(M) \ar[r] & 0.
  }
  \]
\end{theorem}

\begin{proof}
  By Theorem \ref{the:luckrosenberg} we have \(\Eul_M^\dR=\Eul_M^\comb\).
  Since there is no group action, we have
  \(\Eul_M^\comb=\chi(M)\cdot\dim\in\KK_0(C(M),\C)\) if~\(M\) is compact, and
  \(\Eul_M^\comb=0\) otherwise.  Hence~\eqref{eq:PD_Spinor} yields
  \([\Spinor]=0\) if~\(M\) is not compact or if \(\chi(M)=0\) and finishes the
  proof in that case.  Otherwise, \([\Spinor]=(-1)^n\chi(M)\cdot \PD(\dim)\).
  It is easy to see that \(\PD(\dim)=\pnt!\); recall that \(\pnt!\) is the
  class of the map \(C_0(\R^n)\to C_0(M)\) given by a diffeomorphism
  from~\(\R^n\) onto some (small) open ball in~\(M\).  Since any bundle
  restricts to a trivial bundle on this open ball, we get \(x\hot\pnt!=0\) for
  \(x\in\K^1(M)\) and \(x\hot \pnt!=\dim(x)\cdot\pnt!\) for \(x\in\K^0(M)\).
  Hence the kernel and range of the map \(\K^0(M)\to\K^n(M)\) in the Gysin
  sequence are equal to the kernel of \(\dim\) and \(\gen{\chi(M) \pnt!}\),
  respectively.  This yields the desired exact sequences.
\end{proof}

Now we return to the situation of boundary actions and formulate the
non-commutative Gysin sequence.

\begin{theorem}
  \label{the:Gysin_smooth}
  Let \(\bd{X}=\cl{X}\setminus X\) be a boundary action of a locally compact
  group~\(G\) as in Definition \ref{def:boundary_action}.  Suppose that~\(X\)
  is a complete Riemannian manifold on which~\(G\) acts isometrically.
  Suppose also that~\(G\) satisfies the Baum-Connes conjecture with
  coefficients in \(\C\) and \(C(\bd{X})\).  Then there is an exact sequence
  \[
  \xymatrix{
    \K_0(C_0(X)\rcross G) \ar[r]^-{\Eul_X^\dR} &
    \K_0(\Cred G) \ar[r]^-{u_*} &
    \K_0(C(\bd{X})\rcross G) \ar[d]^{\delta} \\
    \K_1(C(\bd{X})\rcross G) \ar[u]_{\delta} &
    \K_1(\Cred G) \ar[l]_-{u_*} &
    \K_1(C_0(X)\rcross G), \ar[l]_-{\Eul_X^\dR}
  }
  \]
  where~\(\Eul_X^\dR\) denotes the Kasparov product with the equivariant
  de-Rham-Euler characteristic \(\Eul_X^\dR\in\KK^G_0(C_0(X),\C)\).
\end{theorem}

\begin{proof}
  The assumption about the Baum-Connes conjecture allows us to replace
  \(\K_*(\Cred G)\) and \(\K_*(C(\bd{X})\rcross G)\) by \(\Ktop_*(G)\) and
  \(\Ktop_*(G,C(\bd{X}))\).  Hence the result follows from Theorem
  \ref{the:Euler_smooth} and Proposition \ref{pro:abstract_Gysin}.
\end{proof}

\begin{proposition}
  \label{pro:Euler_vanishes_odd}
  Let~\(M\) be an oriented Riemannian manifold of odd dimension and suppose
  that~\(G\) acts on~\(M\) by orientation-preserving isometries.  Then
  \(\Eul_M^\dR=0\) in \(\KK^G_0(C_0(M),\C)\).
\end{proposition}

\begin{proof}
  Let \(\mathrm{Vol}\) be the canonical volume form.  Let \(n\defeq \dim M\)
  and write \(n=2k+1\) with \(k\in\N\).  We shall use the \emph{Hodge
    \(\star\) operation} (\cite{deRham:Varietes}*{\S 24--25}).  It is a
  \(C_0(M)\)\brd{}linear map \({\star}\colon L^2(\Lambda^p M)\to
  L^2(\Lambda^{n-p} M)\) for all \(p\in\no\); it is defined by \(\beta\wedge
  \star\alpha = (\alpha,\beta)\cdot \mathrm{Vol}\) for all \(\alpha,\beta\in
  L^2(\Lambda^p M)\), where \((\alpha,\beta)\in C_0(M)\) denotes the pointwise
  inner product induced by the Riemannian metric.  The operator~\(\star\) is
  unitary on \(L^2(\Lambda^*M)\) and satisfies
  \[
  \star\star\alpha = (-1)^{pn+p}\alpha,
  \qquad
  d^*(\alpha) = (-1)^{pn+n+1}\star d\star \alpha
  \]
  for all \(\alpha\in L^2(\Lambda^p M)\) (see \cite{deRham:Varietes}*{\S
    24--25}).  Consider the operator
  \[
  \epsilon\colon L^2(\Lambda^* M)\to L^2(\Lambda^* M),
  \qquad \alpha \mapsto \mathrm{i}^{k+p(p-1)} \star\alpha
  \quad \text{for \(\alpha\in L^2(\Lambda^p M)\).}
  \]
  Straightforward computations show that \(\epsilon^2=1\) and \(\epsilon
  d\epsilon = -d^*\); this implies that~\(\epsilon\) anti-commutes with
  \(d+d^*\).

  Since~\(\epsilon\) is still unitary and odd, it generates a
  grading-preserving representation of the Clifford algebra \(\Cliff(\R)\) on
  \(L^2(\Lambda^* M)\).  It commutes with the representations of~\(G\) and
  \(C_0(M)\) because~\(G\) acts by orientation-preserving maps and~\(\star\)
  is \(C_0(M)\)\nbd{}linear.  Thus \((L^2(\Lambda^* M),d+d^*)\) becomes a
  cycle~\(D_1\) for \(\KK^G_0(C_0(M)\hot\Cliff(\R),\C)\).  We have
  \(\Eul_X^\dR = [u_{\Cliff(\R)}] \hot_{\Cliff(\R)} D_1\), where
  \(u_{\Cliff(\R)}\colon \C\to\Cliff(\R)\) is the unit map.  The Kasparov
  cycle \([u_{\Cliff(\R)}]\) is evidently degenerate.  Hence \(\Eul_X^\dR=0\).
\end{proof}

In the real case, the same argument still goes through if \(\dim M\equiv 1
\pmod 4\) because then \(\mathrm{i}^{k+p(p-1)}=\pm1\) for all~\(p\).

The assumption that the action of~\(G\) be orientation-preserving is necessary
in Proposition \ref{pro:Euler_vanishes_odd}.  For a counterexample, take
\(M=S^1\) and \(G=\Ztwo\) acting on~\(S^1\) by reflection in the
\(x\)\nbd{}axis (equivalently, by complex conjugation).  A straightforward
computation shows \(\Eul_{S^1}^\comb\neq0\) in \(\KK^{\Ztwo}_0(C(S^1),\C)\).

\begin{example}
  \label{exa:semisimple_Lie}
  Let~\(G\) be a connected semi-simple Lie group and let \({K\subseteq G}\) be
  its maximal compact subgroup.  Then the homogeneous space \(X\defeq G/K\) is
  a \(\CAT(0)\) space, and we can take its visibility
  compactification~\(\cl{X}_\vis\) and its visibility boundary~\(\bd{X}_\vis\)
  as in Section~\ref{sec:setup}.  Both \(X\) and~\(\cl{X}_\vis\) are strongly
  contractible.  If~\(G\) is of rank~\(1\), then \(\cl{X}_\vis\) is the
  Fürstenberg boundary of~\(G\); in general, one can find a stratification
  of~\(\cl{X}_\vis\) whose subquotients are suspensions of \(G/P\) for
  parabolic subgroups \(P\subseteq G\); all parabolic subgroups of~\(G\) occur
  as a stabiliser for some point in the visibility boundary.

  It is known that almost connected groups satisfy the Baum-Connes conjecture
  with trivial coefficients (\cite{Chabert-Echterhoff-Nest:Connes-Kasparov}),
  which is also called Connes-Kasparov conjecture in this special case.
  Therefore, \(G\) satisfies the Baum-Connes conjecture with coefficients
  \(\C\) and \(C(G/P)\) for all parabolic subgroups~\(P\).  Using excision and
  the stratification, this implies that~\(G\) satisfies the Baum-Connes
  conjecture with coefficients \(C(\bd{X}_\vis)\).  Therefore,
  Theorem~\ref{the:Gysin_smooth} applies to our situation.

  The crossed product \(C_0(X)\rcross G\) is Morita-Rieffel equivalent to
  \(C^*(K)\), so that \(\K_0(C_0(X)\rcross G)\cong\Rep(K)\) and
  \(\K_1(C_0(X)\rcross G)=0\).  Since~\(X\) is a \(G\)\nbd{}compact universal
  proper \(G\)\nbd{}space, we have
  \[
  \Ktop_*(G) \cong \KK^G_*(C_0(X),\C) \cong \KK^K_*(C_0(X),\C),
  \]
  and \(\Eul_X^\dR\) is the Euler characteristic \(\Eul_\EG\) of~\(G\) (see
  Definition \ref{def:abstract_Euler_group}).

  Since we assume~\(G\) to be connected, \(K\) acts on~\(X\) by
  orientation-preserving maps.  Hence Proposition \ref{pro:Euler_vanishes_odd}
  yields \(\Eul_X^\dR=0\) if \(\dim X\) is odd.  In this case, our Gysin
  sequence splits into two short exact sequences.

  We suppose from now on that \(\dim X\) is even.  Recall that, as a
  \(K\)\nbd{}space, \(X\) is diffeomorphic to a certain linear
  representation~\(\pi_X\) of~\(K\).  Assume that this representation lifts to
  \(\mathrm{Spin^c}\); thus we can form a spinor representation~\(S_X\), which
  we view as an element of \(\Rep(K)\).  We have natural isomorphisms
  \[
  \KK^G_*(C_0(X),\C) \cong \RKK^G_*(X;\C,\C) \cong \KK^K_*(\C,\C) \cong
  \begin{cases}
    \Rep(K) & \text{\(*\) even};\\
    0 & \text{\(*\) odd}.
  \end{cases}
  \]
  The first of these maps is the Poincaré duality isomorphism for the Kasparov
  dual \(C_0(X)\) and maps \(\Eul_X^\dR\) to \([\Spinor]\)
  by~\eqref{eq:PD_Spinor}; the second isomorphism is induced by restriction to
  the point \(K\in X\) and hence maps \([\Spinor]\) to \(S_X\).  Thus the map
  denoted \(\Eul_X^\dR\) in the Gysin sequence is equivalent to the map
  \[
  \Rep(K)\to\Rep(K),\qquad \pi\mapsto \pi\hot[S_X].
  \]

  It is well-known that the ring \(\Rep(K)\) has no zero-divisors.  Hence
  either \([S_X]=0\) or this map is injective.  In the first case, our Gysin
  sequence splits again into two short exact sequences; in the second case, it
  follows that
  \[
  \K_0(C(\bd{X}_\vis)\rcross G) \cong \Rep(K)/(S_X),
  \qquad
  \K_1(C(\bd{X}_\vis)\rcross G) \cong0,
  \]
  where \((S_X)\) denotes the ideal generated by the virtual
  representation~\(S_X\).

  Let~\(\chi_G\) be the character of the representation~\(S_X\).  Then
  \(\abs{\chi_G}^2\) is the character of the representation \(S_X\hot S_X^*\cong
  \Lambda^*X\).  If \(g\in G\), let \(\lambda_1,\dotsc,\lambda_n\) be the
  eigenvalues of~\(g\) acting on~\(X\), counted with multiplicity; then we can
  describe the eigenvalues of~\(g\) on \(\Lambda^*X\) as well and get
  \[
  \abs{\chi_G(g)}^2 = \prod_{j=1}^n (1-\lambda_j)
  = \det (1- \pi_X(g)).
  \]
  Hence \(\chi_G(g)=0\) if and only if~\(1\) is an eigenvalue of
  \(\pi_X(g)\colon X\to X\).
\end{example}

\end{document}